\documentclass[11pt,leqno]{amsart}
\usepackage{amsmath,amsfonts,amsthm,amscd,amssymb,graphicx}
\numberwithin{equation}{section}
\usepackage{fullpage}
\usepackage{epstopdf}

\usepackage{color}

\newcommand{\mB}{{\mathbb B}}

\newcommand{\bpm}{\begin{pmatrix}}
\newcommand{\epm}{\end{pmatrix}}

\def\eps{\varepsilon }

\newcommand{\dif}{\ensuremath{\mathrm{d}}}

\newcommand\R{\mathbb R}

\def\eps{\varepsilon}

\newcommand\br{\begin{remark}}
\newcommand\er{\end{remark}}
\newcommand\bp{\begin{pmatrix}}
\newcommand\ep{\end{pmatrix}}
\newcommand{\be}{\begin{equation}}
\newcommand{\ee}{\end{equation}}
\newcommand\ba{\begin{equation}\begin{aligned}}
\newcommand\ea{\end{aligned}\end{equation}}

\newcommand{\bap}{\begin{app}}
\newcommand{\eap}{\end{app}}
\newcommand{\begs}{\begin{exams}}
\newcommand{\eegs}{\end{exams}}
\newcommand{\beg}{\begin{example}}
\newcommand{\eeg}{\end{exaplem}}
\newcommand{\bpr}{\begin{proposition}}
\newcommand{\epr}{\end{proposition}}
\newcommand{\bt}{\begin{theorem}}
\newcommand{\et}{\end{theorem}}
\newcommand{\bc}{\begin{corollary}}
\newcommand{\ec}{\end{corollary}}
\newcommand{\bl}{\begin{lemma}}
\newcommand{\el}{\end{lemma}}
\newcommand{\bd}{\begin{definition}}
\newcommand{\ed}{\end{definition}}
\newcommand{\brs}{\begin{remarks}}
\newcommand{\ers}{\end{remarks}}

\newcommand{\mat}[1]{\begin{pmatrix}#1\end{pmatrix}}

\newcommand{\mA}{{\mathbb A}}

\newcommand{\RR}{{\mathbb R}}

\newcommand{\CC}{{\mathbb C}}

\newcommand{\const}{\text{\rm constant}}
\newcommand{\Id}{{\rm Id }}

\newcommand{\sgn}{\text{\rm sgn}}
\newtheorem{theorem}{Theorem}[section]
\newtheorem{proposition}[theorem]{Proposition}
\newtheorem{corollary}[theorem]{Corollary}
\newtheorem{lemma}[theorem]{Lemma}

\theoremstyle{remark}
\newtheorem{remark}[theorem]{Remark}
\theoremstyle{definition}
\newtheorem{definition}[theorem]{Definition}

\newtheorem{example}[theorem]{Example}

\newcommand{\RM}{\mathbb{R}}

\newcommand{\beq}{\begin{equation}}
\newcommand{\eeq}{\end{equation}}

\newcommand{\besn}[1]{\begin{equation}\begin{split}#1 \end{split} \notag \end{equation}}

\title{
Convex entropy, Hopf bifurcation, and \\
viscous and inviscid shock stability 
}

\author{ Blake Barker}
\address{ Indiana University, Bloomington, IN 47405}
\email{bhbarker@gmail.com}
\thanks{Research of B.B. was partially supported
under NSF grant no. DMS-0300487}
\author{Heinrich Freist\"uhler}
\thanks{Research of H.F. was partially supported
under DFG Excellence Grant 2007-2012 to the University of Konstanz.}
\address{
Universit\"at Konstanz, 78457 Konstanz, Germany} 
\email{heinrich.freistuehler@uni-konstanz.de}
\author{Kevin Zumbrun}
\address{Indiana University, Bloomington, IN 47405}
\email{kzumbrun@indiana.edu} 
\thanks{Research of B.B. and K.Z. was partially supported
under NSF grant no. DMS-0300487}
\begin{document}

\date{\today}
\maketitle
\begin{abstract}
We consider by a combination of analytical and numerical techniques
some basic questions regarding 
the relations between inviscid and viscous stability
and existence of a convex entropy.
Specifically, for a system possessing a convex 
entropy, in particular for the equations of gas dynamics with a convex
equation of state, 
we ask:
(i) can inviscid instability occur? 
(ii) can there occur viscous instability not detected by inviscid theory?
(iii) can there occur the ---necessarily viscous--- effect of
Hopf bifurcation, or ``galloping instability''? 
and, perhaps most important from a practical point of view,
(iv) as shock amplitude is increased from the (stable) weak-amplitude
limit,
can there occur a first transition from viscous stability to instability
that is not detected by inviscid theory?
We show that (i) does occur for strictly hyperbolic, genuinely nonlinear
gas dynamics with certain convex equations of state.
while (ii) and (iii) do occur for an artifically constructed system with convex 
viscosity-compatible entropy. 
We do not know of an example 
for which (iv) occurs, leaving this as a key open question in viscous
shock theory, related to the principal eigenvalue property of 
Sturm Liouville and related operators.
In analogy with, and partly proceeding close to, the analysis of Smith on (non-)uniqueness of 
the Riemann problem,
we obtain convenient criteria for shock (in)stability
in the form of necessary and sufficient conditions  on 
the equation of state.
\end{abstract}

\tableofcontents

\section{Introduction }\label{s:intro}
Stability of shock waves has been the subject of intensive investigation
over the more than 50 years since the question was opened in the inviscid
setting by Landau, 
Kontorovich, Dy'akov, Lax, Erpenbeck, and others in the early 1960's.
Before that period, shocks seem to have been assumed to be universally
stable, and in fact under most reasonable conditions, they are so.
See \cite{BE} for an interesting general discusion of these early 
investigations, and \cite{S,MP} 
on the related question of uniqueness of Riemann solutions;
we recommend also the original sources \cite{Er1,Er2,G,M1,M2,M3}.
Starting in the mid-1980's, these investigations have been widened to
include also viscous shock stability, or stability in the presence of regularizing transport effects such as viscosity, heat conduction, magnetic resistivity,
and species diffusion.
See the surveys \cite{Z1,Z2,Z3,Z5} for general accounts of progress in this
direction, and \cite{BiB,GMWZ1,GMWZ2} for accounts of progress on the related
inviscid limit problem.
See also the recent investigations
\cite{TZ1,TZ2,TZ3,SS,BeSZ} on the possibility of 
Hopf bifurcation of viscous shocks, 
and stability of resulting time-periodic ``galloping'' shocks.

At a technical level, stability and bifurcation of shock waves is now
fairly well understood, both in one and multi-dimensions.
In particular, nonlinear inviscid shock stability reduces \cite{M1}
to a spectral condition that may be checked by evaluation of 
an associated Evans--Lopatinski function commputable by standard linear
algebraic operations; see \cite{Er1,M1,T,FP} for studies of the
Lopatinski function in various contexts. 
Likewise, determination of nonlinear viscous shock stability/bifurcation 
has been reduced to verification of spectral conditions on the linearized
operator about the wave, that may be readily checked by 
efficient and well-conditioned numerical Evans function computations 
as described in \cite{Br,BrZ,HuZ1,Z4,Z5}
(see, for example, the recent studies  
\cite{HLZ,BHRZ,BHZ1,BLeZ,BLZ,HLyZ1,HLyZ2,Z6}).
(Cf. \cite{HuZ3,PZ,FS1,FS2}, \cite[Thm. 1.4]{HLZ}, \cite[section 4]{BHZ1} for situations 
in which the spectral conditions may be checked analytically.)


However, at a level of basic intuition/understanding, some important questions
remain, even in the one-dimensional setting.
In particular, for all of the above-mentioned investigations of physical
systems, whether inviscid or viscous, shocks were seen to be 
{\it one-dimensionally stable}.
Indeed, as regards physically relevant systems of continuum mechanics,  
one-dimensionally unstable waves have so far been found
only for media permitting phase transitions 
\cite{GZ,Z6,ZMRS,BM,TZ4},
and established wisdom \cite{Er1,BE,MP} seems to say that 
instability is associated with 
effects not modeled in an ideal gas equation of state.

One could ask, therefore, whether there might be some simple and commonly
satisfied structural condition guaranteeing one-dimensional stability.
In particular, in the simplest and most familiar setting of 
inviscid gas dynamics,
{\it could the classical condition of thermodynamic stability, or
convexity of the equation of state 
$$
e=\bar e(\tau, S)
$$
relating internal energy $e$ to specific volume $\tau$ and entropy $S$ by itself
be sufficient to imply one-dimensional stability of shock waves?} 
Remarkably, the answer to this natural question 
up to now does not appear to have been known.\footnote{The suggested counterexample of \cite{G} is incorrect,
as we show below; see Remark \ref{gardnerrmk}.
}
%

A natural generalization pointed out by Lax \cite{L}
of the property of thermodynamic stability to
arbitrary systems of inviscid conservation laws
\be\label{cons}
u_t+f(u)_x=0
\ee
is existence of a convex entropy, i.e., a function $\eta(u)$
satisfying 
$$
d^2\eta>0
\quad\text{and}\quad
d\eta\ \! df=dq
$$ 
for an appropriate entropy flux function
$q(u)$, whence, for smooth solutions,
$$
{\eta(u)_t + q(u)_x=0}.
$$
One could wonder 
{\it whether for general systems \eqref{cons}, 
existence of a convex entropy is sufficient
to imply one-dimensional stability of inviscid shock waves}.
This question includes the above one as a special case, as 
in the case of gas dynamics, one can take $\eta$ as
the negative of the thermodynamical entropy $S(\tau,e)$,
the Legendre transform of $\bar e$.\footnote{$S$ is concave if and only if 
$\bar e$ is convex. See \cite{MP} or Lemma 6.7 in \cite{Z1}.}

For viscous systems of conservation laws
\be\label{vcons}
u_t+f(u)_x=(B(u)u_x)_x,
\ee
shock waves are heteroclinic traveling waves. 
It is known \cite{GZ,ZS} that 
viscous stability is a stronger condition
than inviscid stability.
However, up to now it is not known {\it whether this
logical implication of inviscid by viscous stability
is strict for any interesting general class of constituents $(f,B)$}.
Destabilization by viscous effects of an inviscidly stable shock
wave would be a physically interesting phenomenon if it occurs,
and similar phenomena arising in Orr--Sommerfeld theory and incompressible
flow make this not implausible.
On the other hand, if this could be shown not to occur, that would be
equally interesting, and would greatly simplify the verification of
stability, reducing this to the study of the simpler inviscid Lopatinski
function, a linear-algebraic quantity, rather than the viscous Evans function,
defined in terms of solutions of an associated eigenvalue ODE.
In this connection, we were wondering {\it whether the existence of a 
viscosity-compatible entropy, $\eta$ as above with now also
$$ 
d^2\eta\ \! B\ge 0,
$$ 
might imply equivalence of viscous and inviscid instability}.

As pointed out in \cite{Z1,Z2,TZ2,TZ3}, the same analysis \cite{GZ,ZS}
showing that viscous stability implies inviscid stability, shows 
that{, 
as long as the viscous profile remains transverse as the 
intersection of the invariant manifolds as which it is defined,}
viscous
destabilization of a stable inviscid shock necessarily must occur through
the passage through the imaginary axis into the unstable
(positive real part) half plane of a complex conjugate pair of eigenvalues
of the linearized operator about the wave, a {leading-mode 
Hopf bifurcation}.
Here, we are imagining a transition from viscous stability to instability
as some bifurcation parameter, 
typically shock amplitude, is varied. %
``Galloping instabilites'' arising through leading-mode Hopf bifurcations
 are familiar
in detonation theory \cite{TZ4}, but have up to now not been observed
in the shock wave context.
{\it Could it be that the existence of a viscosity-compatible entropy precludes
complex, or even just purely imaginary eigenvalues? Or could it possibly imply
a ``[strong] principal eigenvalue property'' that
any non-stable eigenvalue, of the linearized operator about the wave,
with largest possible real part must be real [and simple]?} In both cases,   
leading-mode Hopf bifurcation would be impossible.
{\it Could it even be that existence of a viscosity-compatible entropy 
implied both transversality of the shock profile
and impossiblity of leading-mode Hopf bifurcation?} In that case, viscous and inviscid stability would coincide.


\medskip

In this paper, we examine these and related questions
using a combination of analytical and numerical techniques.
 {The following are the main results of this examination:}

\medskip
{\it

 {{\bf Theorem A (on gas dynamics).} (i) There exist equations
of state $e=\bar e(\tau,S),\tau>0,S\in\R,$ that 
(a) satisfy the standard assumptions \eqref{j1}--\eqref{j4}, 
\eqref{g1}--\eqref{g6}, and \eqref{h1}--\eqref{h4} detailed below, 
with in particular, the largest and smallest eigenvalues simple
and genuinely nonlinear, 
entropy and shock speed monotone along the forward 
$1$- and $3$-Hugoniot curves, and  
a global concave entropy function $S(\tau,e)$
defined on $\tau, e > 0$, and 
(b)
admit inviscidly unstable shock waves.

(ii) The equation of state 
\be\label{ceg}
\bar e(\tau,S)= \frac{e^{S}}{\tau}+ C^2 e^{ S/C^2 -\tau/C}, 
\quad C>\! >1,
\ee
is an example for (i). 
}

\medskip

{{\bf Numerical Observation A (on gas dynamics).}  For the equation of state 
\eqref{ceg},
in all cases we investigated numerically, (a) viscous [in]stability is equivalent
to inviscid [in]stability,\\
(b)  the viscous-stability problem has no non-zero imaginary eigenvalues; in particular, 
transitions from stability to instability occur exclusively by real eigenvalues    
passing through the origin, and\\
(c) in situations of instability, the eigenvalue with largest real part is real 
and simple.   
} 

\medskip

{{\bf Numerical Observation B (on general systems).} There exist viscous 
systems \eqref{vcons} of conservation laws, endowed with a compatible entropy, that \\
(a) admit shocks that are inviscidly stable, but viscously unstable,\\
(b) the viscous-stability problem sometimes does have non-zero imaginary eigenvalues, 
while\\
(c) in all situations of instability we investigated numerically, 
the eigenvalue with largest real part is real, and 
transitions from stability to instability occur exclusively by real eigenvalues    
passing through the origin,\\
(d) in some cases, there are an even number of unstable (and all real) eigenvalues, 
and\\
(e) in some cases the eigenvalue with largest real part is not simple.  
} 
}

\medskip

{Our analysis 
divides roughly into an inviscid  and a viscous part, which we now sketch.

\subsection{Inviscid  {analysis}}
In the first part of the paper, we revisit the inviscid stability analysis
of Erpenbeck--Majda \cite{Er1,M1,M2,M3} and reconcile it with the analysis of
Smith \cite{S,MP} on nonuniqueness of Riemann solutions, a related
but slightly weaker condition than instability of an individual
component shock.
Let 
$$
p=\hat p(\tau, e):=-\bar e_\tau(\tau, \hat S(\tau,e))
$$
denote the pressure function determined by inversion of the equation of state 
$e=\bar e(\tau,S)$ for $\tau$ fixed.\footnote{A pressure law 
$p=\hat p(\tau, e)$ may also be 
considered in the absence of an equation
of state, being sufficient by itself to close the equations of gas
dynamics \cite{MP,Sm}. For the observation that any pressure law can to some
(interesting!) extent be interpreted as being associated with an equation of state,
cf.\ Sec.\ 3.3.}

Smith formulates criteria on the equation of state, 
amounting to a {\it weak condition}
\be\label{weakintro}\tag{Weak}
-\frac{\bar e_{\tau s}}{ \bar e_S \bar e_{\tau \tau }}
< - \frac{ 2} {  \bar e_\tau }, 
\;
\hbox{\rm or, equivalently,} \;
\hat p_\tau < \frac{p\hat p_e}{2},
\ee
a {\it medium condition}, 
\be\label{mediumintro}\tag{Medium$_{\rm U}$}
-\frac{\bar e_{\tau S}}{ \bar e_S \bar e_{\tau \tau }}
<- \frac{\frac{\bar e_\tau ^2}{2e \bar e_{\tau \tau}}+1} {  \bar e_\tau },
\;
\hbox{\rm or, equivalently,} \;
\hat p_\tau < \frac{p^2}{2e},
\ee
and a {\it strong condition},
\be\label{strongintro}\tag{Strong}
-\frac{\bar e_{\tau s}}{ \bar e_S \bar e_{\tau \tau }}
< - \frac{ 1} {  \bar e_\tau }, 
\;
\hbox{\rm or, equivalently,} \;
\hat p_\tau < 0,
\ee
and derives in particular that
\eqref{mediumintro} is  
equivalent under 
mild additional hypotheses\footnote{Cf.\ Theorem 1.1,
within which we include Smith's result for the reader's ease.}
to uniqueness of Riemann solutions.

Recall now some standard properties of  
a pressure law
$p=\hat p(\tau,e),\tau, e\in \RM^+$:

\be \label{j1}\tag{J1}
\hbox{\rm $p=\hat p(\tau,e) > 0$.  }
\quad \hbox{\rm (Positivity)}
\ee
\be \label{j2}\tag{J2}
\hbox{\rm $(\partial_\tau -p\partial_e) \hat p<0$.}
\quad \hbox{\rm (Hyperbolicity)}
\ee
\be \label{j3}\tag{J3}
\hbox{\rm $(\partial_\tau -p\partial_e)^2 \hat p>0.$}
\quad \hbox{\rm (Genuine nonlinearity)}
\ee
\be \label{j4}\tag{J4}
\hbox{\rm $ \hat p_e >0.$}
\quad \hbox{\rm (Weyl condition)}
\ee

The following is the technical centerpiece of Part I of this paper.

\bt\label{stabthm}
Under hypotheses \eqref{g1}--\eqref{g6}, \eqref{h1}--\eqref{h2} of Smith (see Sec.\ 3.2 below) 
on the equation of state $\bar e$, or, more generally, assumptions
\eqref{j1}--\eqref{j4} on the pressure law $\hat p$,
positivity of the signed Lopatinski determinant (see Sec.\ 2.1)
(for all shocks) 
is equivalent to
\be\label{mediumstab}\tag{Medium$_{\rm S}$}
-\frac{\bar e_{\tau S}}{ \bar e_S \bar e_{\tau \tau }}
<- \frac{-\frac{\bar e_\tau}{\sqrt{2e \bar e_{\tau \tau}}}+1} {  \bar e_\tau },
\;
\hbox{\rm or, equivalently,} \;
\hat p_\tau < cp/\sqrt{2e}= p \sqrt{  \frac{p\hat p_e-\hat p_\tau}{2e}  },
\ee
while {\rm(\cite{Sm}:)}\! uniqueness of Riemann solutions (for any data) is equivalent to \eqref{mediumintro}.
The four conditions are related by
$$
\hbox{
\eqref{strongintro} $\Rightarrow$
\eqref{mediumintro} $\Rightarrow$
\eqref{mediumstab} $\Rightarrow$
\eqref{weakintro}.
}
$$
In particular, condition \eqref{strongintro} by itself
is sufficient to imply stability of all shocks,  while
violation of \eqref{weakintro} implies existence of unstable ones.
\et

Theorem A is a corollary of Theorem 1.1 and the following finding.

\bt 
Equation of state \eqref{ceg} satisfies \eqref{g1}--\eqref{g6}, \eqref{h1}--\eqref{h2},
and violates \eqref{weakintro}.
\et

Both the 
particular implication 
\eqref{mediumintro} $\Rightarrow$
\eqref{mediumstab} 
in Theorem 1.1 and the general meaning of the signed Lopatinski determinant
are elucidated by the following fact, which holds for arbitrary systems.

\bt
If the signed Lopatinski determinant $\Delta(\alpha)$ of a family 
$(U_-(\alpha),U_+(\alpha)),\alpha \in\R$ of Lax shocks undergoes a sign change
at some critcal value $\alpha_*\in\R$, 
the Riemann problem loses uniqueness near the initial data 
$(U_l,U_r)=(U_-(\alpha_*),U_+(\alpha_*))$.
\et

\br\label{gardnerrmk}
A general polytropic equation of state 
$\bar e(\tau,S)= {e^{S/c_{v}}}/{\tau^{\Gamma} }$,
$\Gamma, c_v>0$,
satisfies \eqref{strongintro}, or $\Gamma^2<\Gamma(\Gamma+1) $.
By Theorem \ref{stabthm}, we see that a local counterexample
$\bar e(\tau,S)={e^{S}}/{\tau}+f(S)$, $f'>\! >1$, proposed by C. Gardner (final
paragraph of \cite{G}) is incorrect, as this equation of state
likewise satisfies \eqref{strongintro}, by essentially the same computation.
Nonetheless, Gardner's larger assertion, that ``... all of the situations we have 
considered are logically possible, without violating the fundamental condition for local thermodynamic stability, namely that $e$ is a convex function of
$\tau $ and $S$,'' turns out to be correct, as shown by example \eqref{ceg}.
\er

\br
Taylor expanding \eqref{ceg} in $C$, 
%
%
we obtain also a simplified local example 
\be\label{lceg}
\bar e(\tau,S)= e^S/\tau +S- C\tau +\tau^2/2,
\quad C>\! >1,
\ee
of a convex equation of state permitting instabilities,
for which Hugoniot curves, monotonicity of $S$ and $\sigma$,
etc. may be computed explicitly, giving a concrete illustration
of the theory.
This does not satisfy \eqref{j1},
but may be treated by similar techniques as used to establish 
Theorem A; see Proposition \ref{locprop} below.
On the other hand, we find that the seemingly similar example
$$
\bar e(\tau,S)= e^S/\tau - C\tau +\tau^2/2,\quad C>\! >1,
$$
does not permit instabilities, despite violating \eqref{weakintro};
see Proposition \ref{simpleprop}.
\label{locexrem}
\er

\br
Example \eqref{ceg}
is the more surprising in view of recent results \cite{LV}
in nearby settings showing that existence of a convex entropy implies stability.
%
Indeed, Theorem A appears at first sight to contradict Theorem 2 of \cite{LV},
which asserts that, under some mild technical assumptions,
existence of a convex entropy,
simplicity of extremal characteristics, and monotonicity of entropy and 
shock speed along forward Hugoniot curves together imply inviscid stability of
arbitrary-amplitude extremal (i.e., $1$- or $3$-) shocks.
This would be very interesting to resolve; at the least, it shows that
there is a very narrow window between negative and positive results.
\er

\br\label{checkrmk}
Smith's weak and strong conditions are phrased in terms of
the gas law $e=\check e(\tau,p)$ obtained by inverting
$p=\hat p(\tau, e)$ with respect to $e$,
under the assumption that $\hat p_e>0$ (\eqref{j4} below).
Our versions are equivalent
to those of Smith if and only if $\hat p_e>0$, 
through the relation $\check e_\tau= -\frac{\hat p_\tau}{\hat p_e}$,
but remain valid also in the general case, as Smith's therefore do not.
Similar observations are made in \cite{MP} regarding the logical ordering
of Smith's conditions, which again requires $\hat p_e>0$.
Thus, our conditions above are in fact extensions of Smith's conditions 
into the realm $\hat p_e\leq 0$, 
or, equivalently (by $\bar p_S= \hat p_e (de/dS)=\hat p_e T$), 
$\bar p_S\leq 0$.
In this case, other aspects of Smith's global theory break down,
in particular, the conclusions of
Proposition \ref{smithprop} and, thus, Theorem \ref{stabthm}.
However, we find that
individual $1$-shocks $(U_-,U_+)$ are stable for $\hat p_e|_{U_+}\leq 0$; 
see Corollary \ref{t:4.1}.
\er

\subsection{Viscous  {analysis}}
In the remainder of the paper, we investigate viscous stability of
the above and other systems via numerical Evans function computations,
seeking viscous instabilities not predicted by inviscid theory, and
in particular  {
purely imaginary eigenvalues}
of the linearized operator about the wave.
The {\it Evans function}, as described for example in
\cite{AGJ,PW,GZ}, is an analytic function $D(\lambda)$ 
associated with a viscous shock wave,
defined on the nonstable complex half-plane $\{\lambda: \, \Re \lambda\ge 0\}$,
whose zeros encode the stability properties of the wave,
in particular, corresponding away from the origin
with eigenvalues of the linearized operator about the wave.
See \cite{Z1,Z2,Z3,Z5,TZ2,TZ3,SS,BeSZ} and references therein for more detailed
discussions of the relation
between the Evans function, spectral, linearized, and nonlinear 
shock stability, and Hopf bifurcation.

As described in \cite{BrZ,HuZ1,Z4,Z5}, the Evans function may be efficiently
computed by shooting methods as a Wronskian evaluated at $x=0$ of
decaying modes at $x=\pm \infty$ of the eigenvalue ordinary
differential equation (ODE) associated with the linearized operator about
the wave, and this has by now been carried out successfully for a number of
interesting systems/waves arising in continuum mechanics using the
standard Runga Kutta 4-5 (ode45) routine supported in MATLAB 
in conjunction with the exterior product or polar coordinate algorithms
\cite{Br,BrZ,HuZ1}
of the MATLAB-based Evans function package STABLAB \cite{BHZ2},
in particular
for the computationally intensive problem of stability of detonation waves
\cite{HuZ2,BZ1,BZ2,BHLyZ2}).

However, the equations of state \eqref{ceg} and \eqref{lceg}, due to the
separation of scales introduced by the large parameter $C>\! >1$ and
exponential dependence on parameters, induce a computational difficulty
far beyond any of the systems so far considered.  In particular, computations
with RK45 were not practically feasible, even on a parallel supercomputing machine.
To carry out Evans function computations for this system, rather, we found
it necessary for the first time to use an ODE solver designed for stiff systems,
namely, the ode15s routine supported in MATLAB, a 
recently-developed algorithm based on modern 
numerical differentiation formula (NDF) methods.

To indicate the relative stiffness of the systems associated with
\eqref{ceg}--\eqref{lceg} as compared to that of previously considered cases, 
the performance of Evans computations with
ode45 and ode15s are similar for all previously computed examples; however,
in the present case, ode15s outperforms ode45 by 2-3 orders of magnitude,
yielding performance that is not only feasible, but in the general range
seen for nonstiff computations.

The results of these numerical computations, gathered above 
as ``Numerical Observation A (on gas dynamics)'' 
hold also for equation of state \eqref{lceg}. 
\par\medskip
The other results, presented above 
as ``Numerical Observation B (on general systems)'', 
are obtained by computations 
on an artificially constructed
$3\times 3$ system 
of viscous conservation laws \eqref{vcons}\footnote{
The same size as the equations of gas dynamics.}
of a form suggested
by Stefano Bianchini \cite{Bi}.
Specifically, as shock amplitude is increased, several 
eigenvalues cross the origin into the unstable half-plane, some of
them eventually coalescing, splitting off of the real axis as a complex
conjugate pair, and crossing back into the stable half-plane through 
the imaginary axis at other points than $0$.  {This is a Hopf bifurcation,
though not a leading-mode one.} Note finally that  
the examples with an even number of unstable eigenvalues exemplify 
a scenario not detected even by signed versions of the inviscid stability 
condition, ---i.e., essentially, the signed Lopatinski determinant--- that find the 
parity of the number of unstable roots as described in \cite{GZ,Z1,Z2}.

\subsection{Discussion and open problems}
In conclusion, we have confirmed the claim of C. Gardner \cite{G} made
almost 50 years ago, but since then apparently not paid much attention to, 
that the equations
of gas dynamics, even under all of the usual thermodynamical and structural
assumptions associated with a ``normal'' gas,
{\it can support unstable inviscid shock waves}.
Moreover, we have shown this in a somewhat stronger sense than that apparently
envisioned by C. Gardner, demonstrating that this can be accomplished not only
locally, but globally, with a convex equation of state
$e=\bar e(\tau,S)$ satisfying expected asymptotics
as $S\to \pm \infty$ or $\tau\to 0,+\infty$.
In the process, we shed new light on the important work of R. Smith \cite{S}
and Menikoff--Plohr \cite{MP} on uniqueness of Riemann solutions.

This in a sense validates the inviscid stability theory of Erpenbeck--Majda,
showing that stability does not hold ``automatically'' for equations of state
satisfying all conditions of a classical gas.
On the other hand, we find so far no distinction between viscous and inviscid
stability for shock waves in gas dynamics.
We do, however, show by example that such a distinction, and in particular
the new phenomenon of Hopf bifurcaton, can occur for
a $3\times 3$ system with global convex entropy.
 {
Whether such phenomena arise also for Lax shocks of physical systems such 
as gas dynamics, magnetohydrodynamics (MHD), or viscoelasticity, must be left to
further investigation.}
We point in particular to the viscous-stability analysis of full (nonisentropic) MHD, 
so far not carried out even for the classical polytropic equation of state,
as a promising candidate and an important problem for further study.

Regarding the initial transition from stability to instability, 
we find no distinction between viscous and inviscid behavior for any of
the systems we consider.
Likewise, we find no counterexample to
the (weak) ``principal eigenvalue property'' that the largest unstable eigenvalue
of the linearized operator about the wave be real.
It would seem most interesting to know whether this
property holds in any generality.

Finally, recalling that shock instability in the physical literature
has often been associated with the phenomenon of phase transition
\cite{Er1,BE,MP}, we repeat that our examples here do not involve phase transition;
rather, the principal mechanism for instability, at a technical level,
i.e., the feature leading to violation of \eqref{weakintro}, is 
{\it stiffness}, or appearance
of multiple scales via the parameter $C>\! >1$.
Indeed, computing
$p(\tau,S)=-\bar e_\tau(\tau,S)= \frac{e^{S}}{\tau^2} + C-\tau$
for example \eqref{lceg} reveals a similarity 
to  a ``stiffened'' or ``prestressed'' equation
of state sometimes used to model shock behavior near the liquid state. 
For example, water is often described by 
\begin{equation}\label{stiffenedeos}
p=\Gamma \rho e - \gamma P_0
\end{equation}
with $\Gamma$ and the base stress $P_0$ determined empirically
\cite{IT,HVPM};
similar techniques are used to model liquid argon, nickel, mercury, etc.
\cite{H,CDM}.
The careful numerical study and categorizaion of behavior
of these and other more exotic 
equations of state such as van der Waals, Redlich-Kwong, etc., both
with and without viscosity,
appears to be another important direction for further study, and
one in which surprisingly little has so far been done.
In particular, a viscous counterpart of the inviscid analyses
of stability carried out for van der Waals gas dynamics in \cite{ZMRS,BM}
would be a valuable addition to the literature.

\bigskip
\centerline{\bf PART I. INVISCID STABILITY FOR GAS DYNAMICS}
\par\medskip
Sections 2 and 3, the first two of the three sections of this first principal part of the paper 
are devoted to establishing a network of algebraic and geometric conditions
on both individual shock waves and Hugoniot curves that alone or together
characterize (un)stable shock waves or imply or preclude the existence of
such waves in the sense of necessary or sufficient conditions. 
Part I culminates in Section 4 that shows the existence of inviscidly unstable 
shock waves in gas dynamics with certain equations of state.

\section{Signed Lopatinski criterion and relations to other conditions}
In this section, we define the signed Lopatinski criterion, evaluate it 
for gas dynamics with very general equations of state, and characterize its
connections with other conditions on the equation of state.

\subsection{The signed Lopatinski determinant}
We begin by recalling the Lopatinski condition of 
Erpenbeck--Majda \cite{Er1,M1,M2,M3} for stability of inviscid shock waves.
Consider a general system of conservation laws
\be\label{gencon}
f^0(w)_t+f(w)_x=0,\quad w\in \RM^n
\ee
and a {\it Lax $p$-shock} \cite{Sm,Se2,Se3}
 $U_\pm$ with speed $\sigma$, satisfying the
{\it Rankine--Hugoniot conditions}
\be
\label{rh}\tag{RH}
\sigma[f^0]=[f],
\ee
{\it hyperbolicity} 
\be\label{hyp}
\sigma(A_\pm) \, \hbox{\rm real},
\quad A_\pm:= (f^0_w)^{-1} f_w (U_\pm), 
\ee
of endstates $U_\pm$ and the {\it Lax conditions}
\ba\label{laxcon}
a_1^-&\leq \cdots \leq a_{p-1}^-<\sigma<a_p^-< \cdots \leq a_n^-,\\
a_1^+&\leq \cdots <a_{p}^+<\sigma<a_{p+1}^+\leq \cdots \leq a_n^+,
\ea
where $a_j^\pm$ are the (ordered) eigenvalues 
of $A_\pm:= (f^0_w)^{-1} f_w (U_\pm) $.

In this setting,
one-dimensional linearized and (under additional mild structural assumptions)
nonlinear inviscid stability is 
{addressed} (see \cite{M1,MZ1,MZ2,BS}) by the {\it Lopatinski condition}
\be\label{lop}
\delta:=\delta(U_-,U_+):=
\det(A^0_-r^-_1,\dots, A^0_-r^-_{p-1}, [f^0],A^0_+r_{p+1}^+,\dots, A^0_+r_n^+)\ne 0,
\ee
where $r_j^-$, $j=1,\dots, p-1$
are the ``outgoing'' modes of $A_-:= (f^0_-)^{-1} f_w (U_-) $
relative to the shock, i.\ e., a basis of the invariant space
associated with eigenvalues $a_j^-$ such that $a_j^--\sigma<0$ and
$r_j^+$, $j=p+1,\dots, n$
are the ``outgoing'' modes of 
$A_+:= (f^0_w)^{-1} f_w (U_+) $, relative to the shock, i.\ e., 
a basis of the invariant space associated with eigenvalues 
$a_j^+$ such that $a_j^+-\sigma>0$,
and $A^0_\pm :=f^0_w$.
Implicit in this description is the 
defining property of a Lax $p$-shock that
the outgoing eigenvectors on 
the left (right) indeed number $p-1$ ($n-p$).

\br
Note that we have allowed here by the inclusion of $f^0$ a somewhat more
general formulation than in \eqref{cons}.
In the simplest case $f^0(w)=w$ occurring in many applications,
the computations considerably simplify.
However, as we shall see in the full gas case, there can in some cases
be an advantage in being able to choose coordinates more flexibly.
A short proof of the case $f(w)=w$ may be found in \cite{ZS}.
Redefining $\tilde w:=f^0$, we obtain
 $\tilde A_+=f_{\tilde w}=f_u (f^0_u)^{-1}$, yielding the result
by the observation that eigenvectors $\tilde r_j^+$ of $\tilde A_+$
are equal to $A^0_+r_j^+$.  See also \cite{Se2,Se3,M1}.
\er

\bd
Consider a given p-shock $(U_-,U_+)$ for a system \eqref{gencon}. Assume that there
is  a continuous one-parameter family $(U^-(\alpha), U^+(\alpha)),0<\alpha \le 1,$
of p-shocks with 
$U^\pm(1)=U_\pm$
and $U^\pm(0+)=U_*$ for some state $U_*$, that the bases $\{r_1^-,\ldots,r_{p-1}^-\},
\{r_{p+1}^+,\ldots,r_n^+\}$ and thus 
$$
\delta(\alpha):=\delta(U^-(\alpha),U^+(\alpha))
$$
are (chosen) continuous in $\alpha$, and that the expression 
\be
\Delta\equiv\delta(1)/\lim_{\alpha\searrow 0}\sgn\delta(\alpha)), 
\ee
makes sense and has 
the same sign for all such homotopies. Then $\Delta=\Delta(U_-,U_+)$
is called   
a {\it signed Lopatinski determinant} 
of $(U_-,U_+)$.
\ed

Obviously, the signed Lopatinski determinant $\Delta$ 
is well defined (--- more precisely, while the Lopatinski determinant $\delta$ is defined up 
to a non-vanishing factor, the signed Lopatinski determinant $\Delta$ is defined up to a 
positive factor ---)
at least whenever
\ba\label{lopwell}
&\hbox{the system's state space is simply connected,}\\ 
&\hbox{$a_p$ is an everywhere simple,
genuinely nonlinear eigenvalue, and}\\ 
&\hbox{all $p$-shocks with arbitrary fixed left state $U_-$ group as 
regular directed curves $S_p(U_-)$.}
\ea 

In the current Part I of this paper we use the signed Lopatinski determinant to detect instability
according to the following evident fact.

\bc
If a given hyperbolic system \eqref{gencon} satisfying \eqref{lopwell} admits 
a p-shock with $\Delta<0$, then it also admits a p-shock with $\delta=0$.
\ec

\begin{proof}
For a corresponding homotopy, obviously
$$
\Delta(U^-(\alpha),U^+(\alpha))>0 \quad\hbox{for sufficiently small }\alpha>0,
$$
and the existence of an $\alpha_*$ with 
$\Delta(U^-(\alpha_*),U^+(\alpha_*))=0$
follows by continuity.
\end{proof}

For the (deeper) meaning of the signed Lopatinski determinant in the viscous context,
see Part II, Sec.\ 6.

\br
As noted in \cite{Fo} in the context of gas dynamics, and \cite{Se1}
more generally, a change in the sign of $\Delta$ may also be 
interpreted in terms
of multi-dimensional inviscid theory, signalling a transition from weak
multi-dimensional stability in the sense of Majda \cite{M1,M2,M3}
to exponential instability; see also \cite{MP,Z1,Z3}.
\er

\subsection{The Lopatinski condition in  gas dynamics}
Specializing to gas dynamics, we now compute the Lopatinski determinant and signed
Lopatinski determinant explicitly.
\par\medskip\noindent
{\bf Isentropic gas dynamics.}
The isentropic Euler equations in Lagrangian coordinates are
\ba\label{ieuler}
\tau_t-v_x&=0,\\
v_t+ p_x&=0,\\
\ea
with a given pressure law $p=p(\tau)$, where $\tau$ denotes specific volume,
$v$ velocity, and $p$ pressure.
From 
$
A_\pm=\begin{pmatrix}
0& -1\\
p_\tau &0
\end{pmatrix}_\pm,
$
we readily find that
$a_j^\pm= -c, c$, $c:=\sqrt{-p'(\tau)}$, and
$r_j^\pm=
\begin{pmatrix}
1\\
-a_j^\pm
\end{pmatrix},
$
so that hyperbolicity corresponds to $p'<0$. 

Moreover, \eqref{rh} gives $\sigma[\tau]=-[v]$, $\sigma[v]=[p]$,
yielding $\sigma^2=-[p]/[\tau]$,
whence, for a $1$-shock,
$\sigma=-\sqrt{-[p]/[\tau]}$,
and so we have 
$
[f^0]=
[\tau]
\begin{pmatrix}
1\\
\sqrt{-[p]/[\tau]}
\end{pmatrix}.
$
Thus,
$$
\delta= [\tau] \det
\begin{pmatrix}
1& 1\\
-\sigma
& -c
\end{pmatrix}
=
[\tau] (\sigma -c)
<0\quad\hbox{and}\quad
\Delta>0
$$
for all Lax shocks. Notice that condition \eqref{lopwell} is satisfied trivially. 
Note that $[\tau]\ne 0$, since otherwise $[v]=0$
by \eqref{rh}, and there is no shock.
\par\medskip\noindent
{\bf Full gas dynamics.}
In Lagrangian coordinates, the full (nonisentropic)
Euler equations are
\ba\label{euler}
\tau_t-v_x&=0,\\
v_t+ p_x&=0,\\
(e+v^2/2)_t+(vp)_x&=0,\\
\ea
with a given pressure law $p=\hat p(\tau, e)$,
where $\tau$ denotes specific volume,
$v$ velocity, $e$ specific internal energy, and $p$ pressure.
We assume here that the pressure law originates from
a thermodynamic equation of state 
\be\label{eos}
e=\bar e(\tau, S),
\ee
where $S$ is entropy, through the relation $p=-\bar e_\tau$,
assuming that temperature $T=\bar e_S$ is positive, so that
\eqref{eos} may be solved for $S=\hat S(\tau, e)$ to obtain
$\hat p(\tau,e)=\bar e(\tau, \hat S(\tau,e))$.

Alternatively, for smooth solutions, we have the simple entropy form \cite{Sm}
\ba\label{ent_euler}
\tau_t-v_x&=0,\\
v_t+ p_x&=0,\\
S_t&=0,.\\
\ea
where 
$p=p(\tau,S)=-\bar e_\tau (\tau,S)$ for equation of state \eqref{eos}.
Choosing $w=(\tau,v,S)^T$, we thus have
$$
A_\pm=
\begin{pmatrix}
0 & -1 & 0\\
p_\tau& 0 & p_S\\
0&0&0
\end{pmatrix}_\pm,
$$
and, essentially by inspection,
\be\label{chars}
a_1=-c,
\; a_2=0,
\;
a_3=c,
\qquad
c:=\sqrt{-p_\tau},
\ee
and
$
r_1^+=
\begin{pmatrix}
1\\ 
c
\\ 0\\
\end{pmatrix}$,
$r_2^+=
\begin{pmatrix}
-p_S\\ 0\\ p_\tau
\end{pmatrix}$,
$r_3^+=
\begin{pmatrix}
1\\ 
-c
\\ 0\\
\end{pmatrix}$,
so that hyperbolicity corresponds to $p_\tau=-\bar e_{\tau \tau}<0$.
Computing
$
A_\pm^0=
\begin{pmatrix}
1 &  0&0\\
0 & 1&0\\
\bar e_\tau &v&\bar e_S
\end{pmatrix}_\pm
=
\begin{pmatrix}
1 &  0&0\\
0 & 1&0\\
-p &v&
T
\end{pmatrix}_\pm$,
we obtain
$$
A_+^0 r_1^+=
\begin{pmatrix}
1\\ 
c
\\ 
-p
+vc
\\
\end{pmatrix},
\quad
A_+^0 r_2^+=
\begin{pmatrix}
-p_S\\ 0\\ 
p
p_S + 
T
p_\tau
\end{pmatrix},
\quad
A_+^0 r_3^+=
\begin{pmatrix}
1\\ 
-c
\\ 
-p
-vc
\\
\end{pmatrix}.
$$

Next, from \eqref{rh}, we have, taking without loss of generality $v_-=0$,
\be\label{gasRH}
[v]=-\sigma[\tau],
\quad
[e+v^2/2]=\sigma^{-1}[vp]=
-p_+[\tau],
\quad
\sigma^2 =-[p]/[\tau],
\ee
so that, for a $1$-shock (recalling the relation $p_\tau= -c^2$),
$
\delta=
[\tau]
\det
\begin{pmatrix}
1& -p_S & 1\\
-\sigma & 0 & -
c
\\
-p & 
p
p_S - 
T
c^2
 & 
-p
-vc
\\
\end{pmatrix},
$
or
\be\label{fullcomp}
\delta
= [\tau]
(p_S c [p] +T
c^2(\sigma-c)),
\ee
where all quantities are evaluated at $U_+$, and
$ c:=\sqrt{ -p_\tau} $ as above denotes sound speed.
Note that $[\tau]\ne 0$, or else $[v]=0$, $[p]=0$,
and thus $[e]=0$, all by \eqref{rh}.
But, $[e]=[\tau]=0$ and $[\tau]=0$ by $\bar e_S>0$ implies that
$[S]=0$ and there is no shock.

Note that, in the decoupled case $p_S=0$, we get
$\delta = [\tau](
Tc^2(\sigma-c))< 0 $, by $T>0$,
essentially reducing to the isentropic case.

{\it We assume, here and henceforth, that 
condition \eqref{lopwell} holds for the Euler equations \eqref{euler}.
Note that whenever Smith's conditions \eqref{g1}--\eqref{g6}, \eqref{h1}--\eqref{h2}
(cf.\ Sec.\ 3.2) hold, this is automatic.}


\bt\label{signlem}
For equations \eqref{euler} with $p_\tau<0$
a Lax 1-shock satisfies   
the signed Lopatinski condition
$\Delta>0$ if and only if 
\be\label{dcond}\tag{Lop}
p_S  [p] < Tp_\tau (\sigma/c-1),
\ee
or, alternatively,
\be\label{lopalt}\tag{Lop$_{\rm alt}$}
\hat p_\tau [p] < -c \sigma,
\ee
where all quantities not in [.] are evaluated at $U_+$.
\et

\begin{proof}
By homotopy from the general case to the decoupled case. 
Equivalence of \eqref{lopalt} and \eqref{dcond} follows by $ p_\tau=-c^2=\hat p_\tau -p\hat p_e$
and $\hat p_e(\tau,e)=\bar p_S(\tau,S)/\bar e_S(\tau,S)$.
\end{proof}


We state the following elementary result mainly in order to illustrate the 
competition between isentropic and nonisentropic effects.

\bc\label{t:4.1}
For a gas with equation of state satisfying
$T, p>0$, $p_\tau<0$, and (``anti-Weyl'' condition) $p_S|_{U_+} \le 0$,
all Lax shocks satisfy the signed Lopatinski condition $\Delta>0$. 
\ec

\begin{proof}
In this case, each of the terms
$p_S c [p]$ and $-Tp_\tau (\sigma-c)$
on the righthand side of \eqref{fullcomp} are strictly negative,
by $T>0$ and $p_\tau<0$ (hyperbolicity),
so that 
$\sgn\ \delta=\sgn[\tau]$.
But $[\tau]<0$ for any 1-shock and $[\tau]>0$ for any 3-shock. 
\end{proof}

\br\label{prmk}
Though we have restricted for clarity to the case that the pressure relation
derives from a complete equation of state $e=\bar e(\tau,S)$, we could equally
well carry out the analysis within the framework \eqref{euler},
with $p=\hat p(\tau,e)$, to obtain the condition
$
\hat p_e[p] < c^2 (1-\sigma/c)
$
or, more simply, \eqref{lopalt}; see Appendix \ref{s:pcond}.
This could alternatively be obtained indirectly, by the observation
that, for any choice of positive temperature function $T=T(\tau,e)$
defined in the vicinity of $U_+$,
we can invert the Thermodynamic relation
$de= TdS- p d\tau$ to obtain a local entropy $S$ and equation of state
$e=\bar e(\tau,S)$; see Section \ref{s:entropy}.
\er

\br
Though the pressure law $p=\hat p(\tau,e)$ is sufficient to close the
Euler equations, the Navier--Stokes equations, involving the Fourier law,
require also the specification of a temperature law $T=\hat T(\tau,e)$.
Both can be obtained from an equation of state $e=\bar e(\tau,S)$,
as described above.
\er
\par\medskip\noindent
{\bf Relation to Majda's condition.}  
We now verify that, \emph{assuming the Weyl condition $p_S>0$}, 
our condition \eqref{dcond} agrees with
Majda's 1D stability condition (more precisely, the
signed version obtained by taking into account multi-dimensiional
considerations; see Remark B.6, p. 523--524, Appendix C, \cite{Z1}):
\be\label{mcond}
\frac{1+M}{\Gamma} >
\Big(\frac{\tau_-}{\tau} -1\Big)M^2,
\ee
where $\Gamma$ is the Gruneisen coefficient
$ \Gamma:= \frac{\tau p_S}{T}, $
and
$ M^2=\frac {(v-\sigma_E)^2}{c_E^2}, $
the subscript $E$ denoting Eulerian values.
Evidently, $\Gamma>0$ if and only if $p_S>0$.
%

Rewriting \eqref{mcond} as $ \frac {\Gamma} {1+M} < \frac{\tau}{-[\tau]M^2} $
using positivity of $\Gamma$
and substituting $\frac{\tau p_S}{T}$ for $\Gamma$,
we have
$$
 \frac {-[\tau] p_S} {T} < \frac{1+M}{M^2} .
$$
We may compare to our condition \eqref{dcond} by expressing
the (Eulerian) Mach number $M$ in terms of Lagrangian quantities.
Evidently, we must verify the identity
$\frac {T}{-[\tau] } \frac{1+M}{M^2} = \frac{ Tp_\tau (\sigma-c)} { c [p]} $, 
or
$$
\frac{1+M}{M^2} = \frac{ c (c-\sigma)} {\sigma^2}  . 
$$

In Eulerian coordinates, we have $\sigma_E [\rho]=[\rho v]$,
or $\sigma_E=\frac{v/\tau}{[1/\tau]}=\frac{-v}{[\tau]}=\sigma$,
in agreement with the Lagrangian value.
Moreover, $v=[v]=-\sigma [\tau]$ gives $v-\sigma=
-\sigma( [\tau]+1)$, so that, in case $\tau_0=1$, we have $v-\sigma=-\sigma\tau$.
On the other hand, $c_E=\tau c$, hence
$
M= \frac{ v-\sigma}{\tau c},
$
and the desired identity is
$
\frac{1+(v-\sigma)/\tau c}{
(v-\sigma)^2/(\tau c)^2 } 
=
\frac{ c (c-\sigma)} {\sigma^2}  ,  $
or
$
\frac{\tau c +(v-\sigma)}{
(v-\sigma)^2 } 
=
\frac{  (c-\sigma)} {\tau \sigma^2}  . 
$
Without loss of generality (by rescaling) taking
$\tau_-=1$, $v-\sigma=-\sigma \tau$, this becomes
$
\frac{\tau c-\sigma \tau }
{\sigma^2\tau^2}
=
\frac{  (c-\sigma)} {\tau \sigma^2}  ,
$
which is evidently true.

\br\label{majrmk}
Similarly as in Remark \ref{checkrmk}, we find that our condition \eqref{dcond}
is an extension of Majda's condition into the realm $p_S|_{U_+}\leq 0$
where Majda's condition no longer applies.
Indeed, for this case $\Gamma\leq 0$ and \eqref{mcond} fails always, whereas
we see by \eqref{dcond} that stability in fact always holds.
\er

\subsection{Sufficient global conditions for (in)stability}
%
The following sufficient conditions for satisfaction and violation, respectively,
of the signed Lopatinski condition are as easy to get as useful.
%
For definiteness, we restrict to the case of $1$-shocks.

\bpr\label{primeprop}
For a Lax $1$-shock with $[\tau]<0$, the signed Lopatinski condition $\Delta>0$
is implied by
\be\label{strong}\tag{Strong'}
-\frac{\bar e_{\tau S}}{ \bar e_S \bar e_{\tau \tau }}|_{U=U_+}
< \frac{ 1} {  [p] }.
\ee
The opposite, $\Delta<0$, is implied by failure of 
\be\label{weak}\tag{Weak'}
-\frac{\bar e_{\tau S}}{ \bar e_S \bar e_{\tau \tau }}|_{U=U_+}
< \frac{ 2} {  [p] } .
\ee
\epr

\begin{proof}
Because of $[\tau]<0$, $[p]$ is positive, and
dividing both sides of \eqref{dcond} by the positive quantity $c[p]$, 
we obtain
\be\label{dcond3}\tag{${\rm Lop}_1$}
-\frac{\bar e_{\tau S}}{ \bar e_S \bar e_{\tau \tau }}
< \frac{-\frac{\sigma}{c}+1} {  [p] }
= \frac{\frac{|\sigma|}{c}+1} {  [p] }. 
\ee
Recalling, \eqref{chars}, 
that the characteristic speeds for \eqref{euler}--\eqref{ent_euler}
are $-c$, $0$, $c$, we have for a Lax $1$-shock that at the right
state $U_+$, $-c<\sigma<0< c$, or
$\sigma<0$ and $|\sigma|<c$. 
Hence,
$
1<-\frac{\sigma}{c}+1<2,
$
from which we obtain
with \eqref{dcond3},
the stated stability conditions.
\end{proof}

\br\label{order}
Conditions \eqref{strong} and \eqref{weak} are global versions of
conditions \eqref{weakintro} and \eqref{strongintro} of \cite{S}.
Evidently, \eqref{weakintro}$\Rightarrow $\eqref{weak} and
\eqref{strongintro}$\Rightarrow $\eqref{strong}.
\er

\subsection{Monotonicity of the Hugoniot curve}
Using the formulation \eqref{dcond}, we show how stability of shocks 
relates to geometric properties of the Hugoniot curve,
recovering both 
the classical observations of Erpenbeck and Gardner \cite{Er1,G}
relating stability to geometric properties of the Hugoniot curve
and the conditions of Smith for uniqueness of Riemann solutions.
It remains to investigate monotonicity of the $1$- (or, equivalently,
the $3$-) Hugoniot curve.
From the assumption $\bar e_S=T>0$, we may solve $e=\bar e(\tau,S)$
to obtain $S=\hat s(\tau, e)$, and thereby 
\be\label{hatp}
p=\hat p(\tau,e):=- \bar e_\tau (\tau, \hat s(\tau, e)).
\ee
Combining the three equations \eqref{rh}, we readily obtain
the single equation
\be\label{hugrel}
\hat H(\tau,e) =[e] + (1/2)(p+p_-)[\tau]=0
\ee
for $(e,\tau)=(e_+,\tau_+)$, implicitly determining,
together with \eqref{hatp},
$e$ as a function of $\tau$ or vice versa.
Alternatively, using \eqref{gasRH} to relate different variables,
we could consider this as determining $p$ as a function of $v$,
or $\tau$ as a function of $s$, whose derivatives
in each case may be computed by a straightforward 
application of the Implicit Function Theorem.
We record the results here.

\bl
We have the change of variables formulae
\be\label{estrans}
\hat s_e=\frac{1}{\bar e_s},
\quad
\hat s_\tau=-\frac{\bar e_\tau}{\bar e_s},
\quad
\hat p_\tau= -\bar e_{\tau \tau} -\bar e_{\tau s}\hat s_\tau=
-\bar e_{\tau \tau} +\frac{\bar e_{\tau s} \bar e_\tau}{ \bar e_s},
\quad 
\hat p_e= -\bar e_{\tau s}\hat s_e=
 -\frac{\bar e_{\tau s}}{\bar e_s},
\ee
and, along the Hugoniot curve $H\equiv 0$, viewing $e$ and $p$
as functions
$e=e_H(\tau)$ and $p=p_H(\tau)$, 
\be\label{eform}
\frac{de_H}{d\tau}= \frac{-\frac12(\hat p_\tau[\tau]+p+p_-)}
{1+\frac12\hat p_e[\tau]},
\quad
\frac{dp_H}{d\tau}= \hat p_\tau + \hat p_e \frac{de_H}{d\tau}=
\frac{\frac12[p]\hat p_e-c^2}
{1+\frac12\hat p_e [\tau]} .
\ee
Alternatively, considering $v$ as a function of $p$ along $H\equiv 0$,
we have, taking $v_-=0$,
\be\label{vp}
-2v \frac{dv}{dp}=
[\tau] \frac {[p]\hat p_e-c^2 -\sigma^2 } {\frac12[p]\hat p_e-c^2}.
\ee
\el

\begin{proof}
Relations \eqref{estrans} and \eqref{eform} follow directly from
the Implicit Function Theorem applied to \eqref{hatp} and \eqref{hugrel},
where the final equality in \eqref{eform} follows from
$
\frac{dp_H}{d\tau}= \hat p_\tau + \hat p_e \frac{de_H}{d\tau}=
\frac{-\frac12(p+p_-)\hat p_e + \hat p_\tau}
{1+\frac12\hat p_e [\tau]} 
$
together with $c^2=-\hat p_\tau + p\hat p_e$.
From $v=-\sqrt{-[p][\tau]}$ (a consequence of \eqref{rh} and $v_-=0$, through
$[v]=-\sigma [\tau]$), we obtain,
differentiating implicitly $v^2=[p][\tau]$,
$
-2v \frac{dv}{dp}=
[\tau] +  [p]\frac{1}{{dp_H/d\tau}}=
[\tau]\frac{dp_H/d\tau -\sigma^2}{{dp_H/d\tau}}=
$,
verifying \eqref{vp} by \eqref{eform}.
\end{proof}

Carrying out these computations, one may obtain various necessary or
sufficient geometric conditions on the Hugoniot curve. 
In particular, we have the following result. 


\bc\label{monprop}
Monotonicity, 
$dv/dp<0$ 
at a point $U_+$ on the $1$-Hugoniot curve 
through $U_-$ is equivalent to
\be\label{medium}\tag{Monotone}
-\frac{\bar e_{\tau S}}{ \bar e_S \bar e_{\tau \tau }}
< \frac{\frac{\sigma^2}{c^2}+1} {  [p] },
\;
\hbox{\rm or, equivalently,} \;
\hat p_e [p]<  \sigma^2+c^2.
\ee
Moreover, we have the string of implications
\be\label{imps}
\hbox{
\eqref{strong}$\Rightarrow$\eqref{medium}$\Rightarrow$\eqref{dcond3}$\Rightarrow$\eqref{weak}.
}
\ee
\ec

\begin{proof}
From \eqref{vp}, we have evidently (by $v,[\tau]<0$) that 
$dv/dp<0$ is equivalent to 
$$
0>\frac{dp_H/d\tau -\sigma^2}{{dp_H/d\tau}}=
2+\frac{c^2-\sigma^2}{\frac12\hat p_e[p]-c^2} \, ,
$$
or, rearranging, $\hat p_e[p]< \sigma^2 + c^2$ as claimed,
yielding \eqref{medium} by $[p]>0$.
The logical  implications \eqref{imps} then follow by the Lax condition $|\sigma|<c$.
\end{proof}

\br\label{monrmk}
The implication 
\eqref{medium}$\Rightarrow$ \eqref{dcond3} is due to C. Gardner \cite{G}.
A slightly less sharp condition pointed out by Erpenbeck \cite{Er1} is
that monotonicity of $\tau$ with respect to $s$ along the forward
Hugoniot curve, 
$ p_S<\frac{2T}{-[\tau]}$
(see Remark \ref{vasprep}), or
\be\label{fwdmon}
-\frac{\bar e_{\tau S}}{ \bar e_S \bar e_{\tau \tau }}
< \frac{\frac{2\sigma^2}{c^2}} {  [p] },
\ee
also implies \eqref{dcond3}, 
by $|\sigma|/c<1$;
likewise, \eqref{fwdmon} implies \eqref{medium}
and uniqueness of Riemann solutions.
For other geometric implications of the various conditions,
see Theorem 4.5 \cite{MP}.
\er

\br
From \eqref{imps}, we obtain 
\eqref{strongintro}$\Rightarrow$\eqref{dcond3}, 
and
\eqref{strongintro}$\Rightarrow$\eqref{strong}.
\er

\section{Local conditions for (in)stability and the assumptions of Smith}
The above-obtained {\it global conditions} require knowledge of
the Hugoniot curve in order to compute, which is in practice quite
limiting.
Adapting the ingenious observations of Smith, we now show that,
under some simple structural assumptions on the Hugoniot curve, 
these global conditions
may be replaced by {\it local conditions} more convenient for analysis.
In this section, we establish Theorems 1.1 and 1.3. While the latter is proved 
in Subsection 3.4, Theorem 1.1 is an immediate consequence of Propositions
3.1 and 3.2 and Corollary 3.6. 

\subsection{Local conditions and their requirements}\label{s:local}
We identify the following properties, satisfied for most standard
equations of state.  
Denoting as the backward $1$-Hugoniot through $U_+$,
$H_1'(U_+)$,
the set of all left states $U_-$ connected to $U_+$ by a Lax $1$-shock,
and assuming that $H_1'(U_+)$ is a directed curve, we distinguish
\be \label{p1}\tag{P1}
\hbox{\rm $[\tau]<0$ on 
$H_1'(U_+)$,}
\ee
\be \label{p2}\tag{P2}
\hbox{\rm $p\to 0$ as $U$ progresses along $H_1'(U_+)$,} 
\ee
\be \label{p3}\tag{P3}
\hbox{\rm $e\to 0$ as $U$ progresses along 
$H_1'(U_+)$,
}   
\ee
\be \label{p5}\tag{P4}
\hbox{\rm $\tau$ is increasing and $p$ decreasing along 
$H_1'(U_+)$.
}
\ee

Recalling now conditions \eqref{strongintro}, \eqref{mediumintro},
\eqref{mediumstab}, \eqref{weakintro} of the introduction, 
we have the following detailed local characterizations of stability.

\bpr\label{stabmain}
(i) Assuming \eqref{p1}, \eqref{strongintro} is sufficient for either stability
of 
shock waves connecting $U_+$ to states along the 
backward $1$-Hugoniot or uniqueness of Riemann solutions involving 
states along the backward $1$-Hugoniot;
(ii) assuming \eqref{p2}, \eqref{weakintro} is necessary for stability or
uniqueness on the backward $1$-Hugoniot.
(iii) Assuming \eqref{p1}--\eqref{p3}, 
\eqref{mediumstab} is necessary and sufficient for stability 
and \eqref{mediumintro} is necessary and sufficient for
uniqueness along the backward $1$-Hugoniot through $U_+$,
with, moreover, 
\be\label{imps2}
\hbox{\rm
\eqref{strongintro} $\Rightarrow$
\eqref{mediumintro} $\Rightarrow$
\eqref{mediumstab} $\Rightarrow$
\eqref{weakintro}.
}
\ee
(iv) Assuming also \eqref{p5}, there is at most one stability transition
as $U$ progresses along the backward $1$-Hugoniot through $U_+$,
that is, if the signed Lopatinski determinant turns negative once, then 
it stays negative.
\epr

\begin{proof}
Assertion (i) follows from \eqref{strong} and the fact that $p>[p]$
for $p_-\ge 0$, assertion (ii) from \eqref{weak} in the limit as $p_-\to 0$.
In assertion (iii), necessity 
follows from \eqref{dcond}--\eqref{lopalt} and \eqref{medium},\footnote{
Here, we are using the standard fact that uniqueness of Riemann solutions
is equivalent to monotonicity of the $1$-Hugoniot curve
(hence, by symmetry, the $3$-Hugoniot)
in terms of pressure vs. velocity \cite{S,Sm,MP}.  }
together with the observation that, by \eqref{hugrel},
$e_+=-\frac12 p_+ [\tau]$ in the limit as $e_-, p_-\to 0$, so that,
for $e_->0$, 
$
\sigma^2=-[p]/[\tau]= p_+^2/e_+.
$

Sufficiency follows, similarly, 
from \eqref{dcond}--\eqref{lopalt} and \eqref{medium}, and the observation
that their righthand sides are monotone increasing with respect
$p_-$ and monotone decreasing with respect to $\tau_-$,\footnote{
Note that dependence on $(\tau_-,p_-)$ is through factors
$\sigma^2/[p]= -1/[\tau]$, $\sigma/[p]=1/\sqrt{-[p][\tau]}$,
and $1/[p]$.
}
so that its minimum is achieved at $p_-=e_-=0$ and 
$\tau_-=\tau_{\rm max}:= 2e_+/p_+ + \tau_+$;\footnote{
Here we are using that, by \eqref{hugrel}, $\tau_-$ is monotone decreasing
with respect to both $p_-$ and $e_-$, so that its maximum is achieved
at the minimum $(p_-,e_-)=(0,0)$.}
the implications \eqref{imps2}
follow readily from Lax condition $|\sigma|< c$.
If also \eqref{p5} holds, then 
the righthand sides of all conditions 
are decreasing along the backward $1$-Hugoniot curve, yielding (iv).
\end{proof}

\subsection{The assumptions of Smith}\label{s:smithcond}

We now recall the assumptions of Smith
and show that they imply \eqref{p1}--\eqref{p5}.  
These include, for $\tau>0$ and $S\in \RM$, the structural conditions:

\be \label{g1}\tag{G1}
\hbox{\rm $\bar e(\tau,S)> 0$.  }
\quad \hbox{\rm (Positive energy)}
\ee
\be \label{g2}\tag{G2}
\hbox{\rm $p=-\bar e_\tau(\tau,S) > 0$.  }
\quad \hbox{\rm (Positive pressure)}
\ee
\be \label{g3}\tag{G3}
\hbox{\rm $T=\bar e_s(\tau,S)> 0$.  }
\quad \hbox{\rm (Positive temperature)}
\ee
\be \label{g4}\tag{G4}
\hbox{\rm $p_\tau=-\bar e_{\tau \tau}(\tau,S) < 0$. 
\quad \hbox{\rm (Hyperbolicity)}
}
\ee
\be \label{g5}\tag{G5}
\hbox{\rm $p_{\tau \tau}=-\bar e_{\tau \tau \tau}(\tau,S) > 0$. 
\quad \hbox{\rm (Genuine nonlinearity)}
}
\ee
\be \label{g6}\tag{G6}
\hbox{\rm $p_S=-\bar e_{\tau S}(\tau,S) > 0$. 
}
\quad \hbox{\rm (Weyl condition)}
\ee
Assuming \eqref{g6}, we may solve $p=-\bar e_\tau (\tau,S)$ 
for $S=\check S(\tau,p)$, to obtain 
$e=\check e(\tau,p):=\bar e(\tau, \check S(\tau, p))$.

Besides the above structural assumptions, Smith
imposes the asymptotic conditions:

\be \label{h1}\tag{H1}
\lim_{s\to -\infty} \bar e(s,\tau)=0,
\quad
\lim_{s\to +\infty} \bar e(s,\tau)=\infty,
\ee
\be \label{h2}\tag{H2}
\lim_{s\to -\infty} -\bar e_\tau (s,\tau)=0,
\quad
\lim_{s\to +\infty} -\bar e_\tau (s,\tau)=\infty,
\ee

We shall not require Smith's further assumptions that 
\be \label{h3}\tag{H3}
\hbox{\rm
$\lim_{\tau\to 0^+} \check e(\tau,p)=\infty$ for $ p>0$,
}
\ee
and
\be\label{h4}\tag{H4}
\lim_{\tau \to +\infty} \bar e(s,\tau)=0,
\ee
which appear to be used for existence and not uniqueness in \cite{S}.

We complete our treatment with the following 
streamlined version of observations of \cite{W,S}.

\bpr[\cite{W,S}]\label{smithprop}
Assumptions \eqref{g1}--\eqref{g6}, \eqref{h1}--\eqref{h2}
imply \eqref{p1}--\eqref{p5}: in particular,
(i) the backward $1$-Hugoniot is $C^1$,
extending till $S\to -\infty$, $p\to 0$, $e\to 0$,
(ii) $\tau$ and $v$ increase and $p$, $S$, and $\sigma^2$ decrease 
along the backward Hugoniot.
\epr

\begin{proof}
Recall that the Lax condition is $c_-^2<\sigma^2<c_+^2$ (with the additional
condition $\sigma<0$ distinguishing $1$-shocks from $3$-shocks).
By standard hyperbolic theory \cite{Sm},
\eqref{g4}--\eqref{g5} imply that in the vicinity of $(\tau_+,S_+)$,
the backward $1$-Hugoniot set of states $(\tau_-,S_-)$
satisfying \eqref{hugrel} and the Lax condition 
consists of a smooth curve on which $[\tau]<0$, $[S]>0$.
Let us focus on the region $[\tau]<0$, therefore.
Rewriting the Hugoniot relation as
$e_+-e= (1/2)(p+p_+)[\tau]$, we see readily that, assuming \eqref{g1}--\eqref{g2},
$
\tau_- \le \tau_{max}:= 2e_+/p_+ + \tau_+.
$
whence $\tau_+ \le \tau_- \le \tau_{max}$.

We shall show that the backward Hugoniot curve in the vicinity of
$(\tau_+,S_+)$ extends as a smooth curve, monotone in
$\tau$ and $S$, globally on
$\tau\in [\tau_+,\tau_{\rm max})$, $ S\in [S_+,-\infty)$ on which
$(\tau,S)$ satisfy both \eqref{hugrel} and the Lax condition.
Moreover, we shall show that this curve contains all of the points satisfying 
\eqref{hugrel} (without the Lax condition) on the region $[\tau]<0$.
From these two latter properties, we find, reversing the roles of
$(\tau_-,S_-)$ and $(\tau_+,S_+)$, that any points on
$[\tau]>0$ satisfying \eqref{hugrel} consist of reversed Lax shocks, so
are not in the backward Hugoniot set, and so we can conclude that the
curve we have constructed is the entire backward Hugoniot set.


Applying \eqref{g5}--\eqref{g6}, we find that
$$
[p]= 
(p(\tau_+,S_+)- p(\tau_+,S_-))
+
(p(\tau_+,S_-)- p(\tau_-,S_-))
\geq
p_{\tau}(\tau_-,S_-) [\tau]
$$
whenever $[S]>0$,
whence $-[p]/[\tau]>p_\tau (\tau_-,S_-)$ whenever $[\tau]<0$, or
$\sigma^2>c_-^2$, verifying the first Lax condition 
so long as $[S]\geq 0$, $[\tau]<0$.

(i) Differentiating \eqref{hugrel}, we have
$ 
H_{S_-}=-e_{S_-}+\frac12 p_{S_-}[\tau]<0 
$
by \eqref{g3} and \eqref{g6}, whence, by the Implicit Function
Theorem, the local backward Hugoniot curve extends as a graph
$S(\tau)$ so long as $S$ remains finite.
Continuing, we compute that
$$
\begin{aligned}
H_{\tau_-}&=
-e_{\tau_-} + \frac12( p_{\tau_-}[\tau] - (p_++p_-))
= \frac12( p_{\tau_-}[\tau] - [p])
=\frac12[\tau](\sigma^2-c_-^2)<0
\end{aligned}
$$
so long as $[S]\ge 0$,
by the (already-verified) first Lax condition, $c_-^2<\sigma^2$,
and thus
$
\frac{dS_-}{d\tau_-}=-\frac{H_{\tau_-}}{H_{S_-}}<0
$
so long as $[S]\geq 0$, from which we find that $\sgn[S]=-\sgn[\tau]>0$
persists along the entire curve, i.e., so long as $S$ remains finite.
Rewriting the Hugoniot relation as
$e_0-e= (1/2)(p+p_0)[\tau]$, we see readily that, assuming \eqref{g1}--\eqref{g2},
\be\label{taumax2}
\tau_- \le \tau_{max}:= 2e_0/p_0 + \tau_0.
\ee
whence $\tau_+ \le \tau_- \le \tau_{max}$, with $S\to -\infty$ as
$\tau \to \tau_{\rm max}$.
Thus, the backward Hugoniot consists of the points $(\tau, S(\tau))$
for $\tau\in [\tau_+,\tau_{\rm max})$
and this curve terminates at $(\tau,S) \to (\tau_{\rm max}, -\infty)$
with $e,p\to 0$, by \eqref{h1}--\eqref{h2}.
Moreover, the global inequality $H_{S_-}<0$ implies
that, as claimed, there are no other solutions of \eqref{hugrel} for $[\tau]<0$.

(ii) We have already shown that $dS/d\tau<0$,
whence we obtain also $dp/d\tau= p_S (dS/d\tau) + p_\tau<0$,
by $p_\tau<0$ and $p_s>0$. This verifies monotonicity of $\tau$, $p$, and $S$.
From $[v]=-\sqrt{-[p][\tau]}$ (computed from \eqref{rh}), we obtain
also the stated monotonicity of $v$.
Next, we compute that
$
\frac{d}{dS_-} \frac{[p]}{[\tau]}=
\frac{ [\tau](d/dS_-)[p]- [p](d/dS_-)[\tau] }
{[\tau]^2}
$
has the sign of
$$
 [\tau](d/dS_-)[p]- [p](d/dS_-)[\tau] =
([p]-p_\tau[\tau])(d\tau/dS_-) - [\tau]p_{S_-}
=[\tau] \big( (c_-^2-\sigma^2)(d\tau/dS) -p_S\big),
$$
or, substituting $(d\tau/dS)=-H_{S_-}/H_{\tau_{-}}$,
$
[\tau]\Big( \frac{(c_-^2-\sigma^2)}{(c_-^2-\sigma^2)}\frac {(-e_S+\frac12p_S[\tau])}
     {\frac12[\tau]} - p_S\Big)=
-2e_S<0.
$
Thus (since $S$ is decreasing) 
$\sigma^2=-[p]/[\tau]$ is monotone decreasing
along the backward Hugoniot.
Noting that the initial value of $\sigma^2$ at $(\tau_+,S_+)$ is $c_+^2$,
we have that $c_+^2>\delta^2$,
verifying the second Lax condition along the backward Hugoniot and
completing the proof.

\end{proof}

\br[Monotonicity of $S$, $\sigma$ along forward Hugoniot]\label{weylrmk}
Symmetric computations show that $H_{\tau_+}=\frac12(\sigma^2-c_+^2)[\tau]>0$
for $[\tau]<0$, where we are here using the fact already shown that the
Lax conditions hold whenever $[\tau]<0$, which yields that $S$ parametrizes
the forward Hugoniot curve.
Similarly, $(d/dS_+)\frac{[p]}{[\tau]}=-2e_{S_+}<0$, so that $S$
and $\sigma^2$ are increasing along the forward Hugoniot curve,
recovering the results of Weyl \cite{W}
under assumptions \eqref{g1}--\eqref{g6}.
\er

\br[Monotonicity of $\tau$ along forward Hugoniot]\label{vasprep}
Monotonicity of $\tau$ does not always hold in the forward direction.
Differentiating \eqref{hugrel} with respect to $+$ variables, we have
$$
\begin{aligned}
H_{\tau_+}&=
e_{\tau_+} + \frac12( p_{\tau_+}[\tau] + (p_++p_-))
= \frac12( p_{\tau_+}[\tau] - [p])
=\frac12[\tau](c_+^2-\sigma^2)< 0
\end{aligned}
$$
by the second Lax condition, $\sigma^2<c_+^2$, and
$ H_{S_+}=e_{S_+}+\frac12 p_{S_+}[\tau] $, hence $\tau$ is
a monotone decreasing function of $S$ along the forward Hugoniot curve
precisely so long as $-H_{S_+}/H_{\tau_+}<0$, or $p_S<-2T/[\tau]$.
Recall that this is Erpenbeck's condition, mentioned in Remark \ref{monrmk},
which is already sufficient to imply \eqref{dcond3} and stability.
Likewise, $v$ and $p$ need not be monotonically related.
\er

\subsection{General pressure laws}\label{s:entropy}
We now consider cases of general pressure laws defined in the
absence of a complete equation of state $\bar e$.
We look at two types of such laws.
\subsubsection{}
We begin with the 
question of when a general pressure law
$$p=\hat p(\tau,e)$$ determines
a full equation of state $e=\bar e(\tau,S)$, or 
integral of the fundamental thermodynamic relation
$$
de=T dS-p d\tau,
$$
with positive temperature $T=\bar e_S$, or, equivalently,
when we may construct an entropy function $S=\hat S(\tau,e)$,
with $S_e=\bar e_s^{-1}>0$, $S(\tau, \cdot):\RM^+\to \RM$,
satisfying the transport equation
\be \label{Seq}
S_\tau - \hat p(\tau,e) S_e=0.
\ee

As noted in Remark \ref{prmk}, for $\hat p\in C^1$,
\eqref{Seq} may always be solved locally
to particular values $(\tau_0,e_0)$, by prescribing an arbitrary
monotone increasing profile $\hat S(\tau_0, \cdot)$ at $\tau_0$, taking
$\RM^+\to \RM$, and
evolving in forward and backward $\tau$ by the method of characteristics,
so long as characteristics continue to cover the full half-line
$e\in \bar \RM^+$ (they do not intersect, by uniqueness of solutions
of the characteristic equation $de/d\tau=-p(\tau,e)$), 
monotonicity being preserved by construction.
Equivalently, we may prescribe an arbitrary positive temperature
profile $ \hat T(\tau_0, \cdot)$ at $\tau_0$, such that $\int (1/T)de$
diverges at both $e\to 0^+$ and $e\to +\infty$.

The following result, the key observation of this section, states that
under the Weyl condition we may obtain also {\it global existence in 
the direction of decreasing $\tau$.}

\bl\label{Slem}
For $\hat p\in C^1(\RM^+\times \bar{\RM^+})$ 
satisfying the Weyl condition $\hat p_e>0$ and $\hat p(\tau, \cdot)\equiv 0$,
and any temperature profile $\hat T(\tau_0, \cdot)>0$
such that $\int (1/T)de$ diverges at both $e\to 0^+$ and $e\to +\infty$,
there exists a unique entropy function $S=\hat S(\tau,e)$
satisfying \eqref{Seq} on $0<\tau\leq \tau_0$, $e\in \RM^+$, with
$\hat T(\tau,e):=(\hat S_e)^{-1}>0$.
Equivaelently, there is a unique equation of state $\hat e(\tau,S)$
defined for $0<\tau\leq \tau_0$ and $S\in \RM$ with $\hat T=(\hat S_e)^{-1}
=\hat e_s$ and $\hat p=-\bar e_\tau$.
\el

\begin{proof}
Under $\hat p_e>0$, we find that characteristics $e(\tau)$,
$de/d\tau=-\hat p(\tau,e)$ diverge in the backward (decreasing $\tau$)
direction, hence continue to cover a half-line $[a(\tau),+\infty)$.
Moreover, the property $\hat p(\tau, \cdot)\equiv 0$ ensures that
the characteristic emanating from $(\tau_0,0)$ remains identically
zero, so that, by continuous dependence of characteristics,
$a(\tau)\equiv 0$ and
the property $\lim_{e\to 0^+}\hat S(\tau,e)=-\infty$ is preserved.
Meanwhile, the property $\lim_{e\to +\infty}\hat S(\tau,e)=+\infty$ 
is preserved by monotonicity of characteristics with respect to $e$
and by the fact that $\hat S$ by construction is constant on characteristics.
\end{proof}

\bc\label{Jprop}
For a $C^2$ pressure law $\hat p$,
assumptions \eqref{j1}--\eqref{j4} imply \eqref{p1}--\eqref{p5}.
\ec

\begin{proof}
Applying Lemma \ref{Slem} for any fixed right state $\tau_+$, 
we obtain an equation of state $\bar e$ satisfying
\eqref{g1}--\eqref{g6} an \eqref{h1}--\eqref{h2} 
for $0<\tau\leq \tau_{\rm max}$, where $\tau_{\rm max}$ is 
defined as in \eqref{taumax2},
which generates the pressure law $\hat p$ through the relation 
$ p:=-\bar e_\tau$.
Noting that only $t_+\leq \tau\leq \tau_{\rm max}$ is relevant for the backward Hugoniot
curve through $U_+$, we find
by Proposition \ref{smithprop} that this implies \eqref{p1}--\eqref{p5}.
\end{proof}

\br\label{existrmk}
Under the slightly extended set of assumptions \eqref{g1}--\eqref{g6},
\eqref{h1}--\eqref{h4}, Smith \cite{S} obtains
not only uniqueness but also existence.
It would be interesting to determine what conditions beyond
\eqref{j1}--\eqref{j4} would be necessary to obtain existence
for a general pressure law $\hat p$.
\er

\subsubsection{}  
Secondly, we briefly consider the situation of a pressure relation
$$
p=\tilde p(\tau, T)
$$ 
defined in terms of temperature $T$ instead of
internal energy $e$, for example an ideal gas law
$$
\hbox{\rm
$p=\tilde p(\tau,T):=RT/\tau$, $R>0$ constant.
}
$$
Using a result of \cite{CF} 
on the form of possible extensions of the pressure law to
a complete equation of state,
(see (6.77)-(6.80) of \cite{CF}),
Smith \cite{S} shows that an ideal gas
necessarily satisfies \eqref{strongintro}, hence the form of the 
pressure law imposes both uniqueness and shock stability,
independent of the particular extension to a full equation of
state $e=\bar e(\tau,S)$.

We note that similar 
computations may be carried out 
for an arbitrary relation $p=\tilde p(\tau, T)$;
namely, conversion to a general pressure law
$p=\hat p(\tau,e)$ amounts to the construction of an energy function
$$
e=\tilde e(\tau,T), \quad \tilde e_T>0,
$$
from which we obtain $T=\hat T(\tau,e)$ by inversion with respect to $e$,
and thus $\hat p(\tau,e):=\tilde p(\tau, \hat T(\tau,e))$.
This yields the relations
\be\label{tilderel}
\hat p_\tau=\tilde p_\tau + \tilde p_T \hat T_\tau,
\quad
\hat p_e=\tilde p_T \hat T_e,
\quad
\hat T_\tau=-\frac{\tilde e_\tau}{\tilde e_T},
\quad
\hat T_e=\frac{1}{\tilde e_T}.
\ee
Recalling the thermodynamic relation $\hat S_\tau -p\hat S_e=0$
and $\hat S_e=\hat T^{-1}$, we obtain
$\hat T_\tau - p\hat T_e + T\hat p_e=0$, hence, substituting from
\eqref{tilderel},
\be\label{Trel}
\tilde e_\tau= T\tilde p_T- p,
\ee
hence $\tilde e_{T \tau }= \tilde e_{\tau T}= T\tilde p_{TT}$, 
and we can solve
globally for a positive temperature $e$ and derivative
 $e_T$ for $\tau\geq \tau_0$ provided 
$\tilde p_{TT}\geq 0$, and $\tilde p(\tau,0)\equiv 0$, which together
imply also $ T\tilde p_T-p\geq 0$.
For further discussion, particularly of relation \eqref{Trel},
see the discussion of Helmholtz energy in Appendix \ref{s:helmholtz}.

In the ideal gas case, we have equality, $\tilde p_{TT}=
T\tilde p_T-p=0$, so find that 
$
\tilde e(\tau, T)=g(T)
$
depends on $T$ alone, where $g$ is any monotone increasing function
one-to-one on $\RM^+$ to $\RM^+$.
We thus have $\hat T_\tau=0$, and, from \eqref{tilderel}, 
$\hat p_\tau=\tilde p_\tau$, so that \eqref{strongintro} is
equivalent to $\tilde p_\tau<0$, independent of the choice of $g$.
Moreover, we have by direct computation $\tilde p_\tau=-RT/\tau^2<0$,
recovering somewhat more simply the result of Smith.

However, note that any admissible temperature function $\tilde e(\tau,T)$
may be modified by the addition of any nondecreasing function
$h(T)$ with $h(0)=0$ to yield another admissible temperature function
$\tilde e_1(\tau,T):=\tilde e(\tau,T)+h(T)$.
Moreover, {\it except in the special situation $\tilde e_\tau\equiv 0$ just
considered, or when $\tilde p_T\equiv 0$}, such modification will result in a change in $\hat T_\tau$ and $\hat p_\tau$, so that satisfaction of
our conditions \eqref{weakintro}, \eqref{strongintro},
 \eqref{mediumstab}, and \eqref{mediumintro}
may be affected by the choice of $\tilde e$, in contrast with the situation
of pressure laws of ``primary'' type $p=\hat p(\tau,e)$.

\bl\label{tauTconditions}
In terms of a pressure law $\tilde p(\tau,T)$, 
hyperbolicity becomes
$ \tilde  p_\tau < \frac{T \tilde p_T^2}{\tilde e_T},$ or
\be\label{Thyp}
\hat p_\tau=\tilde p_\tau - \tilde p_T \tilde e_\tau/\tilde e_T < p\tilde p_T/\tilde e_T=p\hat p_e,
\ee
while \eqref{strongintro} and \eqref{weakintro} become, respectively
\be\label{Tweakstrong}
\hbox{\rm
$\hat p_\tau=\tilde p_\tau - \tilde p_T \tilde e_\tau/\tilde e_T<0$
(strong)
and
$\hat p_\tau=\tilde p_\tau - \tilde p_T \tilde e_\tau/\tilde e_T < 
p\tilde p_T/2\tilde e_T=p\hat p_e/2$
(weak).
}
\ee
\el

\begin{proof}
Immediate, from \eqref{tilderel}, 
\end{proof}

\br
These findings are useful for 
the general van der Waals gas law
$(p+a/\tau^2)(\tau-b)=RT$,
$a,b,R>0,$ 
where
$ \tilde p(\tau,T)= \frac{RT}{\tau-b} -\frac{a}{\tau^2}. $
This will be pursued elsewhere.
\er

\subsection{Link to the Riemann Problem}
The previous subsections have revealed a certain 
parallel
between the transition from 
stability to instability for individual shocks on the hand and 
that from uniqueness to non-uniqueness for the Riemann problem on the other. 
Clearly, these two transitions are not generally equivalent to each other: 
there are equations of state for which all shocks have positive signed Lopatinski 
determinant and yet for certain data, the Riemann problem has more than one solution.
On the other hand, according to Theorem 1.1, the mere existence of shocks with negative signed 
Lopatinski determinant
readily implies the existence of Riemann data admitting several solutions.  
This section serves to make this connection transparent.

As pointed out already by C. Gardner \cite{G} in the context of gas dynamics,
vanishing of the Lopatinski determinant has a natural interpretation in terms 
of bifurcation of Riemann solutions.
\par\medskip

{\bf Proof of Theorem 1.3.} 
For single-shock data, the Lopatinski determinant coincides with the determinant of
the local Jacobian of Lax's wave map \cite{L}. 
The assertion now follows from 
%
the following general fact that we briefly prove as 
we know no reference for it. 

\bl
For any $C^1$ mapping $F:\R^n\to\R^n$, local uniqueness prevents
det$(DF)$ from undergoing a local strict change of sign. 
\el

\begin{proof}
Local uniqueness at a point $p$ implies the existence of an open ball $B$ around $p$
such that $F(B)\cap F(\partial B)=\emptyset$. With $d$ denoting the Brouwer degree,
the function $B\to\R, q\mapsto d(F,B,F(q))$ is thus a constant.
As $d(F,B,F(q))=\text{sgn det}DF(q)$ whenever $\text{det}DF(q)\neq 0$, 
$B$ cannot simultaneously contain points $q^\pm$ with $\pm\text{det}DF(q^\pm)> 0$.
\end{proof}

\br[Stability as selection principle]
The above-described bifurcation analysis bears on the question posed
in \cite{MP} whether shock stability could serve
in the situation of multiple Riemann solutions as a selection principle.
For, generically, the bifurcation
local to the single-shock solution
should be of fold type, with $\Delta$
taking different signs (corresponding to opposite orientations) on either
side of the fold;
see section 6.2, \cite{Z1}, for an analysis of the general $2\times 2$ case,
and \cite{B-G} for an extension to the $3\times 3$ gas dynamical case. 
Assuming that additional small-amplitude waves in other family are stable,
as is expected to be the case, we find that sign of $\Delta$, and thus
stability of the large component shock near $(U_-,U_+,\sigma)$,
determines stability of the associated Riemann patterns, so that, at the
level of Riemann patterns, this is a {\it transcritical bifurcation}
featuring exchange of stability from one branch to another, and indeed
serves as a selection principle.  

%
%
\er


\section{Gas-dynamical examples and counterexamples}

\subsection{Global counterexample}\label{s:global} 
Recall the equation of state 
\be\label{ceg2}
\bar e(\tau,S)= \frac{e^{S}}{\tau}+ C^2 e^{ S/C^2 -\tau/C},
\quad C>\! >1
\ee
asserted to exhibit instability in Theorem A.

\subsubsection{}{\bf Proof of Theorem 1.2.}
The function $\bar e$ is evidently convex. 
Computing, we have 
$$
T=\bar e_S= e^S/\tau + 
e^{ S/C^2 -\tau/C}>0,
$$
$$
p=-\bar e_\tau= e^S/\tau^2 +C 
e^{ S/C^2 -\tau/C}>0,
$$
$$
p_S=-\bar e_{S\tau}= e^S/\tau^2 +
C^{-1}
e^{ S/C^2 -\tau/C}<0,
$$
$$
p_\tau=-\bar e_{\tau \tau}= -2 e^S/\tau^3 -
e^{ S/C^2 -\tau/C}<0,
$$
$$
p_{\tau \tau}=-\bar e_{\tau \tau \tau}= 6 e^S/\tau^4 + 
C^{-1}
e^{ S/C^2 -\tau/C}>0,
$$
whence this model satisfies the basic properties \eqref{g1}--\eqref{g6}
assumed by Smith: 
positivity of $e$, $T$, $p$, genuine nonlinearity, and the Weyl condition.

Continuing, it is readily seen that 
$$ \lim_{S\to -\infty} \bar e(S,\tau)=0,
\quad
 \lim_{S\to +\infty} \bar e(S,\tau)=\infty,
\;\hbox{\rm and }\;
\lim_{\tau\to 0^+} \bar e(S,\tau)=\infty,
\quad
\lim_{\tau \to +\infty} \bar e(S,\tau)=0,$$
validating Smith's asymptotic hypotheses \eqref{h1} and \eqref{h2}.
Boundedness of $p$ as $\tau\to 0^+$ clearly implies
that $S\to -\infty$, or else the first term $e^S/\tau^2$ would blow up.
Noting that for $e^S/\tau^2\le C$, $e^S/\tau\le C\tau \to 0$
as $\tau\to 0^+$, we thus obtain $e\to 0$,
verifying Smith's hypothesis \eqref{h3}:
$ \lim_{\tau\to 0^+} \hat e(p,\tau)=0$.
Likewise, \eqref{h4} is readily verified: 
$\lim_{\tau \to +\infty} \bar e(s,\tau)=0$.
Thus, $\bar e $ satisfies all of Smith's hypotheses
\eqref{g1}--\eqref{g6} and \eqref{h1}--\eqref{h4}.

By Proposition \ref{smithprop} and Remark \ref{vasprep},
we thus obtain monotonicity
of $S$ and $\sigma$ along forward and backward Hugoniot
curves.
Evaluating at $\tau_+=1$, $S_+=0$, and using the assumed 
asymptotics $C>\! >1$, we obtain
$$
\frac{-\bar e_{S\tau}}{\bar e_{\tau \tau}}-\frac{2\bar e_S}{-e_\tau }
= \frac{1+O(C^{-1})}{3+O(C^{-1})}-\frac{2+O(C^{-1})}{C +O(1)} 
\sim \frac{1}{3}>0,
$$
giving failure of Smith's condition \eqref{weakintro}.


\subsubsection{}{\bf Numerical illustration}\label{s:numill}
In Figures \ref{ceghugf} and \ref{ceghugb} below, 
we illustrate the theory 
presented in Theorem A
by numerical computation, displaying the
forward and backward Hugoniot curves
for an unstable $1$-shock of \eqref{ceg},
with endstates $(\tau_+,S_+)=(1,0)$ and $(\tau_-,S_-)=(12.089548257354499,-6)$,
described as plots of
$\tau$, $p$, $v$, $e$, $\sigma$ against $S$.
Here, the Hugoniot curve is obtained numerically
by solving \eqref{hugrel} for $\tau$
as a function of $S$ using the bisection method.
The transition to instability is computed along the
backward Hugoniot using \eqref{mediumstab}.

\begin{figure}[htbp]
 \begin{center}
$
\begin{array}{lcr}
(a) \includegraphics[scale=0.2]{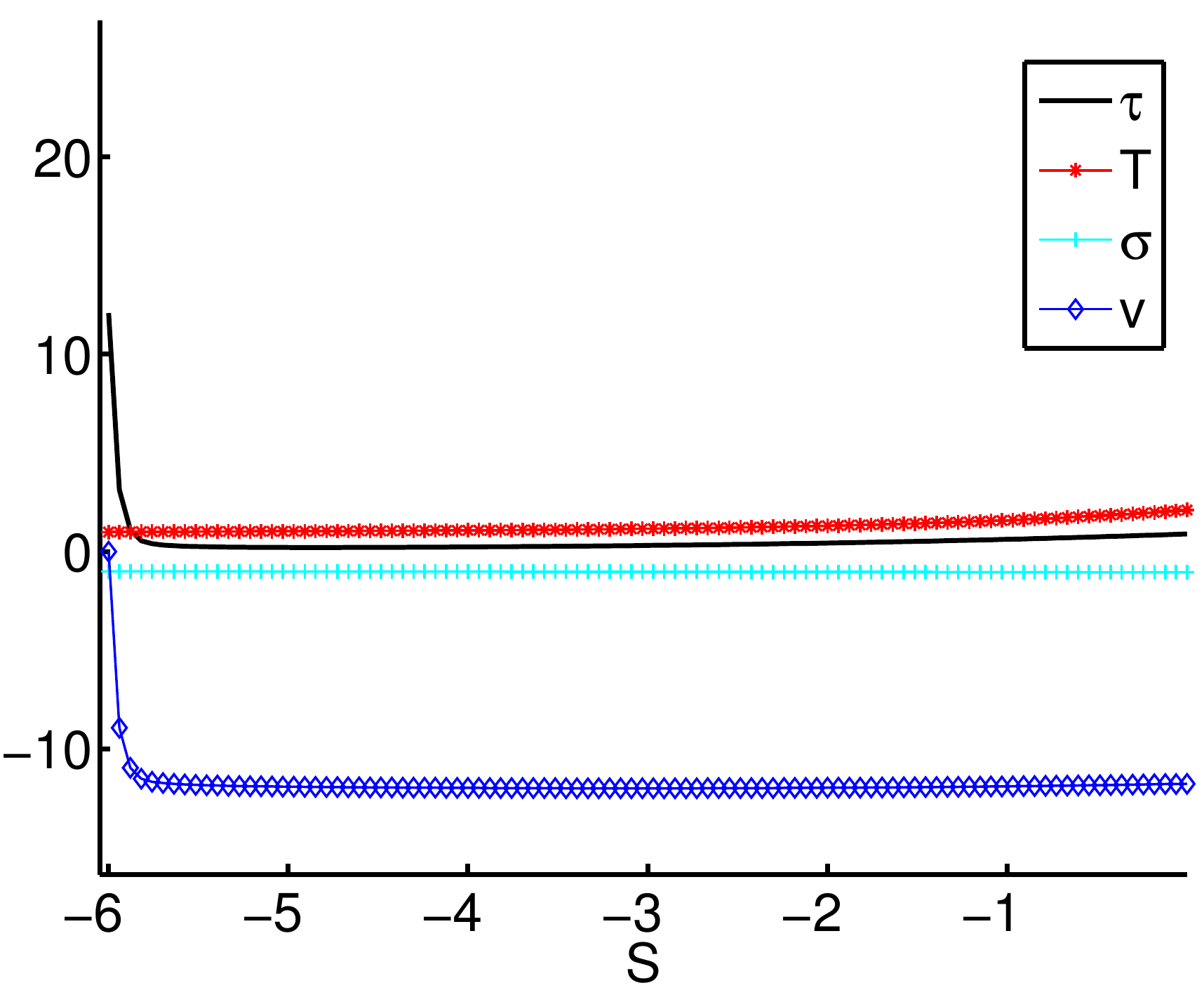}\quad (b) \includegraphics[scale=0.2]{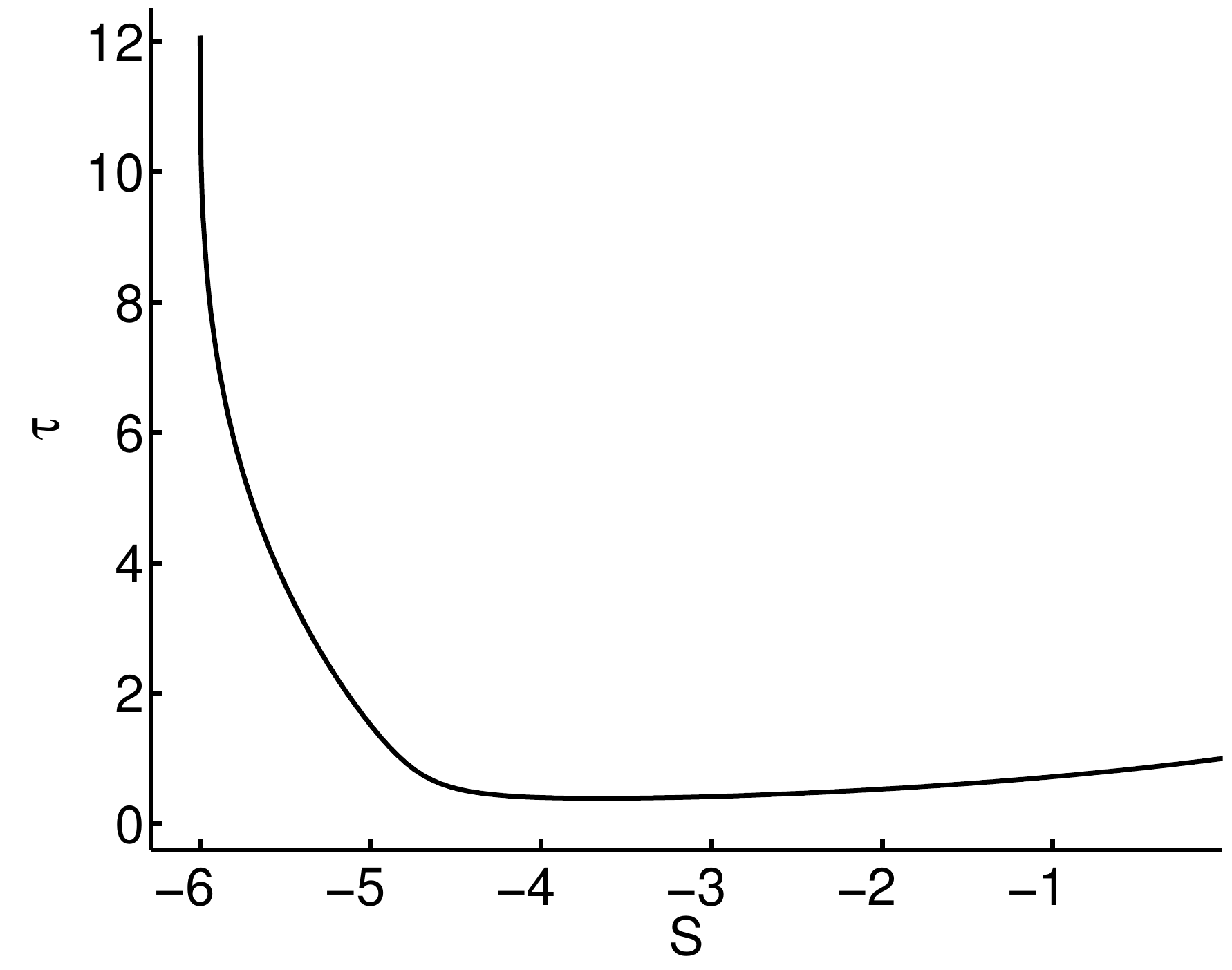}(c) \includegraphics[scale=0.2]{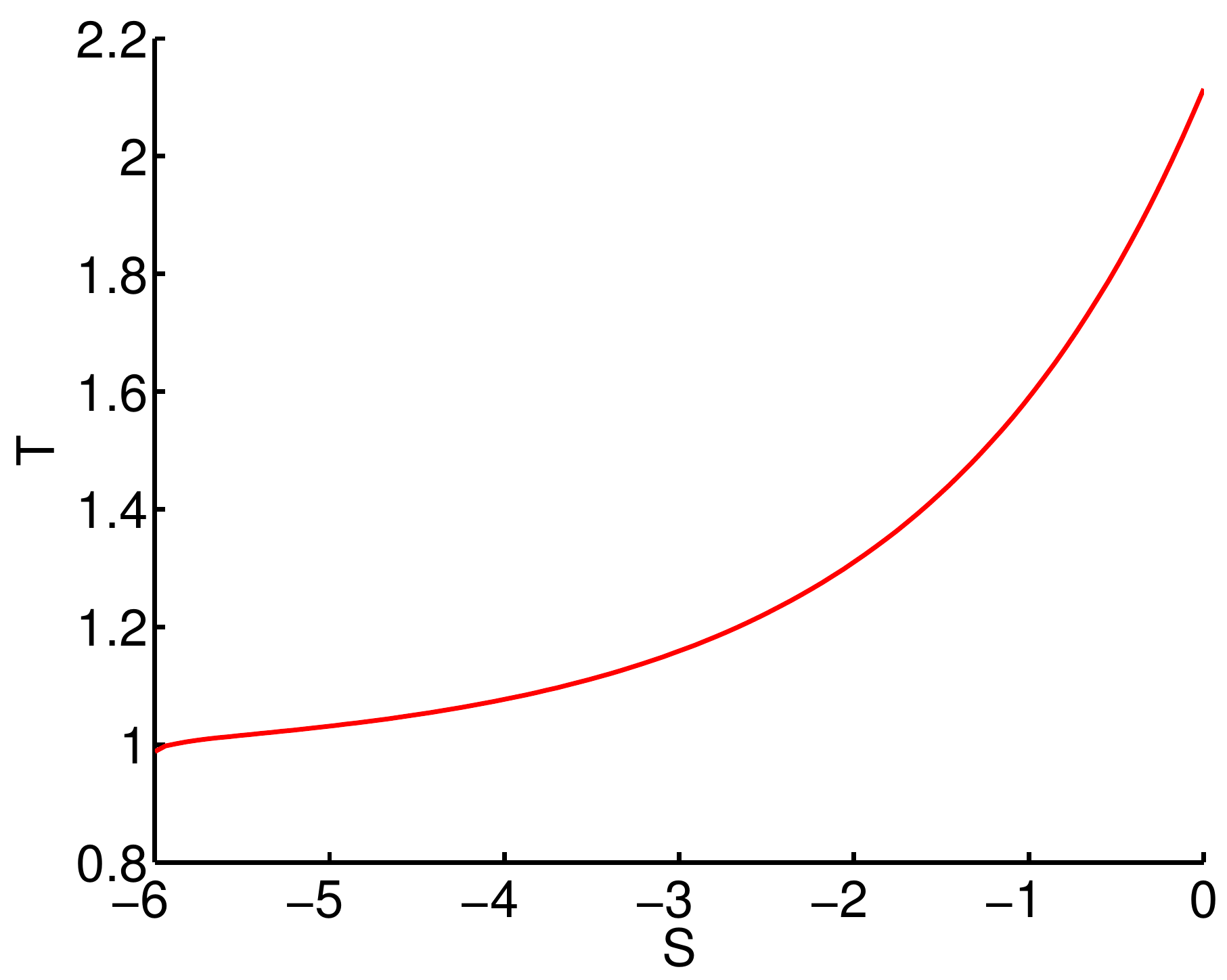}\\
(d) \includegraphics[scale=0.2]{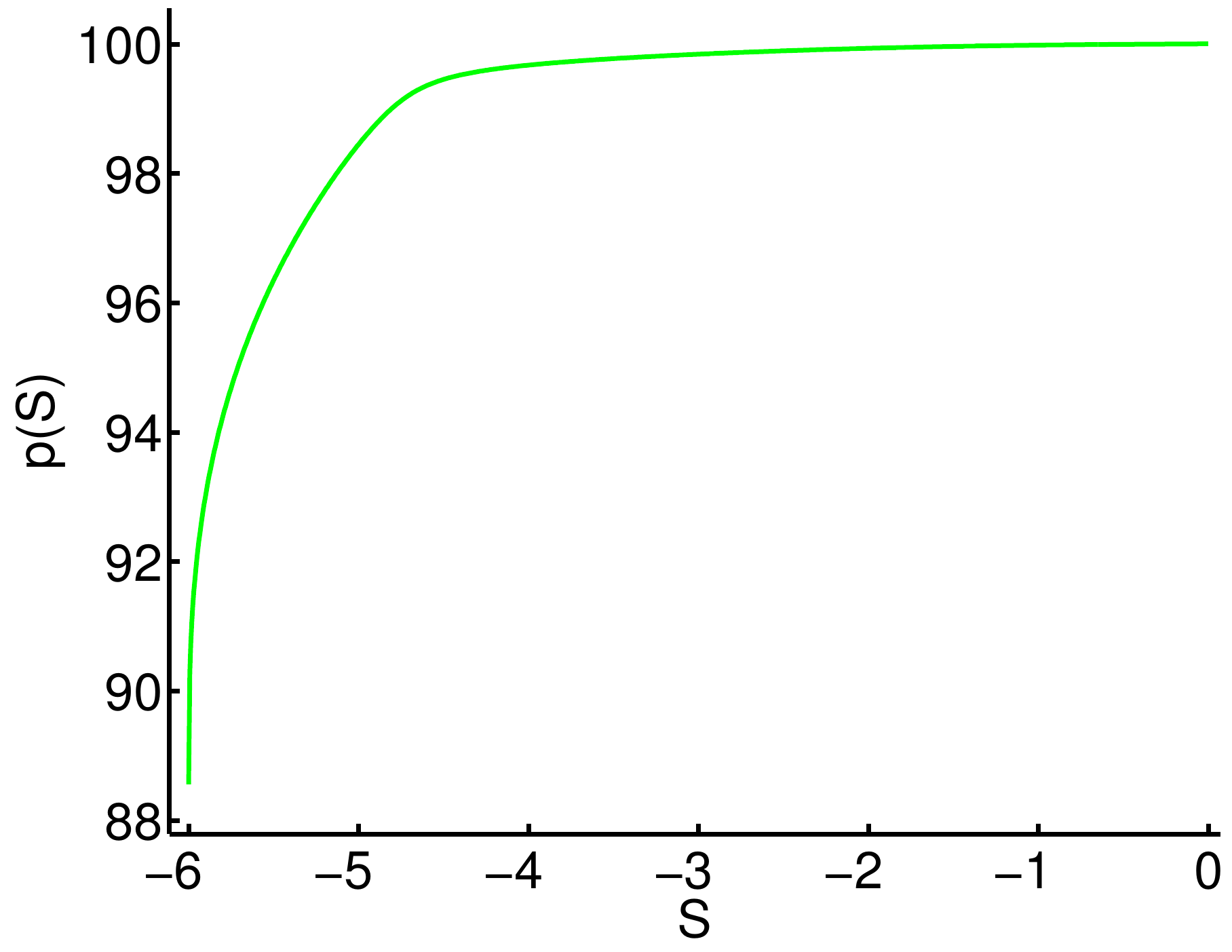}\quad (e) \includegraphics[scale=0.2]{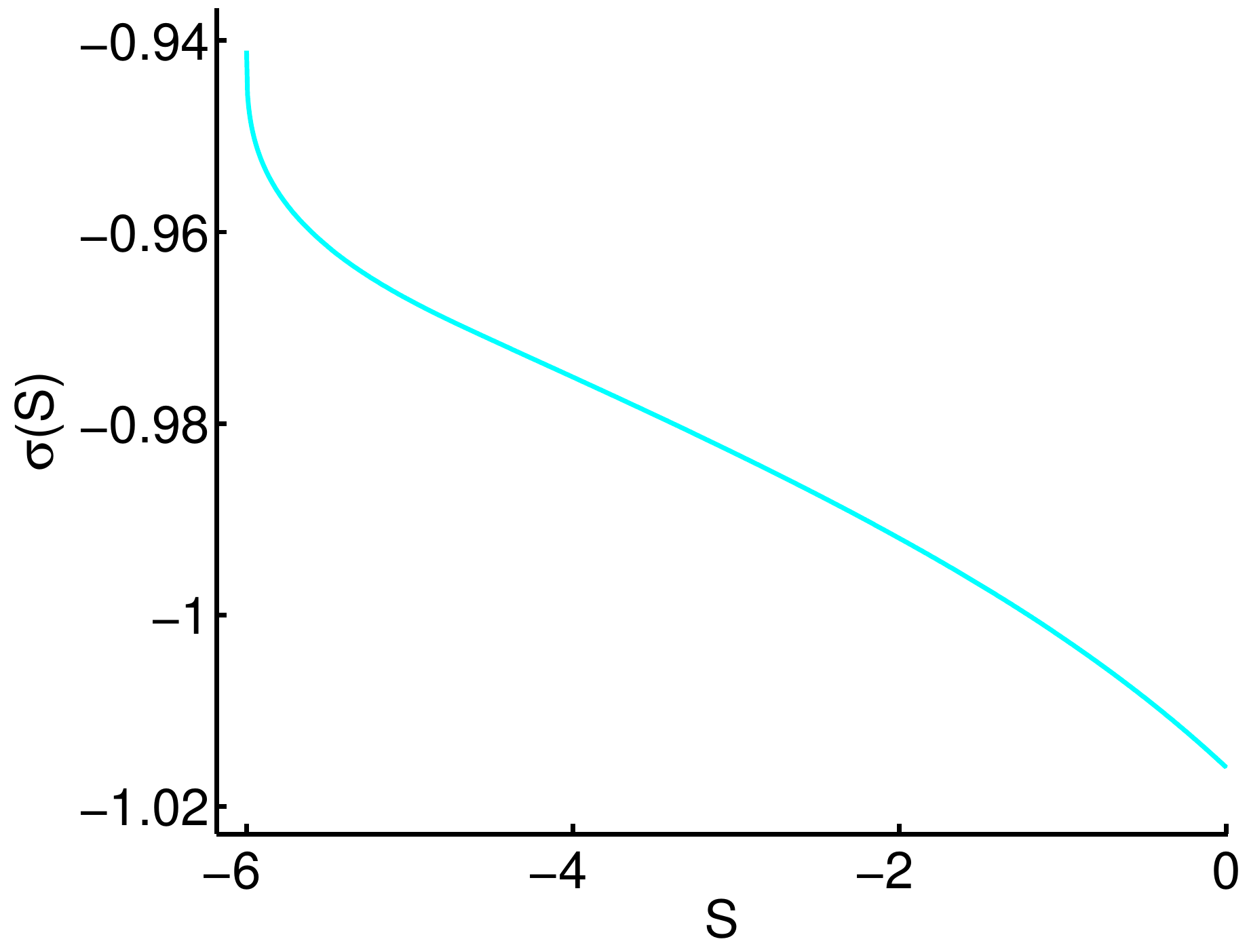}(f) \includegraphics[scale=0.2]{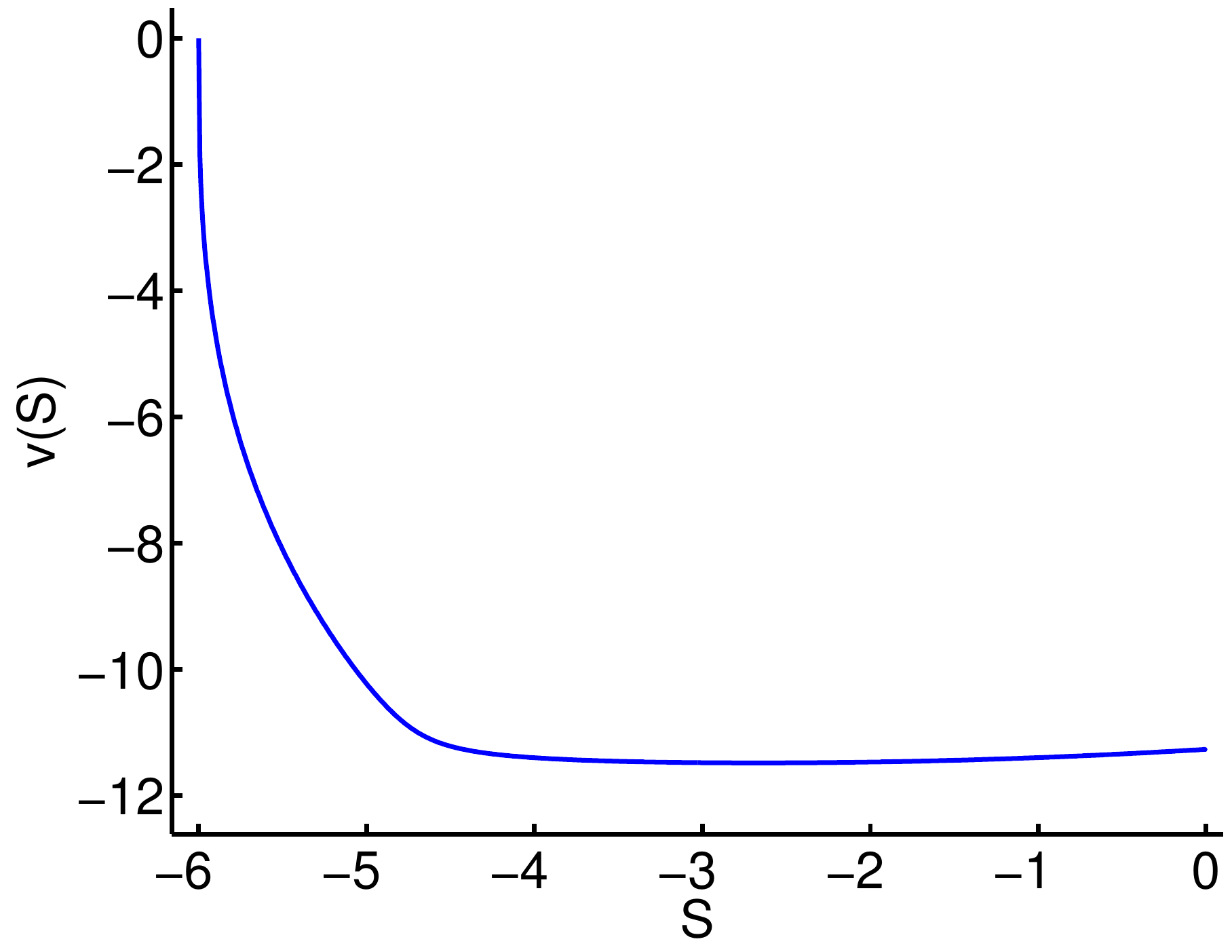}
\end{array}
$
\end{center}
\caption{
(a) The forward $1$-Hugoniot curve 
through $(\tau_-,S_-) \approx (12.09,-6)$
for global model $\bar e(\tau,S)=e^S/\tau + e^{S/C^2-\tau/C}$
of points $(\tau,S)$ to which $(\tau_-,S_-)$ connects by a Lax $1$-shock,
displayed as a graph $(\tau,p,v,e,\sigma)$ over $S$ plotted with respective colors (black, green, blue, red, cyan).  We zoom in to see (b) the Hugoniot curve, (c) $T$ over $S$, (d) $p$ over $S$, (e) $\sigma$ over $S$, and 
(f) $v$ over $S$. Note nonmonotonicity of $\tau$ and $v$ with respect to $S$. }
\label{ceghugf}\end{figure}

\begin{figure}[htbp]
 \begin{center}
$
\begin{array}{lcr}
(a) \includegraphics[scale=0.2]{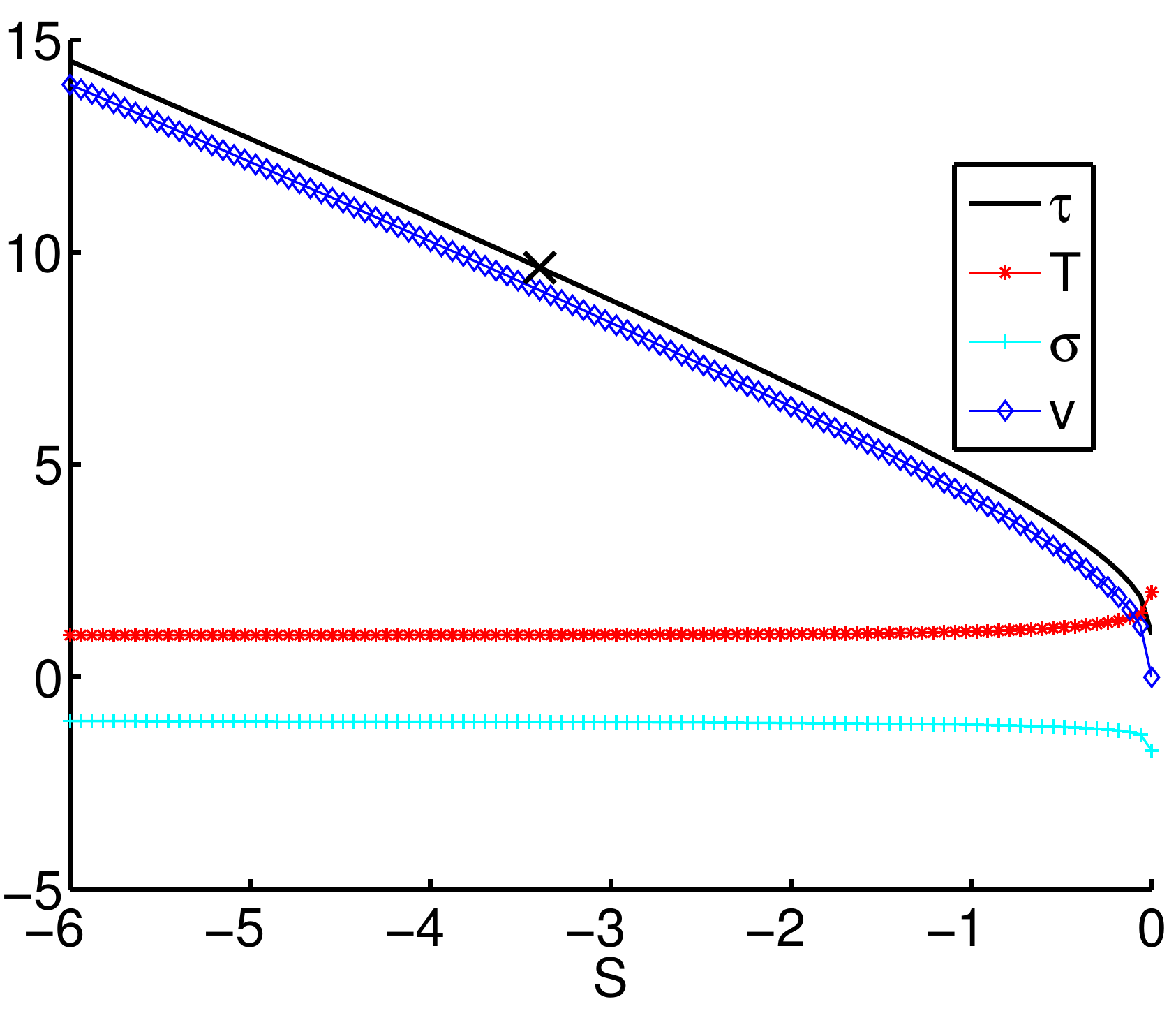}\quad (b) \includegraphics[scale=0.2]{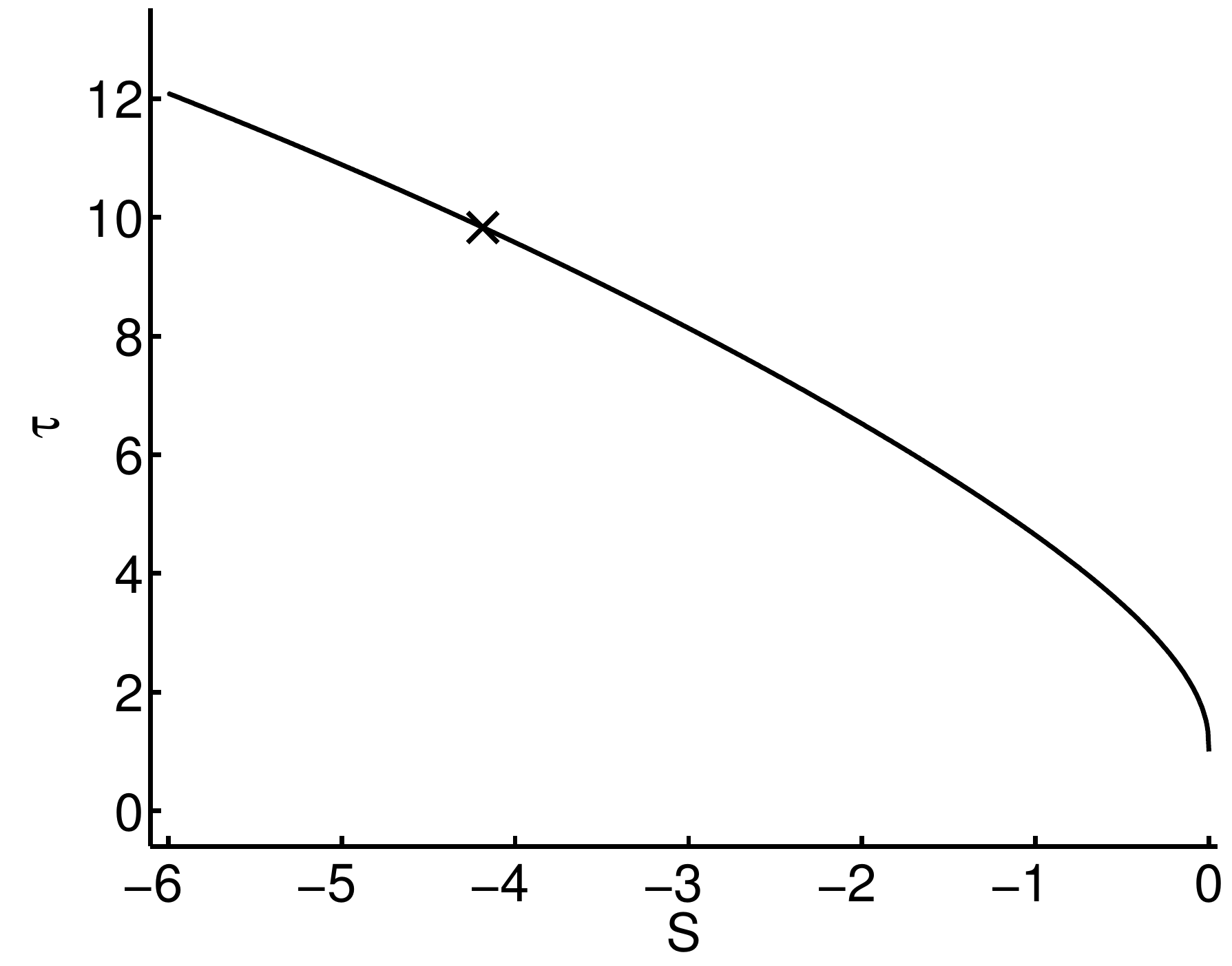}(c) \includegraphics[scale=0.2]{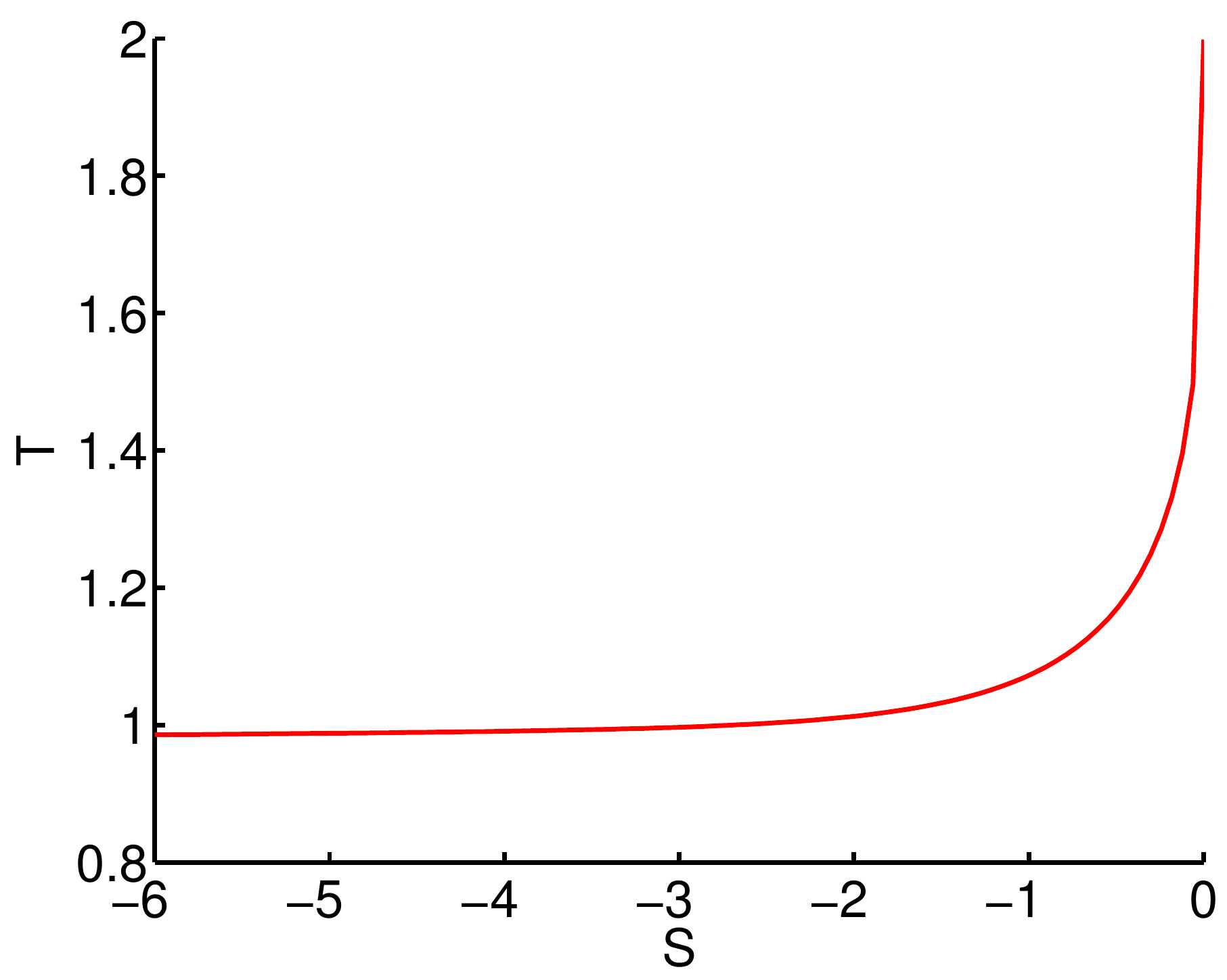}\\
(d) \includegraphics[scale=0.2]{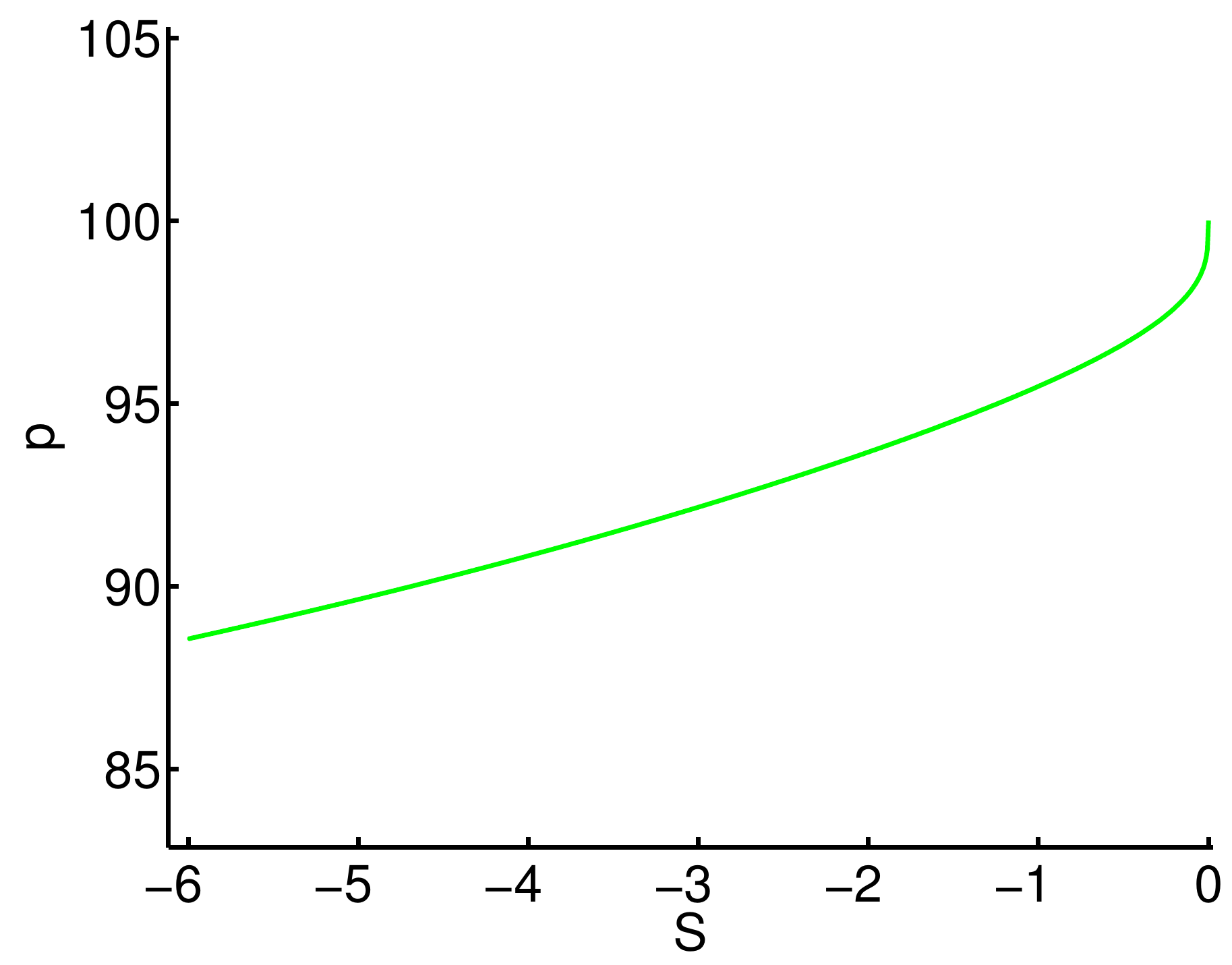}\quad (e) \includegraphics[scale=0.2]{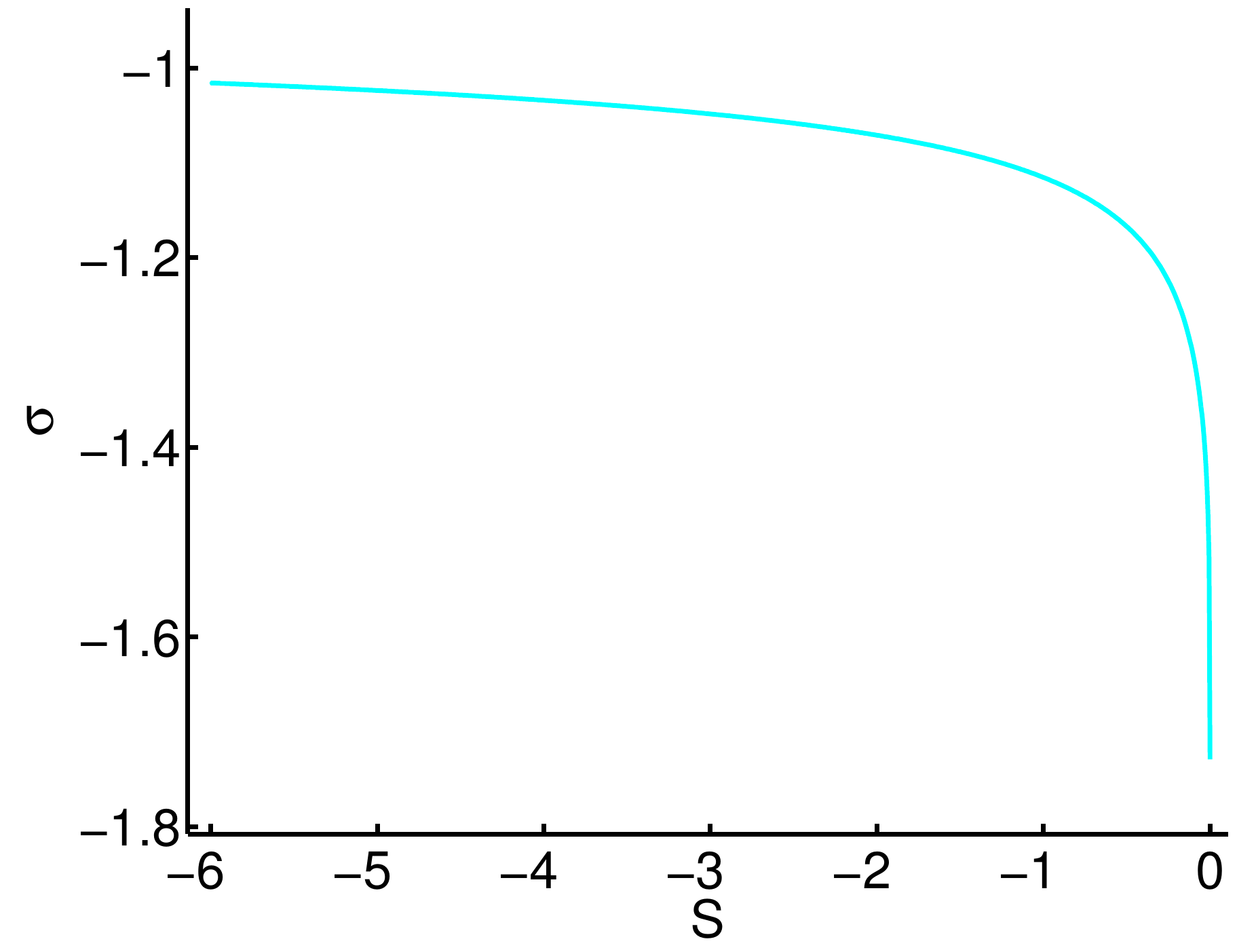}(f) \includegraphics[scale=0.2]{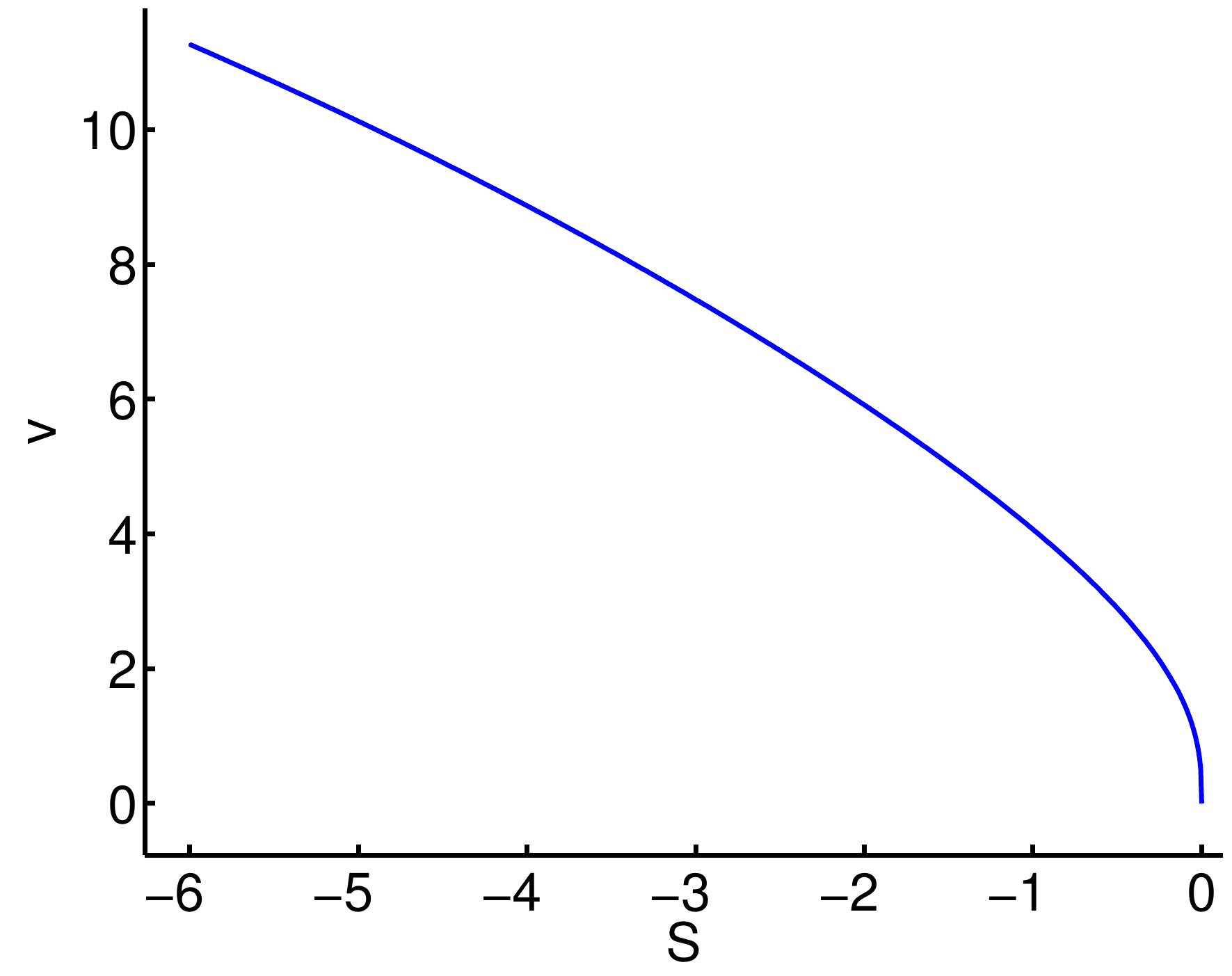}
\end{array}
$
\end{center}
\caption{
The backward $1$-Hugoniot curve through $(\tau_+,S_+)=(1,0)$
for global model $\bar e(\tau,S)=e^S/\tau + e^{S/C^2-\tau/C}$
of points $(\tau,S)$ connecting to $(\tau_+,S_+)$ by a Lax $1$-shock,
displayed as a graph $(\tau,p,v,e,\sigma)$ over $S$ plotted with respective colors (black, green, blue, red, cyan).
We zoom in to see (b) the Hugoniot curve, (c) $T$ over $S$, (d) $p$ over $S$, (e) $\sigma$ over $S$, and (f) $v$ over $S$. The value of $(\tau_-,S_-)$ along the backward Hugoniot curve
at which transition to instability occurs is marked by a black X.
  }
\label{ceghugb}\end{figure}

\begin{figure}[htbp]
 \begin{center}
$
\begin{array}{lcr}
(a) \includegraphics[scale=0.3]{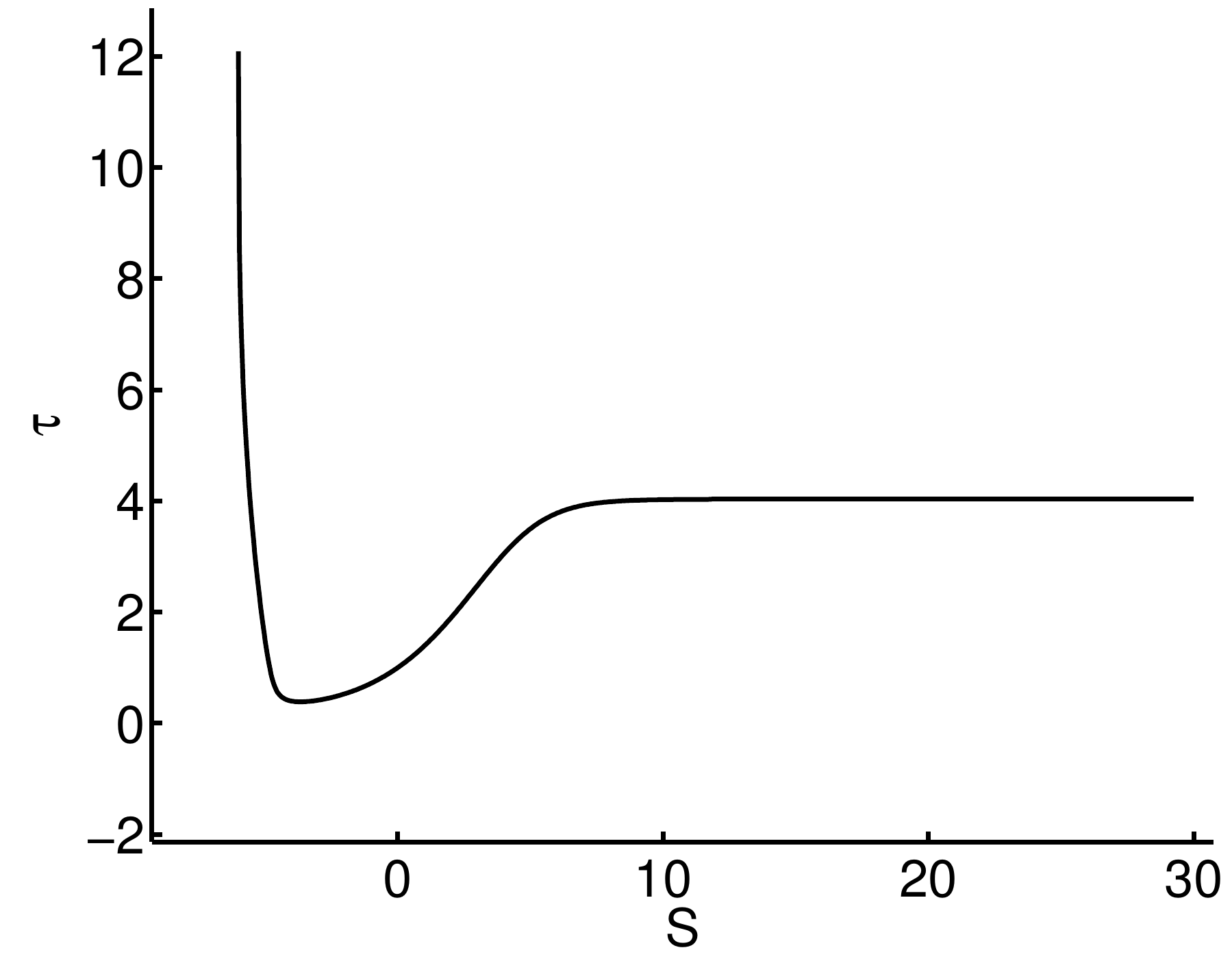}\quad (b) \includegraphics[scale=0.3]{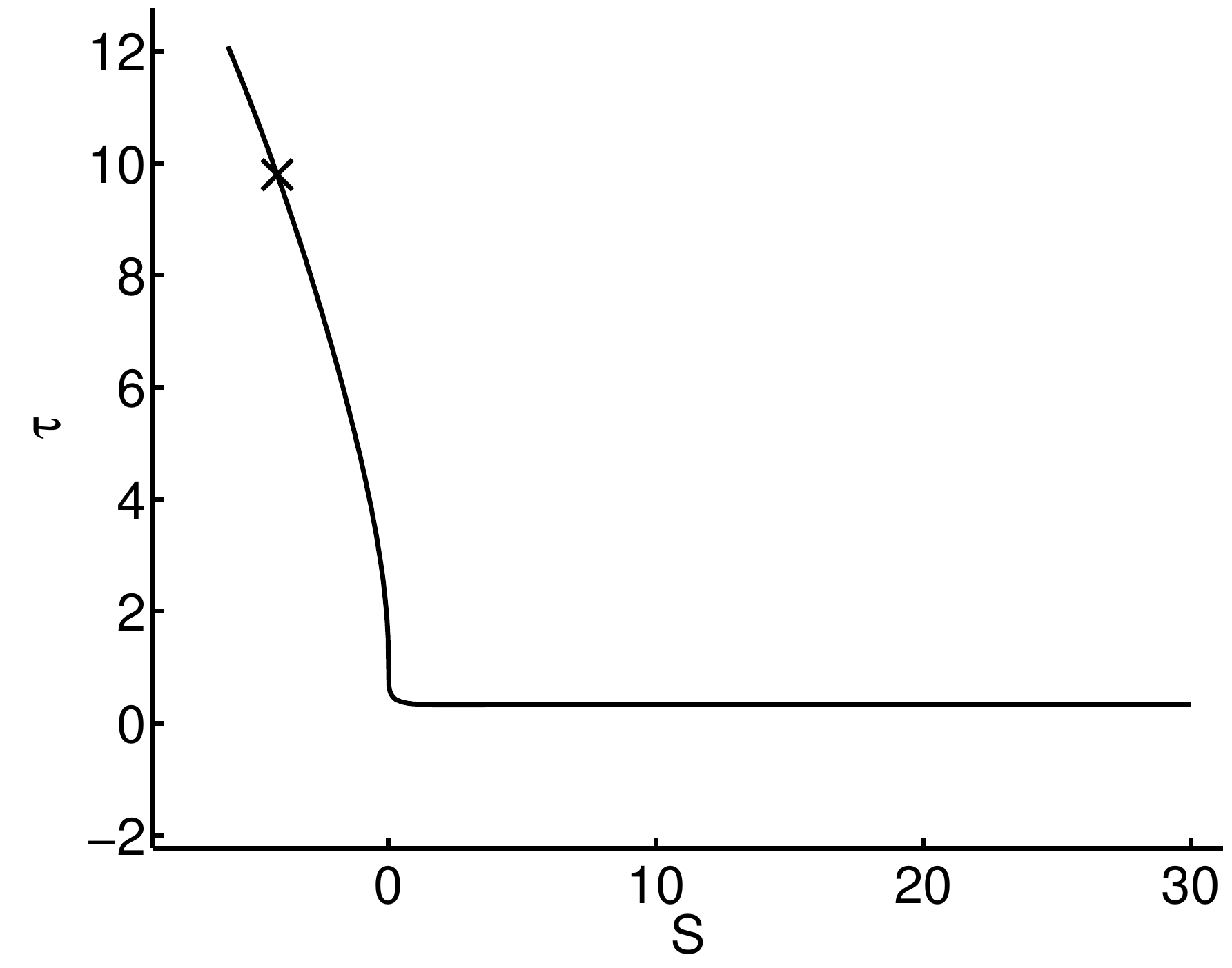}\end{array}
$
\end{center}
\caption{
Extended plots of the (a) forward and (b) backward Hugoniot curves shown in Figures \ref{ceghugf} and \ref{ceghugb}
for global model $\bar e(\tau,S)=e^S/\tau + e^{S/C^2-\tau/C}$.
  }
\label{ceghugfextended}\end{figure}

\subsection{Local counterexample}\label{s:lceg}
Taylor expanding the equation of state \eqref{ceg2}
in the quantitiy $1/C$ about the limiting value
$1/C=0$  as $C\to \infty$, we obtain the simplified,
{\it local  model}
\ba\label{simplified}
\bar e(\tau, S)&= e^S/\tau +  S- C\tau + \tau^2 /2 , \quad C>\! >1,\\
\ea
described in Remark \ref{locexrem}.
It is readily seen that the flow of \eqref{euler} is independent
of the value of $C$, with $C$ only serving to shift pressure
by positive $C$ but not affecting dynamics.
Thus, dropping the assumptions of positive pressure and energy, we
may replace \eqref{simplified} by the canonical model
\be\label{canon}
\bar e(\tau, S)= e^S/\tau +  S + \tau^2 /2 
\ee
for which $C$ no longer appears, capturing the large-$C$ behavior
of \eqref{ceg}.
This is perhaps the simplest model exhibiting a convex entropy,
monotonicity of shock speed and $S$ along the forward Hugoniot,
but also unstable Lax shock waves.

Model \eqref{canon} satisfies \eqref{g1} and \eqref{g3}--\eqref{g6}, but not
\eqref{g2}, \eqref{h1}--\eqref{h2}, or \eqref{p1}--\eqref{p3}.
From \eqref{g3}--\eqref{g6}, we find by the proof of
Proposition \ref{smithprop} that the backward Hugoniot through any
$(\tau_+,S_+)$ extends for all $S\in [S_+,-\infty)$ as a graph
of $\tau$ as a function of $S$, on which the Lax conditions hold
and $p$, $\sigma^2$ are monotone decreasing; likewise, the forward
Hugoniot is parametrized by increasing $S$, with $\sigma$ monotone decreasing
($\sigma^2$ increasing).
However, $p$ and $e$ do not approach zero as $S\to -\infty$, instead
becoming negative.

Indeed, the Hugoniot curve for \eqref{simplified}--\eqref{canon} 
may be solved explicitly as a cubic in $\tau$.
For example, setting $(\tau_+,S_+) = (1,0)$ and noting that
$\bar p(\tau,S) = -\bar e_\tau = {\frac {{{\rm e}^{S}}}{{\tau}^{2}}}-\tau$, 
 we find that $p_+ = 0$ and $e_+ = 3/2$. Substituting into 
$H(\tau,s) = e-e_+ + (1/2)(p+p_+)(\tau-\tau_+)=0$ yields
$$
H(\tau,S) = {\frac {{{\rm e}^{S}}}{\tau}}+S+1/2\,{\tau}^{2}-3/2+1/2\,
 \left( {\frac {{{\rm e}^{S}}}{{\tau}^{2}}}-\tau \right)  \left( 
\tau-1 \right)=0
$$
or, multiplying by $\tau^2$,
\be\label{hcub}
{\tau}^{3}+3\,{{\rm e}^{S}}\tau+(2S-3)\,{\tau}^{2}-\,{{\rm e}^{S}}=0. 
\ee
Thus an explicit solution $\tau(S)$ is available 
through the cubic formula.
Moreover, along the backward Hugoniot, we find as $S\to -\infty$ that
\eqref{hcub} reduces asymptotically to
$\tau^3 +2S\tau^2=0$, so that, asymptotically, $\tau(S)\sim -2S\to +\infty$,
and $e\sim 2S^2 \to +\infty$, $p\sim -\tau \to -\infty$.

\bpr\label{locprop}
For model \eqref{canon} (or, equivalently, \eqref{simplified}),
the backward Hugoniot through any
$(\tau_+,S_+)$ extends for all $S\in [S_+,-\infty)$ as a graph
of $\tau$ as a function of $S$, on which the Lax conditions hold,
$\tau$ is monotone increasing,
and $p$, $\sigma^2$ are monotone decreasing; likewise, $\sigma$ is
monotone decreasing and $S$ monotone increasing along the forward Hugoniot.
Yet, for $(\tau_+,S_+)=(1,0)$, there exist $(\tau_-,S_-)$
along the backward Hugoniot such that the Lax $1$-shock connecting
$(\tau,S)_\pm$ is Lopatinski (i.e., Hadamard) unstable.
\epr

\begin{proof}
By direct computation (see further calculations just below), we obtain
$$
\Big(\frac{-\bar e_{S\tau}}{\bar e_{\tau \tau}}-\frac{2\bar e_S}{[p]}\Big)_{
(\tau,S)=(\tau_+,S_+)}
= \frac{1}{3}-\frac{4}{-p_-} 
\sim \frac{1}{3}>0,
$$
as $S\to -\infty$, $\tau\to +\infty$, and $p_-\to -\infty$, yielding
the result by Proposition \ref{primeprop} and failure of \eqref{weak}. 
\end{proof}

\subsubsection{Estimating the stability transition}\label{localcomp}
From \eqref{canon}, we obtain
\besn{
p&=e^S/\tau^2 -\tau,\quad
e_S=e^S/\tau + 1,\quad
e_{\tau\tau} = 2e^S/\tau^3+1,\quad
e_{\tau S} = -e^S/\tau^2,
}
hence, evaluating at $(S_+,\tau_+) = (0,1)$,
\besn{
p_+ =0,\quad
(e_S)_+=2,\quad
(e_{\tau\tau})_+ = 3,\quad
(e_{\tau S})_+ = -1,
}
so that
$\phi := \left.- \frac{\bar e_S \bar e_{\tau \tau}}{\bar e_{\tau S}} \right|_{(S_+,\tau_+)} = -\frac{(2)(3)}{-1} = 6,$
and
$\theta:= 1+\sqrt{\frac{[p]/[\tau]}{p_{\tau}}}= 1+ \sqrt{\frac{p}{-3(\tau_- - 1)}}$.

Thus, change in stability correponds to change in sign of 
\besn{
\delta&:= p_+-p_- -\theta \phi
= -{\frac {{{\rm e}^{S}}}{{\tau}^{2}}}+\tau-6-2\,\sqrt {- \left( 3\,{
\frac {{{\rm e}^{S}}}{{\tau}^{2}}}-3\,\tau \right)  \left( \tau-1
 \right) ^{-1}}. 
}
Noting that $p\approx -\tau$,  
for $S_-$ large and negative,
so that $[p]/[\tau]\approx \frac{-\tau_-}{\tau_--\tau_+}$, 
and recalling that $c_+^2=-p_{\tau,+}=3$ 
we may well approximate
$\theta\approx 1 +\sqrt{ \frac {\tau_-}{6(\tau_--1})}\approx 1.6$,
so that the stability transition occurs approximately at $\tau_-= 9.6$
and, estimating $S\approx \frac12(3-\tau)$ for $e^{S}<<1$ in
\eqref{hcub}, $S_-\approx -3.3$. In fact, numerically we determine that $\delta = 0$ occurs
for $S_-\in (S_1,S_2) = (-3.3348293,-3.334829)$, $\tau_- \in (\tau_2,\tau_1) =(9.6589763, 9.6589769)$. 
See Figure \ref{fig638} below for a plot of the Hugoniot curve computed by
the Bisection Method, with the stability transition indicated on the curve.

\br
From $c_+^2=3$, $\sigma^2=-\frac{[p]}{[\tau]}\sim 
\frac{\tau_-}{\tau_--\tau_+}\sim 1+ \frac{1}{\tau_-}$,
and $c_-^2 =1+O(e^{S_-})$,
we verify directly the Lax conditions $c_-^2<\sigma^2<c_+^2$
for $S_->\! >1$, with the shock nearly characteristic at $U_-$.
\er

\begin{figure}[htbp]
 \begin{center}
$
\begin{array}{lr}
(a) \includegraphics[scale=0.3]{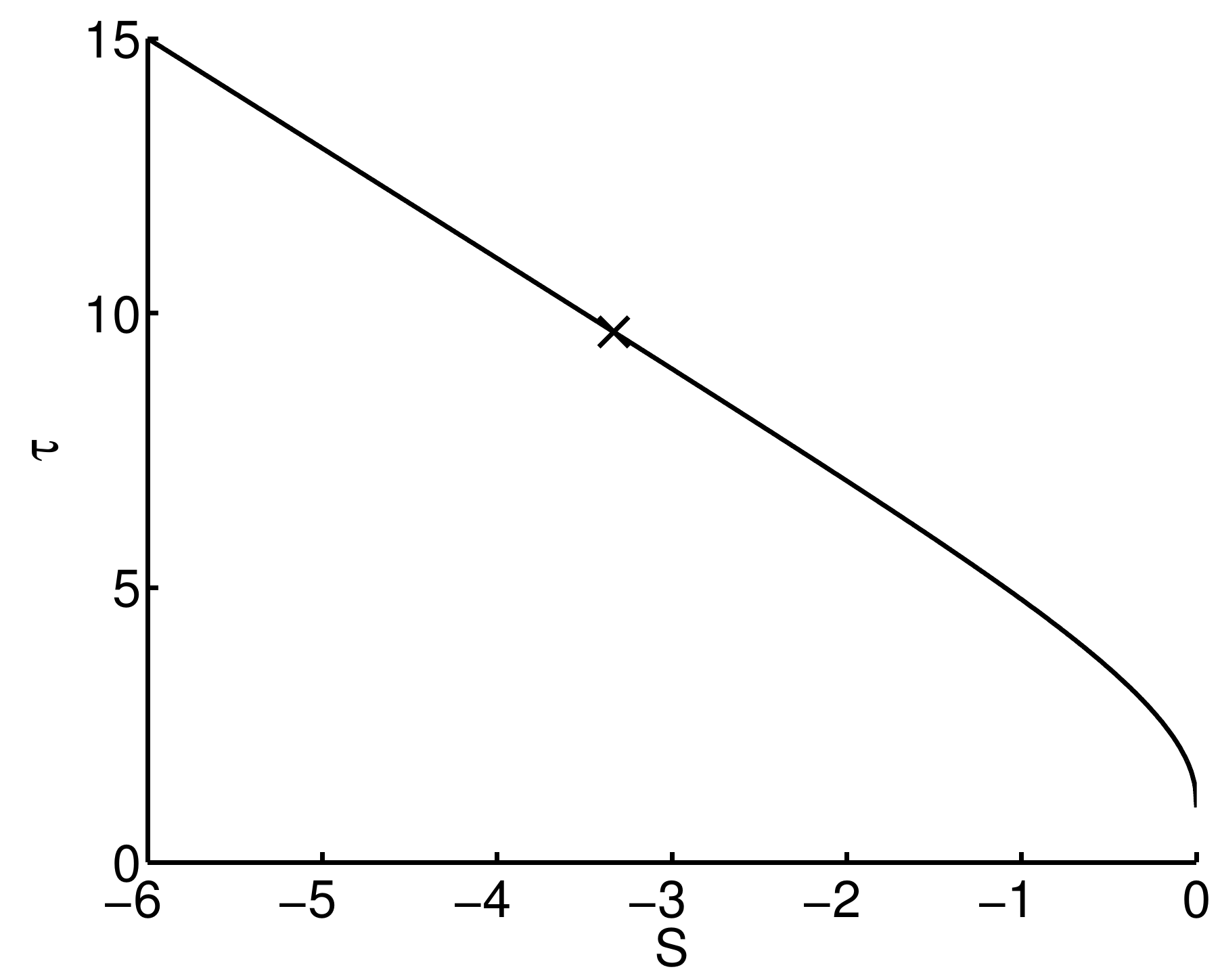}& (b) \includegraphics[scale=0.3]{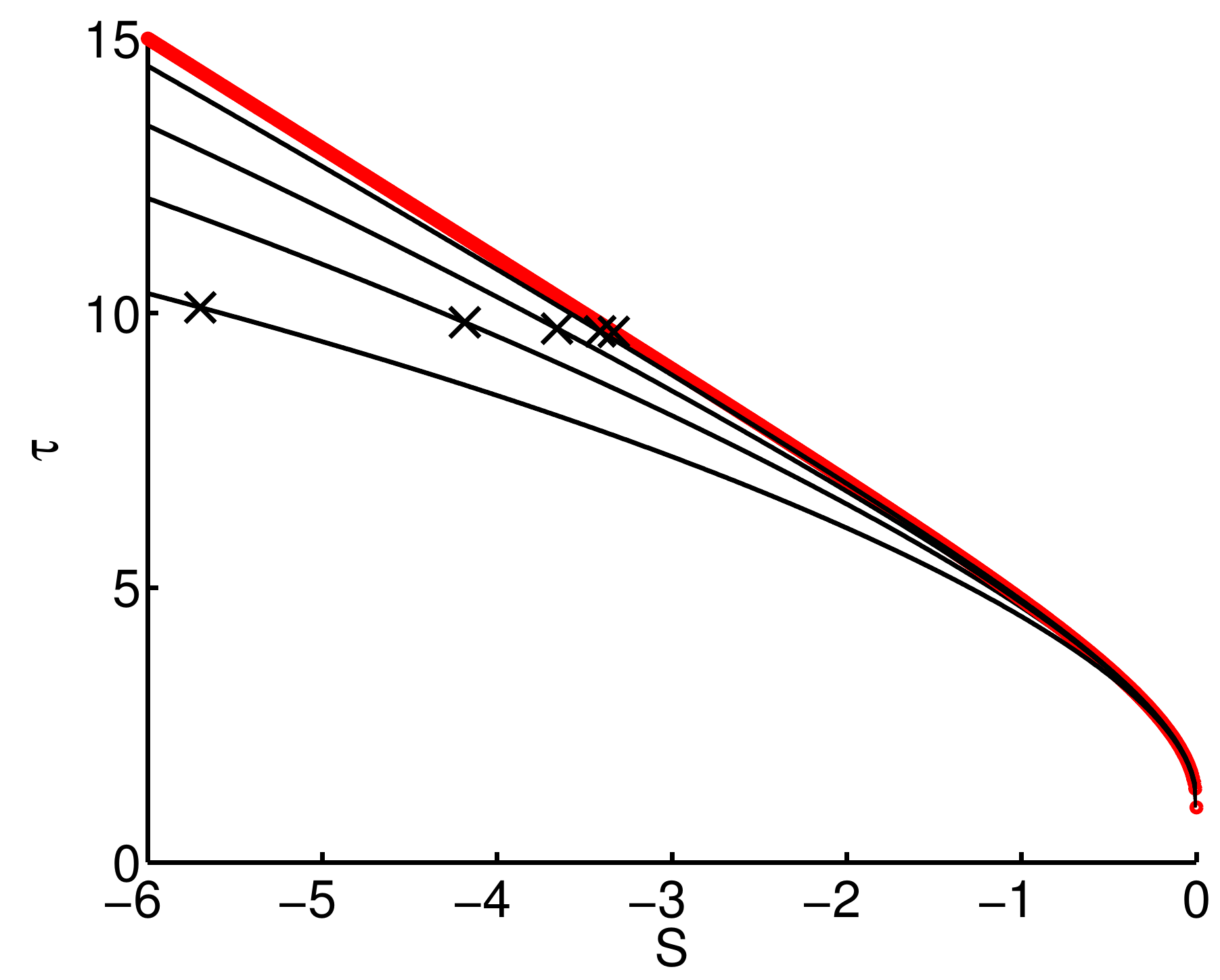}\
\end{array}
$
\end{center}
\caption{We mark transition to instability with a thick X. (a) Plot of the Hugoniot curve for the local model, $\bar e(S,\tau) = e^S/\tau+S+\tau^2/2$, $(\tau_0,S_0) = (1,0)$.  (b) Plot of the the Hugoniot curve for the local model with a thick red line and that of the global model for $C = 40,\ 100,\, 250$. Note that as $C\to \infty$, the Hugoniot curve of the local model matches well that of the global model.}
\label{fig638}\end{figure}


\subsection{Stable example}\label{s:stable}
For comparison, consider finallly the closely related model
\be\label{relmod}
\bar e(\tau,S)= e^S/\tau - C\tau +\tau^2/2,\quad C>\! >1.
\ee
Here, likewise the flow of \eqref{euler} is independent
of the value of $C$, so that it is equivalent to consider the
simpler version
\be\label{simplermod}
\bar e(\tau, s)= e^S/\tau +  \tau^2 /2 .
\ee

\bpr\label{simpleprop}
For model \eqref{simplermod} (or, equivalently, \eqref{relmod}),
the Hugoniot through any
$(\tau_+,S_+)$ extends for all $S\in (-\infty,-\infty)$ as a monotone graph
of $\tau$ as a function of $S$, on which the Lax conditions hold
precisely for $[\tau]<0$.
Moreover, all Lax shocks have positive signed Lopatinski determinant.
\epr

\begin{proof}
Again, from \eqref{g3}--\eqref{g6}, we find by the proof of
Proposition \ref{smithprop} that the backward Hugoniot through any
$(\tau_+,S_+)$ extends for all $S\in [S_+,-\infty)$ as a graph
of $\tau$ as a function of $S$, on which the Lax conditions hold
and $p$, $\sigma^2$ are monotone decreasing; by similar reasoning, the forward
Hugoniot is parametrized by increasing $S$, with $\sigma$ monotone decreasing
($\sigma^2$ increasing).
Indeed, the Hugoniot curve for \eqref{simplermod}
can be determined explicitly by solving
$$
\begin{aligned}
0=H(\tau,S)&= 
e^S/\tau +  \tau^2 /2 -e_+ 
-\frac12(e^S/\tau^2 - \tau + p_+)(\tau_+-\tau)\\
&=
e^S/\tau  -e_+ -\frac12(e^S/\tau^2  + p_+)(\tau_+-\tau)+ \frac12\tau \tau_+\\
\end{aligned}
$$
or, using the identities $p_++\tau_+=\frac{e^{S_+}}{\tau_+}$
and $e_++ \frac{p_+\tau_+}{2}= \frac32 \frac{e^{S_+}}{\tau_+}$,
$e^{S-S_+}= \Big(\frac{3\tau_+-\tau}{3\tau -\tau_+}\Big)\frac{\tau^2}{\tau_+^2}$,
for
$$
S-S_+=\log \Big(\frac{3\tau_+-\tau}{3\tau- \tau_+}
\Big)\frac{\tau^2}{\tau_+^2}.
$$
Computing $(d/d\tau)(S-S_+)= -\frac{6(x-1)^2}{x(3-x)(3x-1)\tau_+}<0$
for $x:=\tau/\tau_+$ and $1/3<x<3$, we find that $S$ is a monotone
function of $\tau$ along the entire (forward and backward) Hugoniot
curve.
But, monotonicity along the forward Hugoniot implies
stability by Remarks \ref{monrmk} and \ref{vasprep}.
\end{proof}

\begin{figure}[htbp]
 \begin{center}
$
\begin{array}{lcr}
(a) \includegraphics[scale=0.3]{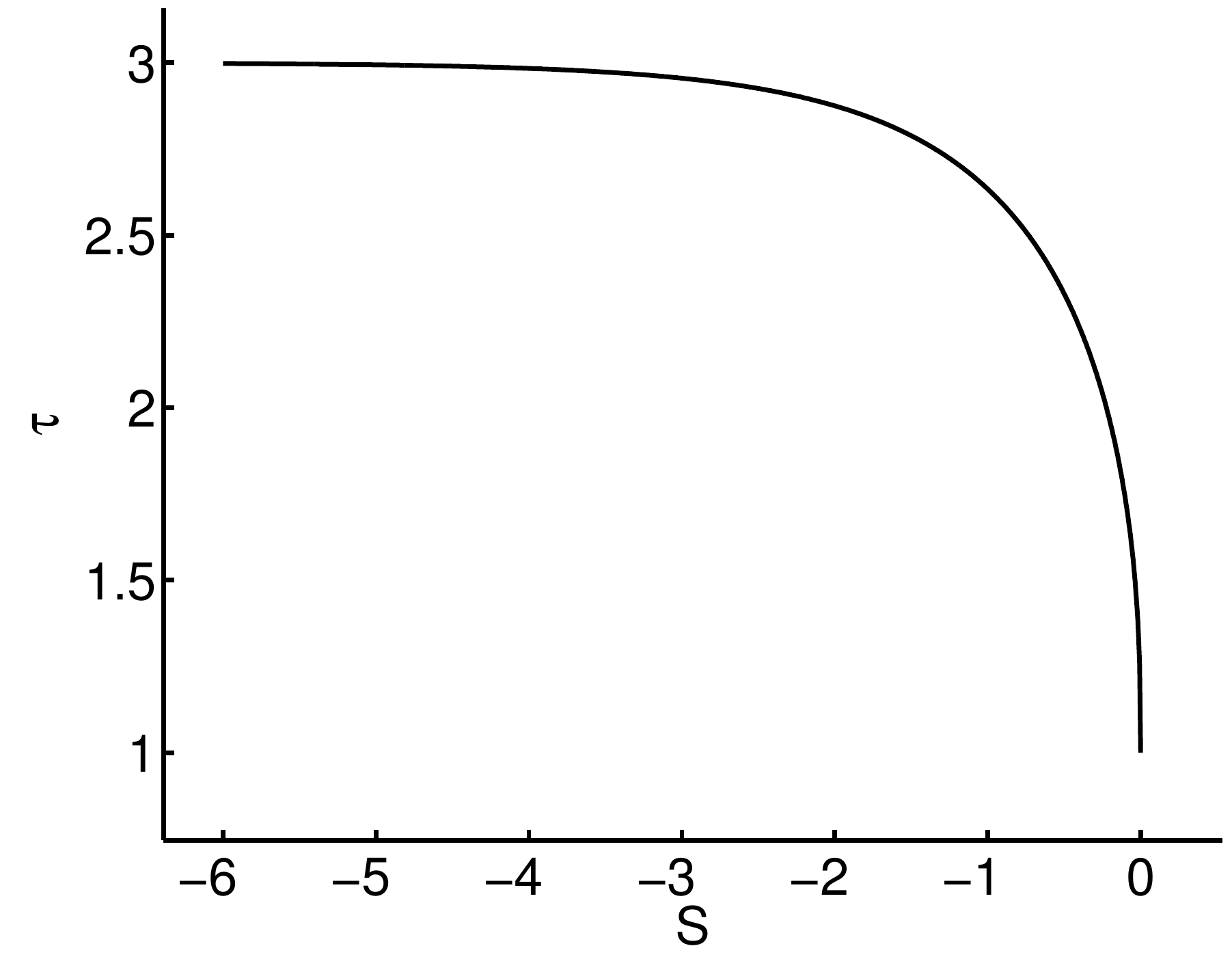} \quad (b)\includegraphics[scale=0.3]{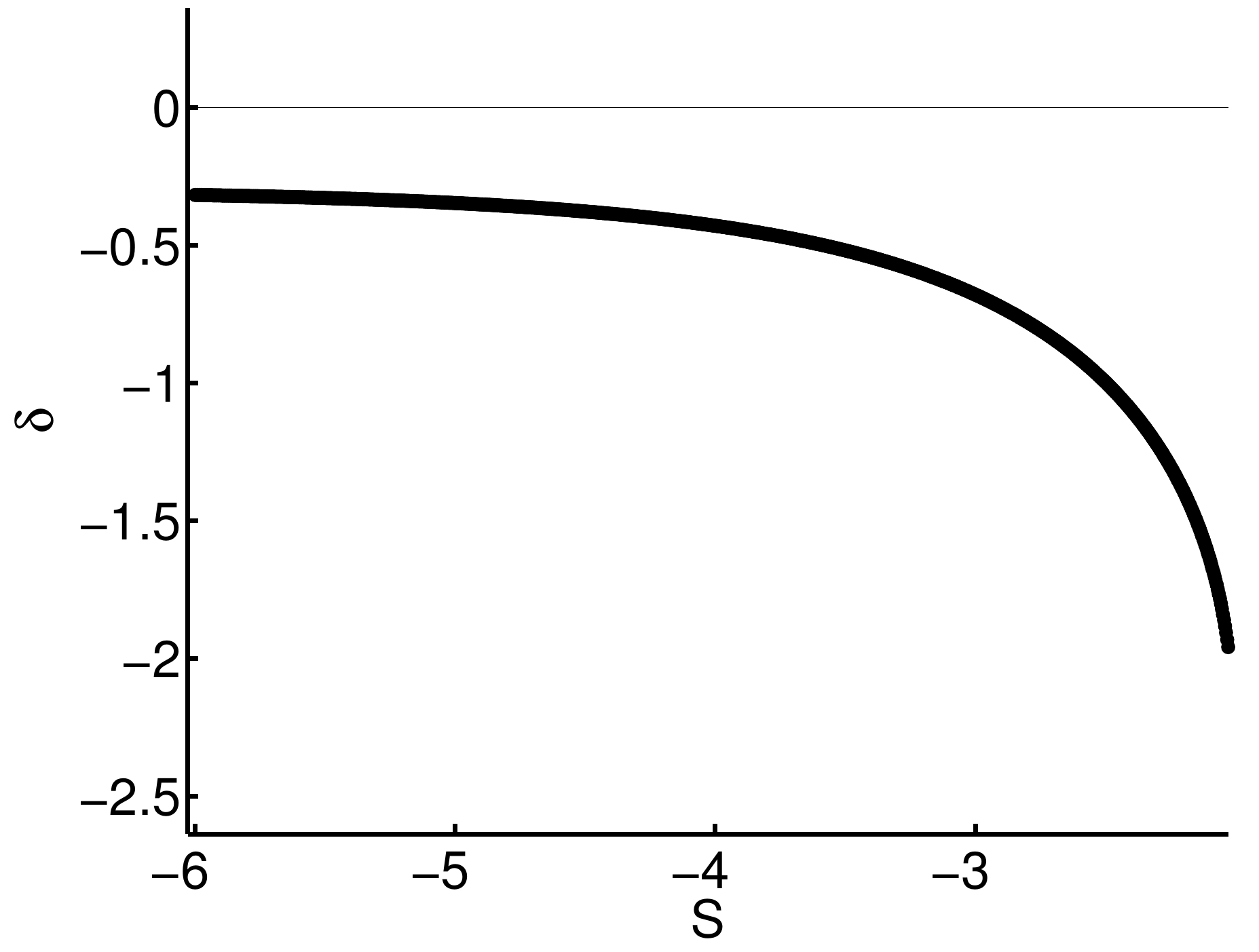}
\end{array}
$
\end{center}
\caption{(a) Plot of the Hugoniot curve for the local model $\bar e(S,\tau) = e^S/\tau-C\tau+\tau^2/2$ with $C = 100$, $(\tau_0,S_0) = (1,0)$. (b) Plot of $\delta$ against $S$ where $\delta < 0 $ corresponds to stability.
}
\label{fig640}\end{figure}

\bigskip


\centerline{\bf PART II. VISCOUS STABILITY}
\par\medskip
Stepping outside the inviscid setting, we can extract further interesting information
from the sign of the signed Lopatinski determinant. This fact and the Evans function, to
which it refers, are  discussed in Section 5. Section 6 presents the findings that
were summarized as Numerical Observation A in Sec.\ 1. The phenomenologically most interesting
piece of Part II might be  Section 7, where the findings of Numerical Observation B are detailed,
on the rich stability behavior of shock waves in an 
artifically designed system. 

\section{The Evans function}\label{s:evans}
We now turn to the study of viscous stability, or stability
of continuous traveling shock fronts
\be\label{prof}
w(x,t)=\bar w(x-\sigma t),\quad \lim_{z\to \pm \infty} \bar w(z)=U_\pm
\ee
of a hyperbolic--parabolic system of conservation laws
\be\label{vgencon}
f^0(w)_t+f(w)_x=(B(w)w_x)_x,
\ee
corresponding to an inviscid shock $(U_-,U_+,\sigma)$
of \eqref{gencon}, that is, a solution of the traveling-wave ODE
\be\label{twode}
B(\bar w)\bar w'= f(\bar w) -\sigma f^0(\bar w) -(f(U_-)-\sigma f^0(U_-)).
\ee

Recall \cite{GZ,ZH,MaZ1,MaZ2,Z1,BLeZ,Z5} that, under rather general circumstances,
shared by most equations of continuum mechanics-- in particular, the 
equations of gas dynamics, MHD, and viscoelasticity--
Lax-type traveling-wave profiles $\bar w$, when they exist, 
(i) are generically smooth and converging at exponential rate 
as $z\to \pm \infty$ to their endstates $U_\pm$,
and (ii) are nonlinearly orbitally stable as solutions of \eqref{vgencon}
if an only if they satisfy an \emph{Evans function condition}:
\be\label{evcon}
\hbox{\rm 
$D(\lambda)$ has no zeros on $\Re \lambda\ge 0$ except for a
multiplicity one zero at $\lambda=0$,}
\ee
where the Evans function $D$ is a Wronskian at $x=0$ of 
appropriately chosen bases
of solutions of the linearized eigenvalue equations about $\bar w$ that decay
at plus and minus spatial infinity, respectively, 
that is analytic in $\lambda$ on a neighborhood of $\Re \lambda\ge 0$,
with complex symmetry $\overline{D(\lambda)}=D(\bar \lambda)$,
whose zeros agree on $\{\Re \lambda \ge 0\}\setminus \{0\}$
in location and multiplicity with the eigenvalues of the 
linearized operator about the wave.

\subsection{Viscous vs. inviscid stability}\label{s:vs}
Like inviscid stability, viscous stability holds always 
for Lax type shocks in the small-amplitude limit \cite{ZH,HuZ3,PZ,FS1,FS2}.
Moreover, as found in \cite{GZ,ZS}, we have the fundamental relation
\be\label{ZSagain}
 D(\lambda) = \lambda \nu \delta +O(\lambda^2)
\ee
for $\lambda$ sufficiently small, where $\nu$ is a Wronskian associated
with the linearized traveling-wave ODE, with nonvanishing of $\nu$ corresponding
to transversality of the connecting viscous profile, and $\delta$ is
the Lopatinski determinant \eqref{lop} associated with the inviscid 
stability problem.
Without loss of generality, normalize $\nu$ and $\delta$ 
so as to be real.
Then, moving along the Hugoniot curve in the direction of increasing amplitude,
we find that, so long as there persists a transversal traveling-wave connection,
so that $\nu $ cannot change sign, that real zeros of the 
Evans function can cross through the origin into the unstable half-plane 
\emph{if and only if} the Lopatinski determinant $\delta$ changes sign, that is,
there is a corresponding stability transition at the inviscid level
\cite{GZ,ZH,ZS,Z1,Z3}.

More precisely, we see that, assuming persistence of transversal connections,
inviscid stability transitions correspond
at viscous level to the passage of an odd number of unstable roots through
the origin, so can never occur without a corresponding viscous stability
transition.
On the other hand, there can occur viscous stability transitions corresponding
to passage of an even number of real roots, or of
one or more complex conjugate pairs of unstable roots through the imaginary
axis, which are not detected at the inviscid level.
Thus, \emph{viscous stability is a logically stronger condition than
inviscid stability} \cite{ZS}.
Our goal in this part of the paper is to investigate whether this
logical implication can be strict for systems \eqref{vgencon} possessing 
an associated convex entropy, and in particular whether complex conjugate
pairs of roots can cross, signaling Hopf bifurcation to time-periodic
``galloping'' behavior \cite{TZ1,TZ2,TZ3,BeSZ}.

\subsection{Numerical stability analysis}\label{s:numevans}

The Evans function, or, for that matter, the underlying profile
\eqref{prof}, is only rarely computable analytically.
However, as described, e.g., in \cite{Br,BrZ,HuZ1,Z4},
these may be approximated numerically in a well-conditioned and 
efficient manner.
Here, we carry out numerical experiments using MATLAB and
the MATLAB-based Evans function package STABLAB \cite{BHZ2}
developed for this purpose by Barker, Humpherys, Zumbrun, and collaborators. 
More precisely, we study the {\it integrated Evans function} 
$\tilde D(\lambda)$, a Wronskian associated with the integrated eigenvalue
equations, whose zeros 
on $\{\Re \lambda \ge 0\}\setminus \{0\}$ likewise agree
in location and multiplicity with the eigenvalues of the 
linearized operator about the wave, but which does not have a zero
at the origin, obeying a small-frequency expansion
\be\label{ZSint}
 \tilde D(\lambda) = \nu \delta +O(\lambda)
\ee
in place of \eqref{ZSagain}.

Descriptions of the underlying 
algorithms and other information useful in reproducing
these computations may be found in Appendix \ref{s:protocol}.

\section{Gas dynamical examples}\label{s:vgas}

The compressible Navier--Stokes equations, augmenting \eqref{euler}
with the transport effects of viscosity and heat conduction,
appear in Lagrangian coordinates as
\ba\label{NS}
\tau_t-v_x&= 0,\\
v_t+ p_x&= \big(\frac{\mu v_x}{\tau}\big)_x ,\\
(e+v^2/2)_t+(vp)_x&=
\big(\frac{\mu vv_x}{\tau}\big)_x+\big(\frac{\kappa T_x} {\tau}\big)_x,
\\
\ea
where $\tau$ denotes specific volume,
$v$ velocity, $e$ specific internal energy, $p$ pressure,
and $T$ temperature,
and $\mu$ and $\kappa$ are coefficients of viscosity and heat  
conductivity, here taken constant. 

To close these equations requires \emph{both} a pressure law
$p=\hat p(\tau,e)$ and a temperature law $T=\hat T(\tau,e)$,
in contrast to the case of the Euler equations \eqref{euler},
that is, a complete equation of state.
Indeed, the most convenient formulation for our purposes here
is in terms of the variable $w:=(\tau, v,T)$ replacing $e$ with $T$.
We will use  equations of state given in the form
\be\label{nslaws}
p=\check p(\tau,T),\quad e=\check e(\tau,T)
\ee
arising through the Helmholtz energy formulation described
in Appendix \ref{s:helmholtz}, for which the equations
take the convenient block structure described in Appendix \ref{s:egen}.
Note that besides hyperbolicity, $\bar p_\tau<0$, 
the Navier--Stokes equations \eqref{NS} require
also parabolicity, 
$$
\bar e_{SS}=(\partial T/\partial e)(\partial e/\partial S)>0, 
$$
in order to be well-posed, which is the condition needed to
obtain laws \eqref{nslaws} from a given energy function $e=\bar e(\tau,S)$
of the form considered in previous sections.

Substituting \eqref{prof} into \eqref{vgencon} and 
solving the first equation of \eqref{twode} for 
\be\label{tausolve}
 w_1=\tau =\tau_-+ \sigma^{-1}(v_--v) 
\ee
as described in Appendix \ref{s:egen},
we obtain for a profile corresponding to an inviscid shock
$(U_-,U_+,\sigma)$ of \eqref{euler} the traveling-wave equations
\ba\label{twodens}
\tau^{-1}\bp \mu & 0 \\\mu v & \kappa\ep
\bp v\\T\ep'=
\bp
-\sigma(v-v_-) + p-p_-
\\
-\sigma(e+v^2/2-e_--v_-^2/2) + pv-p_-v_-
\ep.
\ea

\bpr[\cite{Gi}]\label{Giprop}
For a $C^2$ energy function $\bar e(\tau,S)$ satisfying
\eqref{g3}--\eqref{g6}, every Lax type shock of \eqref{euler}
has a unique continuous shock profile 
\eqref{prof} of \eqref{NS}--\eqref{nslaws}, which is a transversal
connection of \eqref{twodens} converging exponentially to endstates $U_\pm$ 
as $x\to \pm \infty$ in up to two derivatives.
\epr

\begin{proof}
More generally, Gilbarg shows that, assuming the Weyl
condition $\bar p_S>0$, a Lax-type inviscid $1$-shock 
for which the pressure at the righthand state $U_+$ is 
minimal
among all connections to $U_-$ with the same shock speed $\sigma$
possesses a unique viscous profile \eqref{prof}.
Since \eqref{g3}--\eqref{g6} imply uniqueness of inviscid Lax 
shocks for a given speed, we thus obtain existence of a unique 
profile; transversality then follows by the observation that
this connection must be of saddle-node type, so transversal whenever
it exists and exponentially decaying by standard stable/unstable manifold
theory.
See also \cite{W} and \cite[Appendix C]{MP}.
\end{proof}

\bpr[Kaw,KSh]\label{kawprop}
For a convex $C^2$ energy function $\bar e(\tau,S)$ satisfying
\eqref{g3}--\eqref{g6}, \eqref{NS} possesses
the convex viscosity-compatible entropy $\hat S(\tau,e)$,
globally on $\tau, e>0$.
\epr

\begin{proof}
See, e.g., \cite{Kaw,KSh,MP,Z1,Z2}.
\end{proof}

\subsection{Numerical details}\label{s:details}
Despite the familiarity of the gas dynamical equations \eqref{NS},
and the straightforward analytical treatment under the assumptions
we impose, we find the numerical Evans function analysis of our
gas dynamical examples to be quite challenging, due to the 
presence of multiple scales.
Indeed, this turns out to be the most computationally intensive problem that
we have so far considered, even including the famously intensive
problem of detonation stability considered in \cite{Er3,LS,BHLyZ1,BZ2}.
We thus find it necessary to include some extra details on how
the computations are carried out, despite that this is not
our main emphasis here.

\subsubsection{Computation of profiles}\label{s:gasprof}
To approximate numerically the profiles guaranteed by
Proposition \ref{Giprop} for a given equation of state
$\bar e(\tau,S)$,
we see two basic possible approaches.  The more direct, but
numerically somewhat cumbersome approach is
to evaluate the functions $\check p(\tau,T)=\bar p(\tau, \check S(\tau,T)$ 
and $\check e(\tau, T)=\bar e(\tau, \check S(\tau, T)$
by numerically inverting the relation $T=\bar e_S(\tau,S)$ 
to solve for $S=\check S(\tau,T)$, using $\bar e_{SS}>0$.
As computation of the profile is a one-time cost, the associated
cost in function calls is perhaps not too serious, but
this still seems preferable to avoid.

The approach we follow here is, rather, to compute the profile in $(\tau,S)$ coordinates
using the change of variables formula
$ \bp \tau \\S\ep'=
\bp   -\sigma & 0\\ \bar e_{S\tau} &\bar e_{SS}\ep^{-1} \bp v\\T\ep'
$
deriving from $T'=\bar e_{S\tau}\tau'+ \bar e_{SS}S'$ and
the relation $\tau'=-\sigma^{-1} v'$ obtained from \eqref{tausolve},
yielding, finally,
\be\label{ftwode}
 \bp \tau \\S\ep'=
\tau
\bp   -\sigma & 0\\ \bar e_{S\tau} &\bar e_{SS}\ep^{-1} 
\bp \mu & 0 \\\mu v & \kappa\ep^{-1}
\bp
-\sigma(v-v_-) + p-p_-
\\
-\sigma(e+v^2/2-e_--v_-^2/2) + pv-p_-v_-
\ep.
\ee
From $(\bar \tau, \bar S)(x)$, we can recover all other variables
through \eqref{tausolve} and the equation of state.

A final detail is that in the large-amplitude limit $S\to -\infty$,
the traveling-wave ODE \eqref{ftwode} features a wide separation
of scales.
Specifically, the linearized equations about endstates $U_\pm$,
\ba\label{linode}
\bp \tau \\S\ep'&=
\bp -\frac{1}{\sigma}& 0\\ \frac{\bar e_{S\tau}}{\bar e_{SS}} & 
\frac{1}{\bar e_{SS}}\ep_\pm
\bp \frac{1}{\mu}& 0\\ -\frac{v}{\kappa}& \frac{1}{\kappa}\ep_\pm
\bp \sigma^2+ c^2 & \bar p_S\\
v(2\sigma^2 -c^2) & -\sigma T + \bar p_{S}v\ep_\pm 
\bp \tau \\S\ep,
\ea
for global model \eqref{ceg2}
in the regime $1<<|S_-|<<C$ considered in the crucial stability
transition regime, is asymptotic at $U_-$ ($x\to -\infty$) to
(computing $|\sigma|\sim 1$, $c=-\bar e_{\tau \tau}\sim 1$,
$\bar e_{S}\sim 1$, $\bar e_{S\tau}\sim C^{-1}$, $\bar e_{SS}\sim C^{-2}$)
$$
\bp \tau \\S\ep' \sim
\bp \frac{2}{\mu} & 0\\ -\frac{v_- C}{\kappa} & \frac{C^2}{\kappa}\ep
\bp \tau \\S\ep,
$$
which exhibits growth modes with exponential rates $\sim 1/\mu$ and 
$ C^2/\kappa$ differing by two or more orders of magnitude.
Similar but more extreme behavior is inherited by the local model 
\eqref{canon} in the regime $\tau_-\to +\infty$, 
$S_-\sim \frac12(3-\tau)\to -\infty$ roughly corresponding to $C\to \infty$
for \eqref{ceg2}, with growth rates 
$\sim 1/\mu$ and $\tau_- e^{\tau_-/2}/\kappa$.

We find that this gives trouble for the standard method
(described further in Appendix \ref{s:num}) of solving
for profiles with a boundary-value solver with projective boundary conditions;
indeed, we are unable to compute by this method 
beyond approximately the range $ |C|\le 75$ for \eqref{ceg2}
or $\tau_-\le 19$ for \eqref{canon}.\footnote{
This is with initial guess given by a simple $\tanh$ interpolation.
We can go slightly farther by continuing along the Hugoniot curve,
using previous profile as initial guess; however, the computations are still
extremely slow.}

Above these values, we rely on a more primitive, shooting method, starting at a 
point distance $10^{-8}$ from the saddle $U_+$ in the direction of
the direction of the stable subspace of \eqref{linode}, and integrating
backward in $x$ the traveling-wave ODE \eqref{ftwode} 
using an appropriate stiff numerical ODE solver until the solution
comes within $10^{-8}$ of the repellor $U_-$.
This simple solution gives excellent results for essentially arbitrary
values of $S_-$.

%

\subsubsection{Evans computations}\label{s:gasevans}
In computing the Evans function, we 
face the same numerical stiffness issues that we faced for the
profile equation
(recall that at freqency $\lambda=0$, the eigenvalue equation reduces
to the linearized traveling-wave ODE).
Though it is standard practice for us to compute the Evans function by
shooting, we find it necessary to carry out this shooting algorithm
using a stiff ODE solver rather than the standard RK45.
For further discussion, see Appendix \ref{s:num}.

A related issue is that for $C> 10$ or so in \eqref{ceg2}, or
for $\tau_- > 6$ in \eqref{canon}, the standard Evans function experiences 
exponential growth/decay arising from spatial variation in the coefficient
matrix of the Evans (i.e., first-order eigenvalue) system so great
as to go out of numerical range- in this case, returning a value
that is to numerical approximation identically zero.
Similar phenomena have been seen to arise in
other numerically difficult Evans computations \cite{BZ1,BZ2,BHLyZ1}.
We follow the simple solution here (mentioned also in the cited references) 
of turning off growth; that is, computing using the polar coordinate method
of \cite{HuZ1}, we suppress variation in the radius and simply evolve
orthonormal bases of the subspaces of decaying solutions at plus and minus
spatial infinity (the ``polar no radial'' 
setting in the STABLAB Evans package). 
This amounts to the continuous orthogonalization method of Drury \cite{Dr},
which preserves winding numbers/zeroes of the Evans function, but loses
the desirable property of analyticity \cite{Z4}.
(It is in fact $C^\infty$ in $\lambda$.)

The same numerical stiffness issues prevent degrade our usual
estimates excluding unstable eigenvalues outside a semicircle of radius
$R$. 
In this section, therefore, we often compute for a fixed radius $R$
without guarantees that there are no unstable eigenvalues of larger
modulus.

\br
The presence of the very slowly decaying
term $ e^{S/C^2-\tau/C}$ in the Evans function coefficient arixing in
our global counterexample \eqref{ceg2}
degrades the standard convergence estimates described in, e.g., \cite{Z4,Z5},
wich depends on uniform exponential convergence of the coefficient matrix
of \eqref{Zeq} to its limits as $x\to \pm \infty$. 
However, uniform-in-$C$ estimates are easily recovered
similarly as (but much more
simply than) in \cite{HLZ}, by the incorporation
of an additional ``tracking'' argument in the slowly-varying far-field
regime.
%
A similar argument, though we do not carry it out here,
may be used to rigorously verify convergence as $C\to \infty$
of the Evans function for the global counterexample (suitably normalized)
to the Evans function
for the local example of Section \ref{s:local}.
\er

\subsection{Global counterexample}\label{s:vglobal}
Returning to the equation of state 
\be\label{ceg3}
\bar e(\tau,S)= \frac{e^{S}}{\tau}+ C^2 e^{ S/C^2 -\tau/C},
\quad C>\! >1
\ee
shown to exhibit inviscid instability in Theorem A,
we report the results of our numerical Evans function investigations.
%
%
In our study we examined the parameter set $S_- = \{-1,-1.1,...,-50\}$, 
$(\tau_+,S_+)=(1,0)$,
 for $C = 10$. We solved the profile using Matlab's boundary value solver bvp5c. In addition, we computed the Evans function on a smaller radius for $(\tau_+,S_+)=(1,0)$
and several values of $S_-$ including $S_- = -10^{5}$ for $C = 100$ with $R = 4$, solving for the profile using Matlab's implicit ODE solver, ode15s. 

The results of these investigations were, in all cases, that, as shock amplitude was increased (i.e.,
$S_-$ decreased), there was, similarly as in the inviscid stability analysis,
a single stability transition consisting of a real root of multiplicity one passing
through the origin into the unstable complex halfplane $\Re \lambda>0$, where it remained thereafter
as the only unstable eigenvalue. 
In particular, there were seen no complex unstable eigenvalues at any point in the process.
That is, in this case, viscous and inviscid stability predictions appear to agree.

For $C=10$,
the computed inviscid stability transition occurs approximately at $S_- = -23.2$, 
whereas the computed viscous stability transition occurs for $S_- \in( -22.1,-22.2)$;
recall that, by the analytical theory, these two values should agree,
as parameters for which a real root passes through the origin.
Except for this difference in location of the transition point,
the results of inviscid and viscous numerical stability analyses agreed.
The relatively large size of this discrepancy we ascribe
to the stiffness of the underlying eigenvalue equation both for viscous and inviscid cases;
compare for example to the excellent agreement found in the investigations of the nonstiff system
considered in Section \ref{s:designer}.

The study was further challenged by the slow rate at which the Evans function converges to its high frequency behavior $D(\lambda)\sim C_1e^{C_2\sqrt{\lambda}}$, as described
in \eqref{hfas}, Appendix \ref{s:num}. 
In tables \ref{tb_error1} and \ref{tb_error2} we provide data about curve fitting the Evans function with $C_1 e^{C_2 \sqrt{\lambda}}$ as $\lambda \to \infty$. As seen in the tables, the needed radius for convergence for large $C$ exceeds the range where numerical computations can be carried out accurately. Consequently, we arbitrarily chose $R = 250$ for the Evans function study. Figure \ref{fig664} demonstrates typical output from our study. 

\begin{table}[!b]
\begin{tabular}{|c|c|c|c|}
\hline
 $R$&error& $C_1$&$C_2$\\
\hline 
\hline
2& 0.8726 & 0.03233&2.427 \\
 4&1.693& 0.002496& 2.996\\
   8&  2.103& 3.386e-05& 3.639\\
  16&  6.33& 3.659e-08& 4.281\\
   32&  12.47& 1.476e-12& 4.816\\
  64&  17.38& 1.784e-18& 5.108\\
 128&  4.211& 3.922e-26& 5.171\\
 256& 0.9991& 7.957e-36& 5.051\\
  512& 0.9949& 3.594e-48&4.828\\
 1024&  0.9994& 8.181e-65& 4.611\\
 2048& 0.9945&2.664e-88& 4.456\\
4096&  0.9805& 5.719e-122& 4.362\\
  8192& 0.948& 3.166e-170& 4.312\\
1.638e+04& &  NaN& NaN\\
\hline
\end{tabular}
\caption{For $C = 10$, $S_- = -15$ in the global model we compare how well $C_1 e^{C_2\sqrt{\lambda}}$ approximates $D(\lambda)$ for $|\lambda| = R$. Here $R$ is the radius of the semicircle on which we compute the Evans function, $C_1$ and $C_2$ are the curve fit returned, and the convergence error is the maximum relative error between $D(\lambda)$ and $C_1 e^{C_2\sqrt{\lambda}}$ evaluated at $\lambda = R,\ Ri$.}
\label{tb_error1}
\end{table}

\begin{table}[!b]
\begin{tabular}{|c|c|c|c|}
\hline
 $R$&error& $C_1$&$C_2$\\
\hline 
\hline
2 & 0.3373 & 0.711 &  0.2412\\
4 &  0.3457 &  0.9109 &  0.04668\\
8 &  0.287 & 1.268 &  -0.08408\\
16 &  0.2102 &  1.876 &  -0.1572\\
  32 &  0.1479 &  3.003 &  -0.1944\\
 64 &  0.1055 & 5.439 &  -0.2117\\
 128 &  0.0779 &  12.07 &  -0.2202\\
256 &  0.05837 &  36.77 &  -0.2253\\
\hline
\end{tabular}
\caption{For $C = 2$ in the global model we compare how well $C_1 e^{C_2\sqrt{\lambda}}$ approximates $D(\lambda)$ for $|\lambda| = R$. Here $R$ is the radius of the semicircle on which we compute the Evans function, $C_1$ and $C_2$ are the curve fit returned, and the convergence error is the maximum relative error between $D(\lambda)$ and $C_1 e^{C_2\sqrt{\lambda}}$ evaluated at $\lambda = R,\ Ri$.}
\label{tb_error2}
\end{table}


\begin{figure}[htbp]
 \begin{center}
$
\begin{array}{lcr}
(a) \includegraphics[scale=0.25]{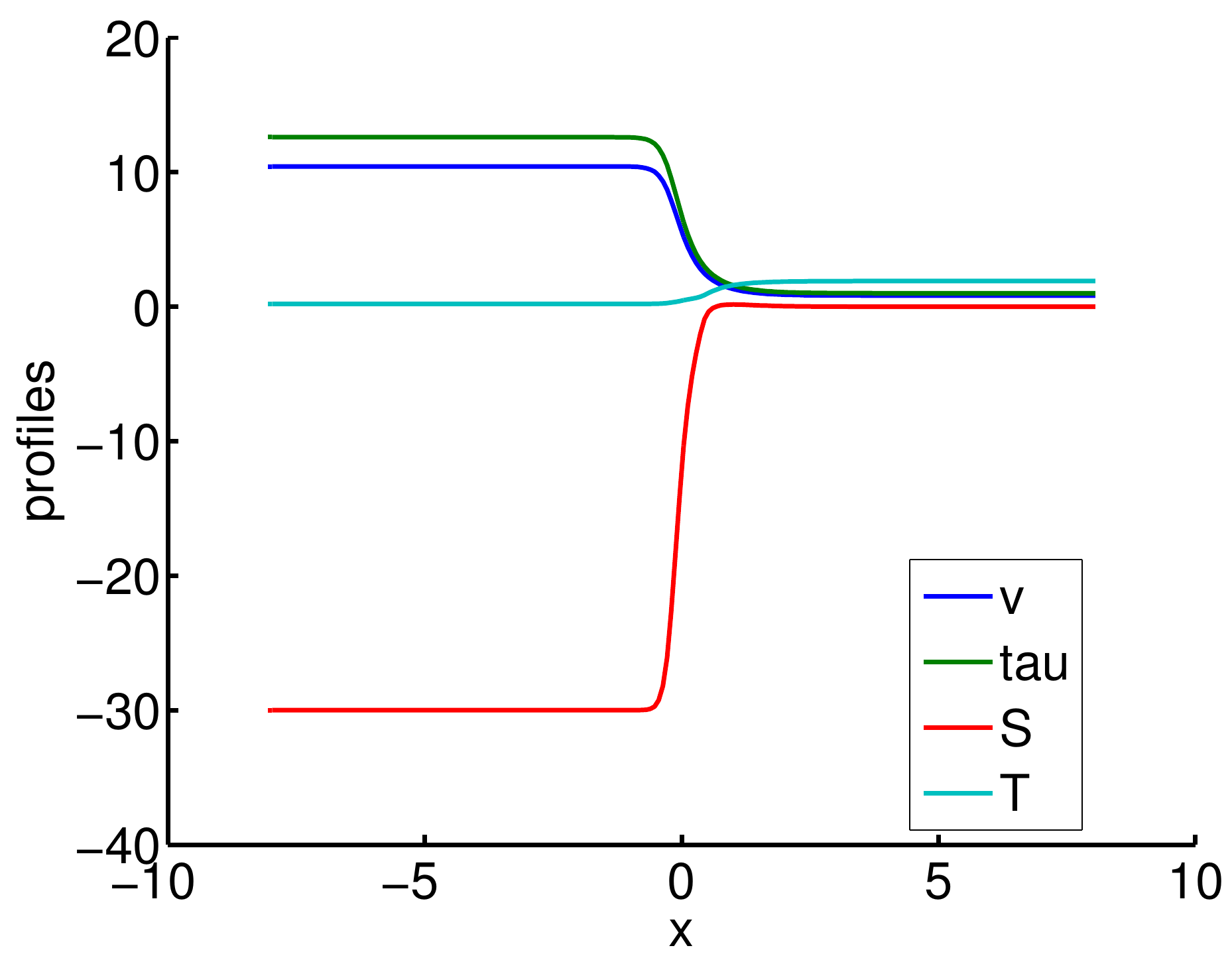}& (b)  \includegraphics[scale=0.25]{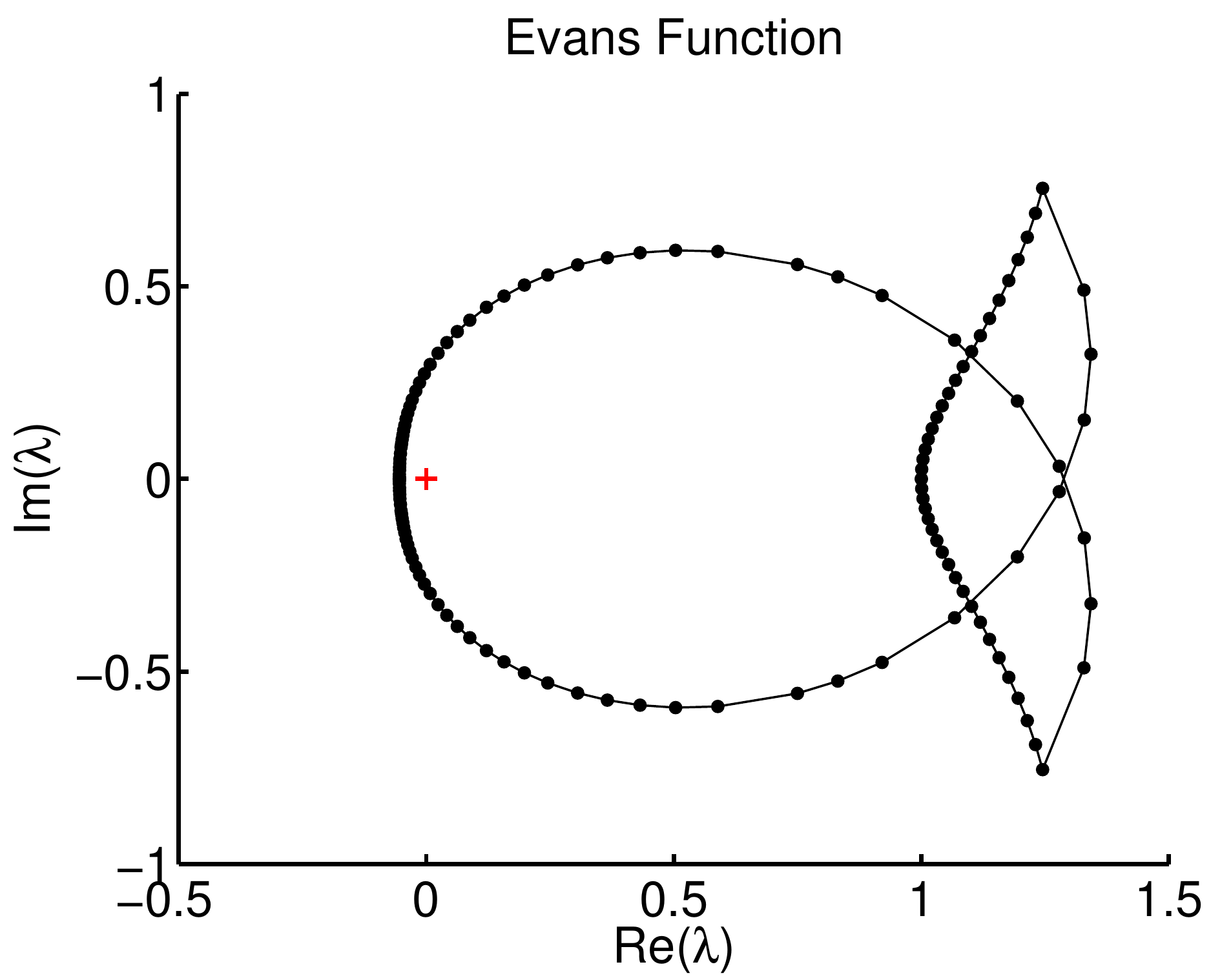}& \includegraphics[scale=0.25]{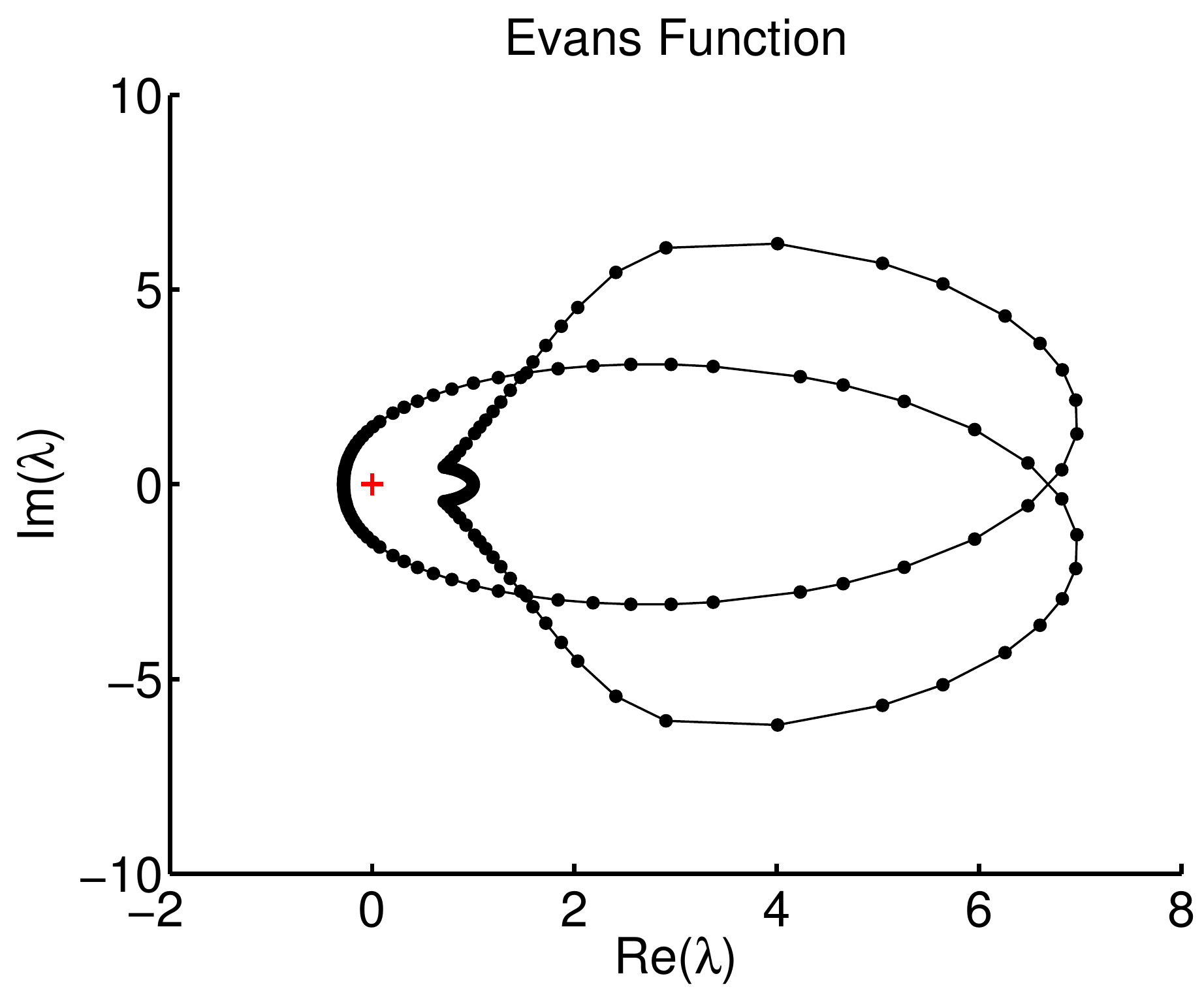}
\end{array}
$
\end{center}
\caption{Viscous study results for $\bar e(S,\tau) = e^S/\tau + C^2e^{C^{-1}(C^{-1}S-\tau)}$ with $S_- = -30$, $C = 10$, $\mu= \kappa = 1$, $(\tau_0,S_0) = (1,0)$. (a) Viscous profile. (b) 
Integrated Evans function computed with the adjoint polar coordinate method, evaluated on a semi-circle with radius $R =1$. Winding number is 1. (c) As in (b) but with $R = 250$.}
\label{fig664}\end{figure}

%
%

\subsection{Local counterexample}\label{s:vlocal}
Next, we turn to the local model \eqref{canon} obtained by expansion
of the global model in the $C\to \infty$ limit,
describing the results of a corresponding numerical Evans function investigation in that case.
In our study we examined the parameter set $S_- = \{-1, -1.1,...,-8.1\}\backslash \{-3.3\}$. 
Numerically, we found that both the inviscid and viscous stability transitions 
occurred around $S_- = -10/3$. 
For all parameters covered in our study, the viscous
Evans function had one root or no root according as the shock was inviscid unstable or stable.
As with the global model, we could not carry out a numerical high frequency study, and so we arbitrarily chose $R = 250$ for the contour radius in our Evans function computations. Figure \ref{fig666} demonstrates typical output from our study. 
In conclusion, the results of our numerical study 
indicate in this case too that viscous and inviscid stability agree.

\begin{figure}[htbp]
 \begin{center}
$
\begin{array}{lcr}
(a) \includegraphics[scale=0.25]{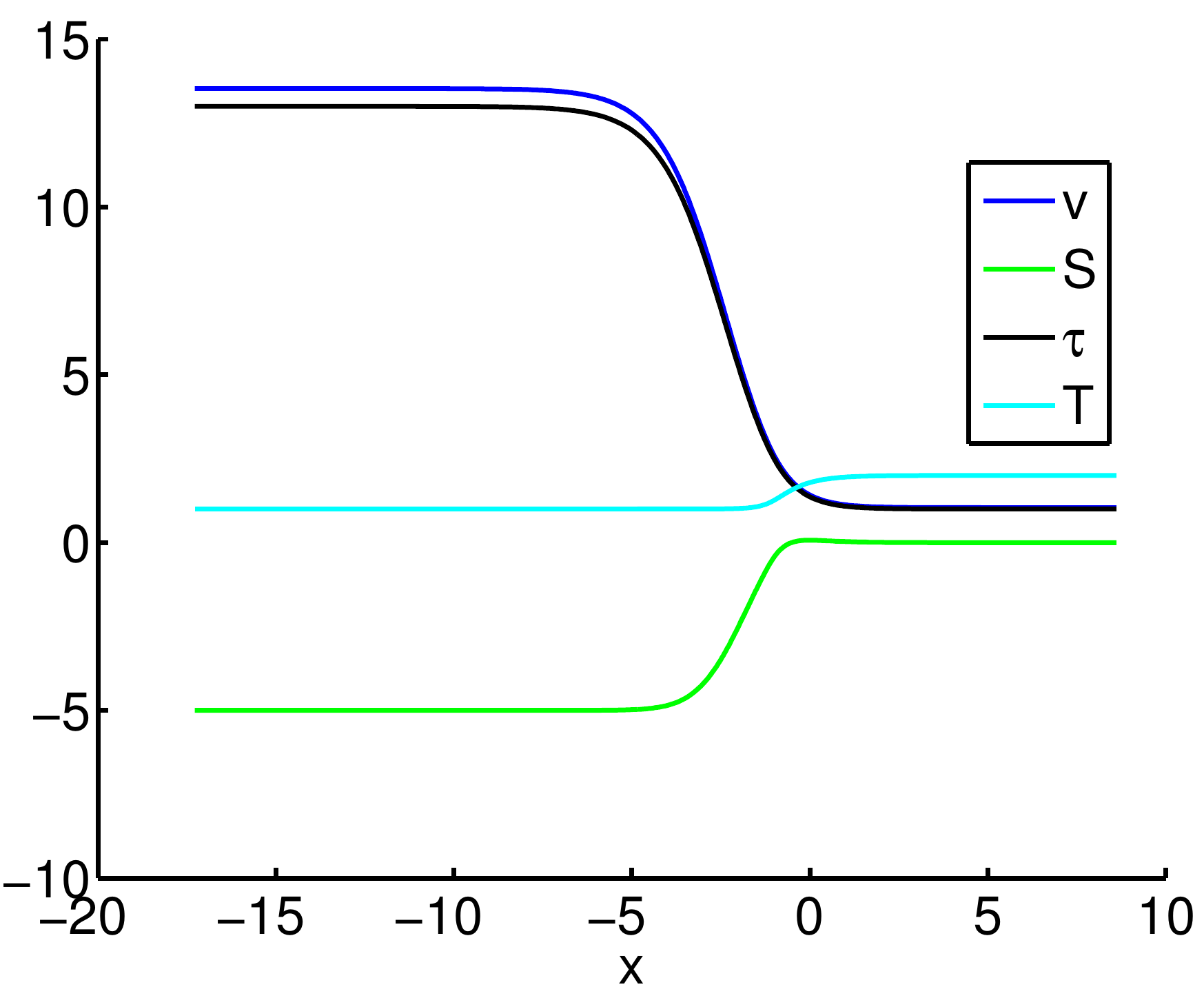}& (b) \includegraphics[scale=0.25]{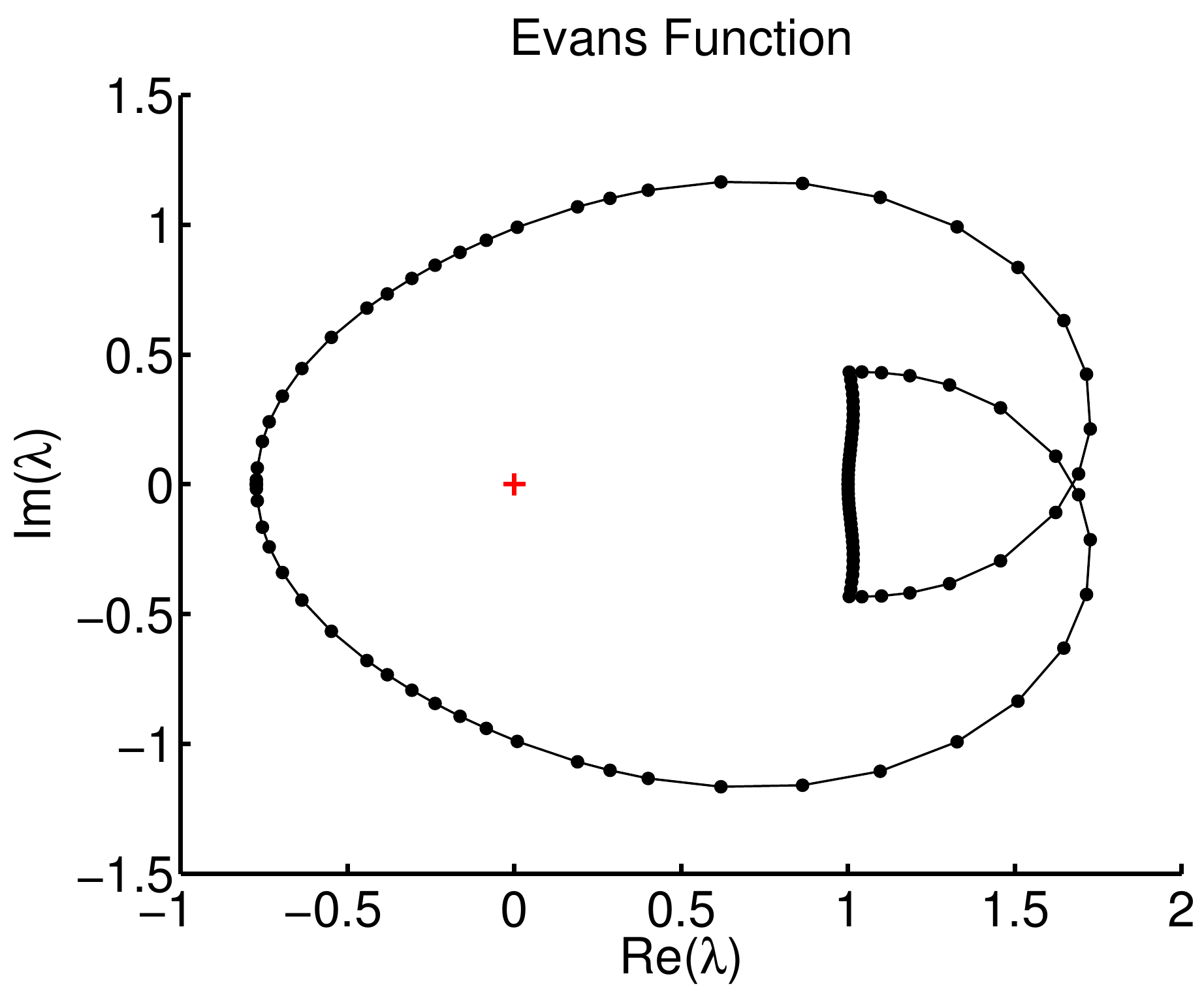}& (c) \includegraphics[scale=0.25]{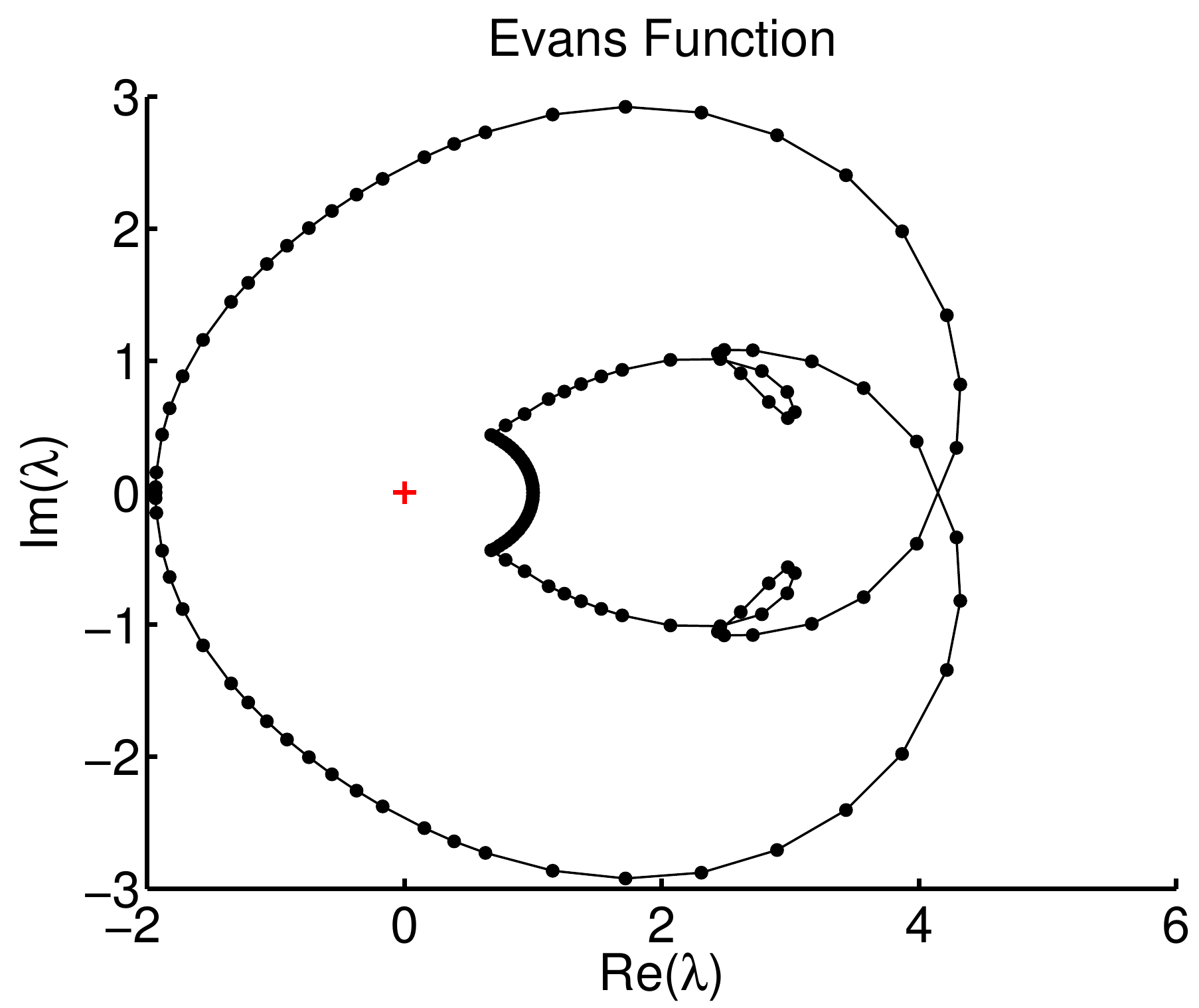}
\end{array}
$
\end{center}
\caption{Viscous study results for $\bar e(S,\tau) = e^S/\tau + S + \tau^2/2$, $S_- = -5$, $\mu=\kappa = 1$, $(\tau_0,S_0)= (1,0)$. (a) Viscous profile. 
(b) Integrated Evans function computed with the adjoint polar coordinate method, evaluated on a semi-circle with radius $R =1$. Winding number is 1. (c) As in (b) but with $R = 250$.
 }
\label{fig666}\end{figure}

%
%

\subsection{Stable example}\label{s:vstable}
Finally, we briefly discuss the case of example \eqref{simplermod}
seen to be inviscidly stable.
In our study we examined the parameter set $S_- \in \{-0.1, -0.2,..., -11.5\}$. For a few values the boundary value solver could not solve the profile equation,
returning a ``singular Jacobian'' error.
For these, we solved the profile using MATLAB's implicit ODE solver, ode15s.
To ensure that
we computed the Evans function on a contour enclosing any possible unstable eigenvalues, we performed a 
high-frequency study curve-fitting with $C_1 e^{C_2\sqrt{\lambda}}$, requiring that
relative error be no greater than 0.2 between the Evans function and the approximating function. 
Table \ref{tb_error3} demonstrates high-frequency convergence for a typical parameter. 
In addition we used an adaptive $\lambda$ mesh to ensure that the change in argument was no greater than 0.2 between successive $\lambda$ points.  We found for all parameters studied that the Evans function has no root in the right half complex plane, consistent with stability. 
The maximum relative 
change
in $\lambda$ for each contour computed was 0.1005. The average value of $R $ was 4 and the average number of points computed was 40 (80 points after using conjugate symmetry of the Evans function to reflect
the values computed on the upper half-plane $\Im \lambda \geq 0$ to the lower half-plane). 
The average time to solve a profile was 12.8 seconds and the average time to solve the Evans function was 6.76 seconds using 8 cores on a Mac Pro. Evans function output for a typical point in parameter space is given in Figure \ref{fig669}.
In all cases, we found a winding number of zero, indicating viscous stability in agreement
with our inviscid findings.


\begin{table}[!b]
\begin{tabular}{|c|c|c|c|}
\hline
 $R$&error& $C_1$&$C_2$\\
\hline 
\hline
2 & 0.04501 & 0.6019 &  0.359 \\
4 &  0.08243 &  0.4554 &  0.3933 \\
 8 & 0.1244 &  0.2953 &  0.4312 \\
16 & 0.147 &  0.1556 &  0.4651 \\
 32 &  0.1416 & 0.06249 &  0.4901 \\
64 & 0.1179 &  0.01757 &  0.5052 \\
128 & 0.09174 & 0.002994 &  0.5136 \\
 256 &  0.06888 &  0.0002471 &  0.5191 \\
 512 & 0.05042 &  7.366e-06 &  0.5223 \\
\hline
\end{tabular}
\caption{For $S_- = -5$ in the stable model we compare how well $C_1 e^{C_2\sqrt{\lambda}}$ approximates $D(\lambda)$ for $|\lambda| = R$. Here $R$ is the radius of the semicircle on which we compute the Evans function, $C_1$ and $C_2$ are the curve fit returned, and the convergence error is the maximum relative error between $D(\lambda)$ and $C_1 e^{C_2\sqrt{\lambda}}$ evaluated at $\lambda = R,\ Ri$.}
\label{tb_error3}
\end{table}


\begin{figure}[htbp]
 \begin{center}
$
\begin{array}{lcr}
\includegraphics[scale=0.25]{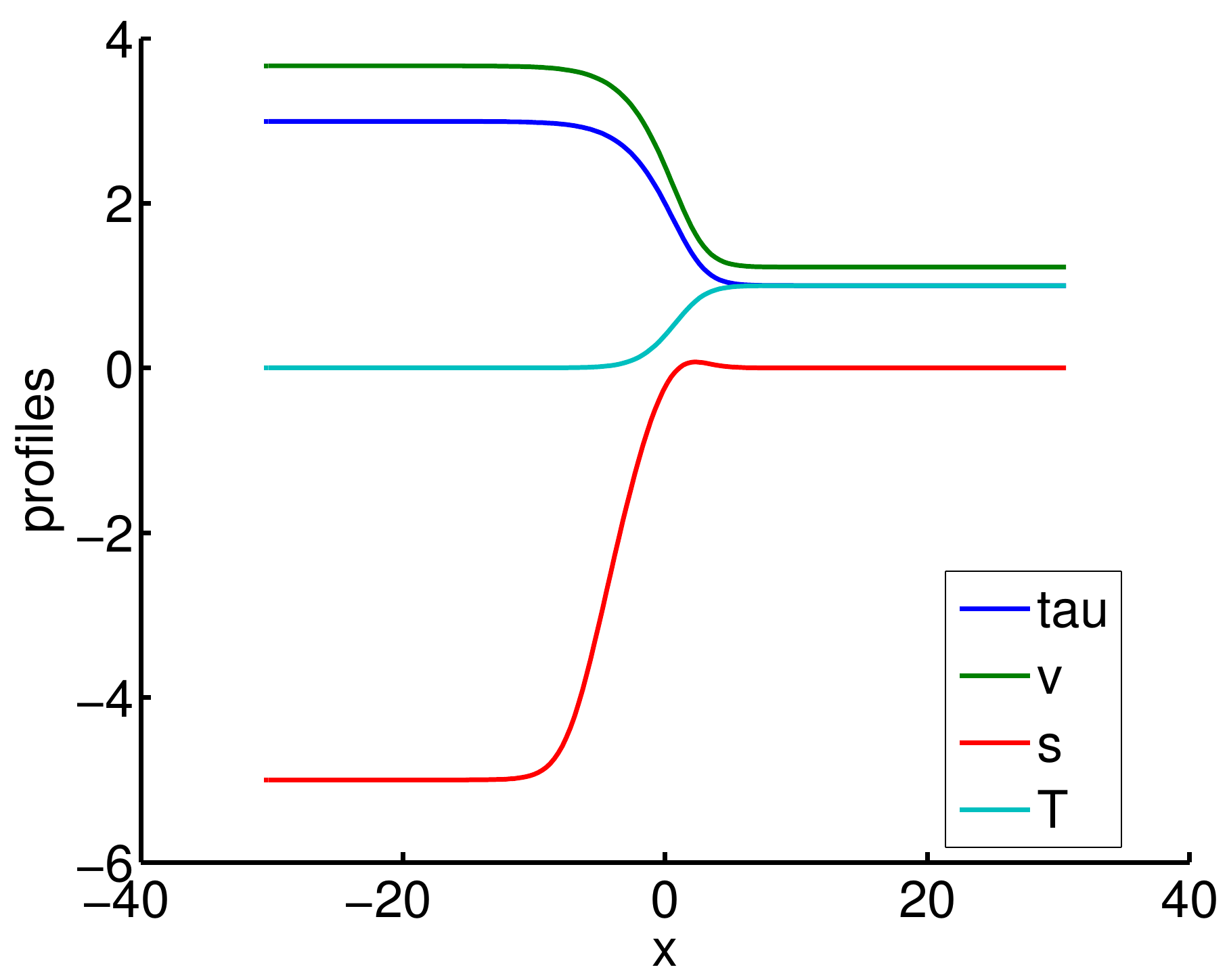}& \includegraphics[scale=0.25]{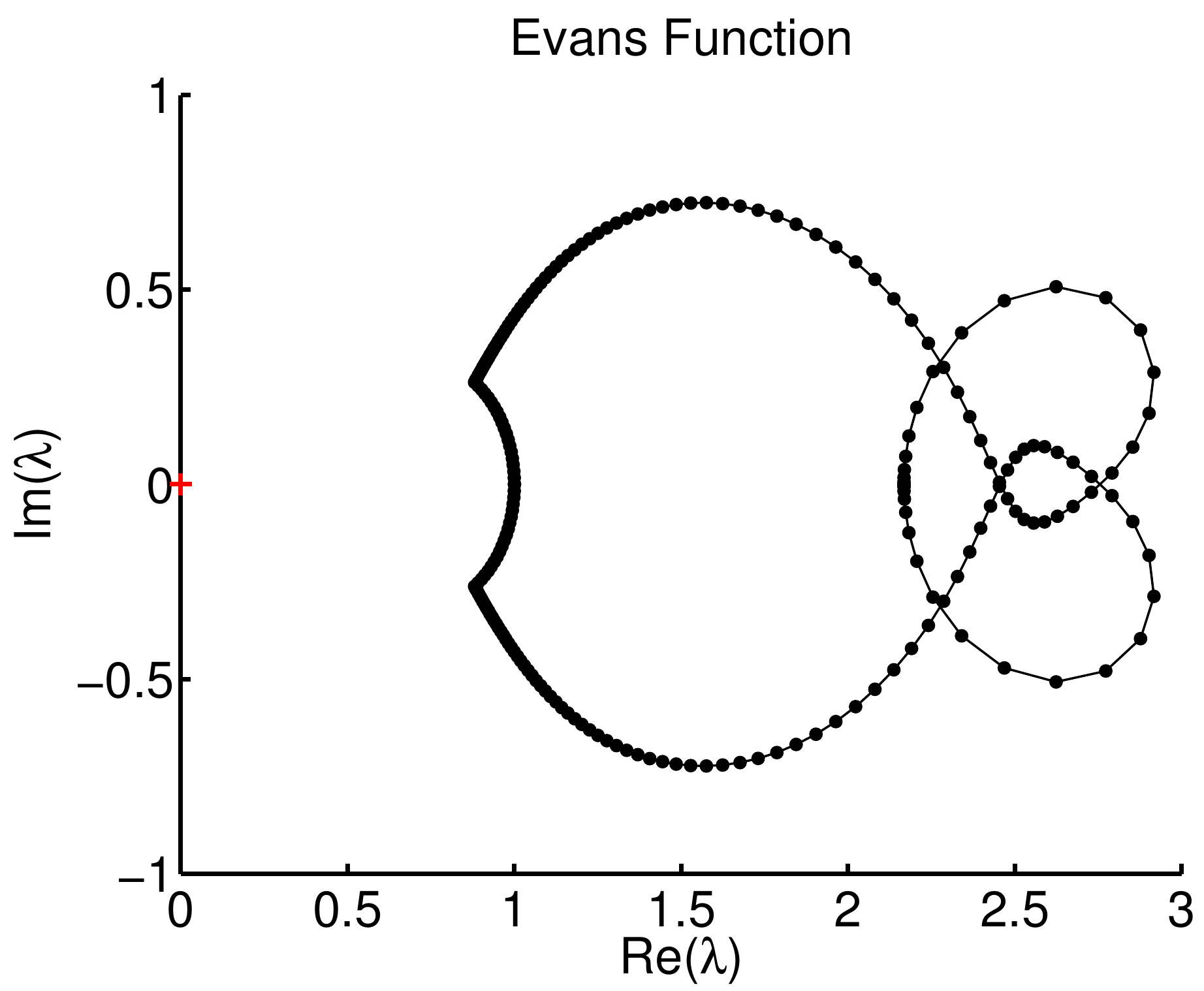} & \includegraphics[scale=0.25]{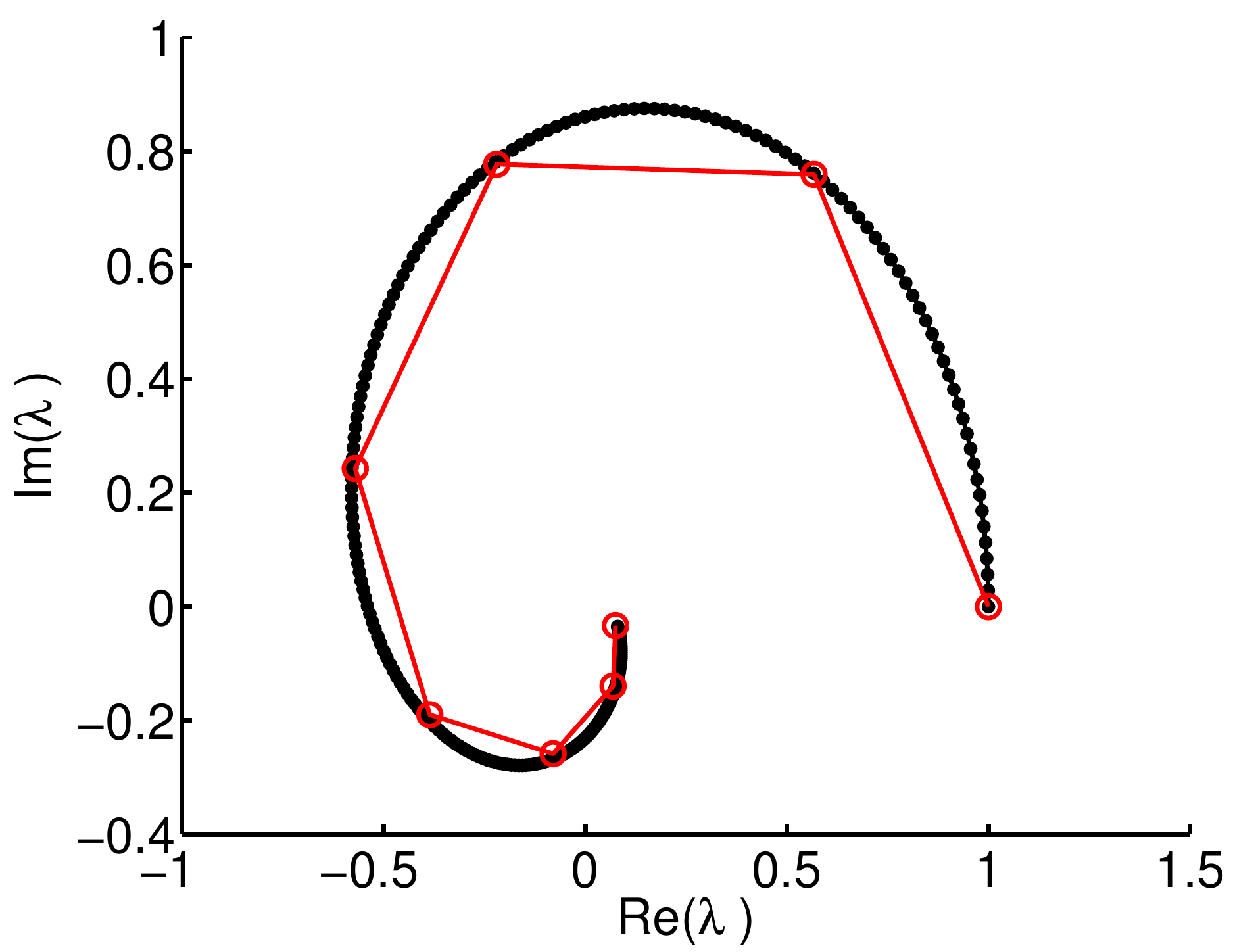}
\end{array}
$
\end{center}
\caption{Viscous study results for $\bar e(S,\tau) = e^S/\tau + \tau^2/2$ with $(\tau_0,S_0) = (1,0)$, and $S_- = -5$. (a) Viscous profile. 
(b)Integrated Evans function computed with the adjoint polar coordinate method, evaluated on a semi-circle with radius $R =250$. Winding number is 1. (c) The Evans function evaluated at 8 points along a semi-circle of radius R = 256 is plotted with open circles. The best curve fit of the Evans function with $C_1 e^{C_2 \sqrt{\lambda}}$ is plotted with closed dots. The relative error between the two is 0.069.
 }
\label{fig669}\end{figure}

\section{Designer systems}\label{s:designer}
We conclude by investigating the range of possible behaviors for general systems
possessing a convex entropy.
Let $A(x)$ be a given real symmetric $n\times n$ matrix-valued function,
and $\mathcal{B}(x)$ a given real positive definite $n\times n$ matrix-valued function,
in the sense that $\Re \mathcal{B}:=\frac12(\mathcal{B}+\mathcal{B}^T)>0$, each converging exponentially
as $x\to \pm \infty$ to limits $\mathcal{A}_\pm$ and $\mathcal{B}_\pm$,
and consider the eigenvalue problem
\be\label{eig}
\lambda u= \partial_x ( \mathcal{B}(x) \partial_x -\mathcal{A}(x)) u,
\quad
u\in \CC^n,
\ee
resembling those found in the stability analysis of viscous shock waves.
This might be considered as a model for the type of eigenvalue
equation that would arise from a symmetrizable system,
in particular a system possessing a convex entropy.

Our first observation is that, by an ingenious construction of 
Bianchini \cite{B}, {\it any} such system can in fact be realized as the
eigenvalue problem for a viscous shock wave solution of a
system of conservation laws with a viscosity-compatible strictly convex
entropy.
Specifically, consider the system
\ba\label{2bsys}
u_t + \big(A(v)u\big)_x&= \big(B(v)u_x\big)_x,\\
v_t + \big(\frac12 u\cdot A'(v)u + \frac12 v^2\big)_x&=v_{xx},
\ea
where $v\in \RM^1$.

This contains a trivial standing viscous shock solution
$
(\bar u, \bar v)=(0,\bar v(x)),
$
where $\bar v(x):=-\tanh(x/2)$ is a shock of Burgers equation
$v_t + (v^2/2)_x=v_{xx}$.
Linearizing \eqref{2bsys} about $(\bar u,\bar v)$,
we obtain
\ba\label{full_lin}
\bp u\\v\ep_t + \bp A(\bar v)& 0\\
0& \bar v  \ep 
\bp u\\v\ep_x=
 \Big( \bp B(\bar v) & 0\\ 0& 1 \ep
 \bp u \\ v\ep_x\Big)_x\, ,
\ea
which decouples into a (stable) linearized Burgers equation in $v$ and
an equation 
\be\label{ulin}
u_t + (A(\bar v(x))u)_x=(B(\bar v(x))u_x)_x,
\ee
that has the associated eigenvalue equation \eqref{eig} 
provided that we choose 
$ A(v):= \mathcal{A} (2\arctan (-v))$
and
$ B(v):= \mathcal{B} (2\arctan (-v))$.
%
Moreover, it is easily verified that \eqref{2bsys} possesses
a viscosity-compatible stricly convex entropy
$
\eta(u,v):= \frac12|u|^2 +\frac12 |v|^2,
$
with associated entropy flux
$q(u,v):=  \frac12\langle u, A((v)u\rangle 
+ \frac13|v|^3 + \frac12 v \langle u, A'(v) u.$
That is, $d\eta dF=dq$ 
and $d^2\eta \bp B&0\\0&1\ep= \bp B&0\\0&1\ep $ is positive definite,
where 
$
F(u,v)=\bp A(v)u \\ \frac12 u\cdot A'(v)u +\frac12 v^2\ep.
$
This means that we may search for stability phenomena of
symmetrizable systems 
possessing a viscosity-compatible strictly convex entropy, 
by considering the linear problem \eqref{eig}, which is free
to our specifications.

\subsection{The rotating model}\label{s:rot}
A natural choice in seeking examples of instability is the 
{\it rotating model} consisting of \eqref{2bsys} with
\be\label{Adef}
A(v)=  R_{\theta(v)} A_m R_{-\theta(v)},
\quad
B(v)\equiv \Id,
\ee
where 
$R_\theta:=\bp \cos \theta & -\sin \theta\\
\sin \theta & \cos \theta\ep,$
$A_m=\bp 1&0\\0&-1\ep,$
and 
\be\label{theta}
\theta(v)= M\pi v
\ee
$M$ an arbitrary real number,
considering the family of stationary shocks
$\bar v_\gamma(x):=-\gamma \tanh (\gamma x/2)$.

For, it is readily seen that these are Lax $2$-shocks, with
a Lopatinski determinant 
$$
\delta_{\gamma,M}=\det \bp \cos{-M\pi \gamma} & -\sin(M\pi \gamma) & 0\\
\sin{-M\pi\gamma} & \cos( M\pi \gamma) & 0\\
0&0 &-2\gamma\ep=
-2\gamma (\cos^2 (M\pi\gamma) -\sin^2(M\pi \gamma))
$$
that changes sign from negative to positive and back infinitely
often as $\gamma$ increases from the small-amplitude limit $\gamma=0$,
with the first stability transition occuring at $M\pi \gamma_*=\pi/4$,
or $\gamma_*= 1/4M$.
This in passing gives another example of a $3\times 3$ system 
with convex entropy exhibiting inviscidly unstable shock waves, 

\subsubsection{Evans system}\label{s:Evans}

For 
profile $\bar u\equiv 0$, $\bar v=-\gamma\tanh(\gamma x/2)$,
the integrated Evans system for the decoupled, $u$ equation \eqref{ulin}
is
\be\label{e:evans}
W'=\mathcal{A}(x,\lambda)W,
\qquad
\mathcal{A}=
\begin{pmatrix}
0&I\\\lambda B^{-1}& B^{-1}A
\end{pmatrix}(\bar v(x)).
\ee

\br\label{explicit}
Setting $W=TY$, where 
$
T:=\bp R_\theta & 0\\
0 & R_\theta \ep,
$
we may convert \eqref{e:evans} to $Y'=\mB Y$, where
\be\label{mB}
\mB=
 T^{-1}\mA T - T^{-1}T'   
=
\bp 0 & I\\
\lambda & A_m\ep
-
M\pi \bar v'(x) \bp J & 0\\ 0 & J\ep,
\ee
and $J:=\bp 0 & -1\\1& 0\ep$.
Involving only exponential functions in $x$ entering through
the scalar multiplier $\bar v'$, this seems possibly
amenable to exact solution, an interesting direction for further study.
\er

\subsubsection{Transversality of profiles}\label{s:dindex}
Again appealing to the triangular form of the equations, we find that
transversality of profiles is equivalent to nonexistence of asymptotically
decaying solutions of the $u$-component 
$ u'=A(\bar v(x))u$ of the linearized standing-wave equation.
More, the transversality coefficient $\nu$, like the Lopatinski
determinant, factors into $-2\gamma$ times a transversality coefficient
for this decoupled $u$-component.

\bl\label{dproflem}
For $|M\gamma|$ bounded, profiles are transversal for $\gamma$ 
sufficiently small.
More, when appropriately normalized, 
(i) $\nu/(-2\gamma) \to 1$ as $\gamma\to 0^+$ and (ii)
$\nu/(-2\gamma)\to \delta/(-2\gamma)$ as $\gamma\to +\infty$,
while (iii) $\nu/(-2\gamma)$ changes sign at least $[4(M\gamma -1)]$
times as $\gamma$ increases from $0^+$ to $+\infty$,
where $[\cdot]$ denotes the greatest integer function.
\el

\begin{proof}
(i) (Basic tracking argument; see \cite{Z5}.)
Setting $u=R_\theta z$, similarly as in Remark \ref{explicit},
we obtain 
$$
\hbox{\rm $z'=A_mz + M\bar v'(x)J$, with $J:=\bp 0 & -1\\1& 0\ep$.}
$$
Observing that $|M\bar v'(x)|\le M\gamma^2\le C\gamma\to 0$ for
$\gamma\to 0$ and $|M\gamma|\leq C$, we obtain the result
by the Tracking Lemma of \cite{ZH}, a standard estimate for
slowly-varying-coefficient systems, after factoring out the (real)
exponential growth rates in growing and decaying modes.

(ii) Noting, for \eqref{Adef}--\eqref{theta}
that Evans system \eqref{e:evans} satisfies 
$|\mathcal{A}(x)-\mathcal{A}_\pm|\le CM\gamma e^{-\gamma|x|}
\leq C_2e^{-\gamma|x|}$
for $x\gtrless 0$, we obtain by the Convergence Lemma of \cite{PZ}
that the Wronskian $\nu$ converges up to an error of order
$\|C_2e^{-\gamma|x|}\|_{L^1}=O(C/\gamma)\to 0$
to the determinant of a (smoothly chosen real) stable
eigenvector of $\mathcal{A}_+$ and a (smoothly chosen real) 
unstable eigenvector of $\mathcal{A}_- $, i.e., $\delta/(-2\gamma)$.

(iii) Reviewing the previous arguments more closely, we find that for $\gamma$
sufficiently small the decaying
mode as $x\to -\infty$ rotates angle $2M\gamma \pi$, while for $\gamma$
sufficiently large it rotates at most $\pi$, with the directions of
stable and unstable subspaces as $x\to \pm \infty$ held
fixed for all $\gamma$.  Thus, as $\gamma$ goes from $0$ to $+ \infty$,
the decaying mode as $x\to \infty$ crosses the decaying mode as $x\to +\infty$
at least $[4( M\gamma -1)]$ times, each time 
signaling an associated change in sign.
\end{proof}

From Lemma \ref{dproflem}, we find that, fixing $K$ sufficiently
large, taking $\gamma_*$ sufficiently small and $M:=K/\gamma_*$,
and varying $0\leq \gamma \leq \gamma_* $, the profiles $\bar v_\gamma$
remain transversal, while the Lopatinski determinant changes sign 
$\approx 2K$ times, signaling (recall \eqref{ZSint}) $2K$ passages 
of an eigenvalue through the
imaginary axis, hence appearance of up to $2K$ unstable roots.
(Our numerical experiments indicate that these roots all cross in
one direction, similarly as in standard Sturm-Liouville theory.)

On the other hand, fixing $M\gamma=\mathrm{constant}$, and 
increasing $\gamma$ from $0^+$, we find that $\delta$ stays constant,
while $\nu$ changes sign at least $\approx 4(M\gamma -1)$ times. 
More, varying $M\gamma$ within a bounded set for $\gamma$ sufficiently large, 
we find that changes of sign in $\delta$ and $\nu$ may be made to occur arbitrarily close together, increasing the chance of collision of roots and subsequent splitting into a complex conjugate pair.
(Recall that changes of sign of $\delta$ and $\nu$ signal crossings through the origin of
roots of the Evans function $D$.)
This simple rule of thumb guides our strategy in searching for complex roots and Hopf bifurcation,
and indeed gives a reasonable a posteriori fit to the data that we see.
Why some colliding roots split into the complex plane while others appear to stay real 
is an interesting question to which we do not have an answer.

\subsubsection{High-frequency bounds}\label{s:dhf}
We start by bounding the modulus of unstable eigenvalues.
\bl\label{deslem}
For system \eqref{full_lin} with the choices \eqref{Adef}--\eqref{theta}, 
there are no unstable eigenvalues with modulus 
$
|\lambda|\geq 
4.
$
\el

\begin{proof}

Let $A$ be as in  \ref{Adef} and consider the integrated eigenvalue problem,
$
\lambda u + A u'  = u''.
$
Multiplying by the conjugate $\bar u$ and integrating over $\R$ yields
\be
\begin{split}
\lambda ||u_j||^2 + \int_{\R} a_{j1} \bar u_j u_1 'dx + \int_{\R} a_{j2} \bar u_j u_2'dx&= \int_{\R} \bar u_1 u_1''dx, \quad j = 1,2,
\end{split}
\label{bound1}
\ee
where $a_{ji}$ is the $j$-$i$ entry of $A$.
Integrating the 
righthand side by parts and rearranging gives
\be
\begin{split}
\lambda ||u_j||^2 + ||u_j'||^2&=- \int_{\R} a_{j1} \bar u_j u_1 ' dx- \int_{\R} a_{j2} \bar u_j u_2'dx, j = 1,2.
\end{split}
\label{tool2}
\ee
Taking the real part of \eqref{tool2} for $\Re(\lambda) \geq 0$, and applying Cauchy's inequality, we have 
\be
\begin{split}
\Re(\lambda) ||u_j||^2 + ||u_j'||^2 & \leq C\big( \int_{\R} | u_j | |u_1 '| dx +  \int_{\R} | u_j|| u_2'|dx\big)
\leq C\big(  2\eps ||u_j||^2+ \frac{1}{4\eps}||u_1'||^2  + \frac{1}{4\eps} ||u_2'||^2\big),
\end{split}
\label{tool3}
\ee
for $j = 1,2$ where (using 
$\sin 2\theta=2\cos \theta \sin \theta$) $C = \max_{x\in \R}{|a_{j,i}(x)|}=1$. Similarly,
\be
\begin{split}
|\Im(\lambda)| ||u_j||^2 
&\leq C\big(  2\eps ||u_j||^2+ \frac{1}{4\eps}||u_1'||^2  + \frac{1}{4\eps} ||u_2'||^2\big),
\end{split}
\label{tool4}
\ee
Summing \eqref{tool3} and \eqref{tool4} for $j = 1,2$ we have
$$
(\Re(\lambda)+|\Im(\lambda)|)( ||u_1||^2+||u_2||^2) + ||u_1'||^2+||u_2'||^2 
\leq C\big(  4\eps ||u_1||^2+4\eps||u_2||^2+ \frac{1}{\eps}||u_1'||^2  + \frac{1}{\eps} ||u_2'||^2\big).
$$
Taking $\eps = C$, we obtain
$(\Re(\lambda)+|\Im(\lambda)|)( ||u_1||^2+||u_2||^2)
\leq 4C^2(||u_1||^2+||u_2||^2)$,
so that for $\Re(\lambda) \geq 0$, $|\lambda| \leq 4C^2 \leq 4$.
\end{proof}

\subsection{Numerical results}\label{s:dnum}
Here we describe our numerical studies for the designer system \eqref{2bsys}.

\subsubsection{Comparison of $D(0)$, $\nu$, $\delta$}
For system \eqref{e:evans} with $\lambda=0$, explicit orthonormal initializing bases $R_{\pm}$ 
for the Evans function, i.e., bases for the stable (unstable) subspaces of 
$\mathcal{A}(-\infty)$ ($\mathcal{A}(_\infty)$), are given by
\begin{small}
\besn{   
R_- = \left\{\frac{1}{2}\mat{\cos(M\gamma \pi)\\\sin(M\gamma \pi)\\\cos(M\gamma \pi)\\\sin(M\gamma \pi)},\mat{-\sin(M\gamma \pi)\\\cos(M\gamma \pi)\\0\\0}\right\}, \quad R_+ = \left\{\frac{1}{2}\mat{\sin(M\gamma \pi)\\\cos(M\gamma \pi)\\-\sin(M\gamma \pi)\\-\cos(M\gamma \pi)},\mat{\cos(M\gamma \pi)\\-\sin(M\gamma \pi)\\0\\0}\right\}.
}
\end{small}
Likewise, initializing bases $r_\pm$ for the linearized traveling-wave ODE
$w'=Aw$ of the unstable (stable) subspaces of $A(\pm \infty)$ are given by
$r_- = \frac{1}{2}\mat{ \sin(M\gamma \pi)}$,
and $r_+ = \frac{1}{2}\mat{ -\sin(M\gamma \pi)\\-\cos(M\gamma \pi)}$.
We compute $D(0)$ and $\nu$ with these choices, for which \eqref{ZSint} becomes essentially
a tautology.\footnote{Note that vectors $(*,*,0,0)^T$ are stationary solutions of
\eqref{e:evans}, so that the Wronskian of the associated modes simplifies to $\delta$
times a $2\times 2$ block involving third and fourth components, i.e., the Wronskian $\nu$
of $w'=Aw$.}

Figure \ref{fig662} demonstrates the effects of the rotation in solving the Evans function ODE.

\begin{figure}[htbp]
 \begin{center}
$
\begin{array}{lccr}
(a) \includegraphics[scale=0.15]{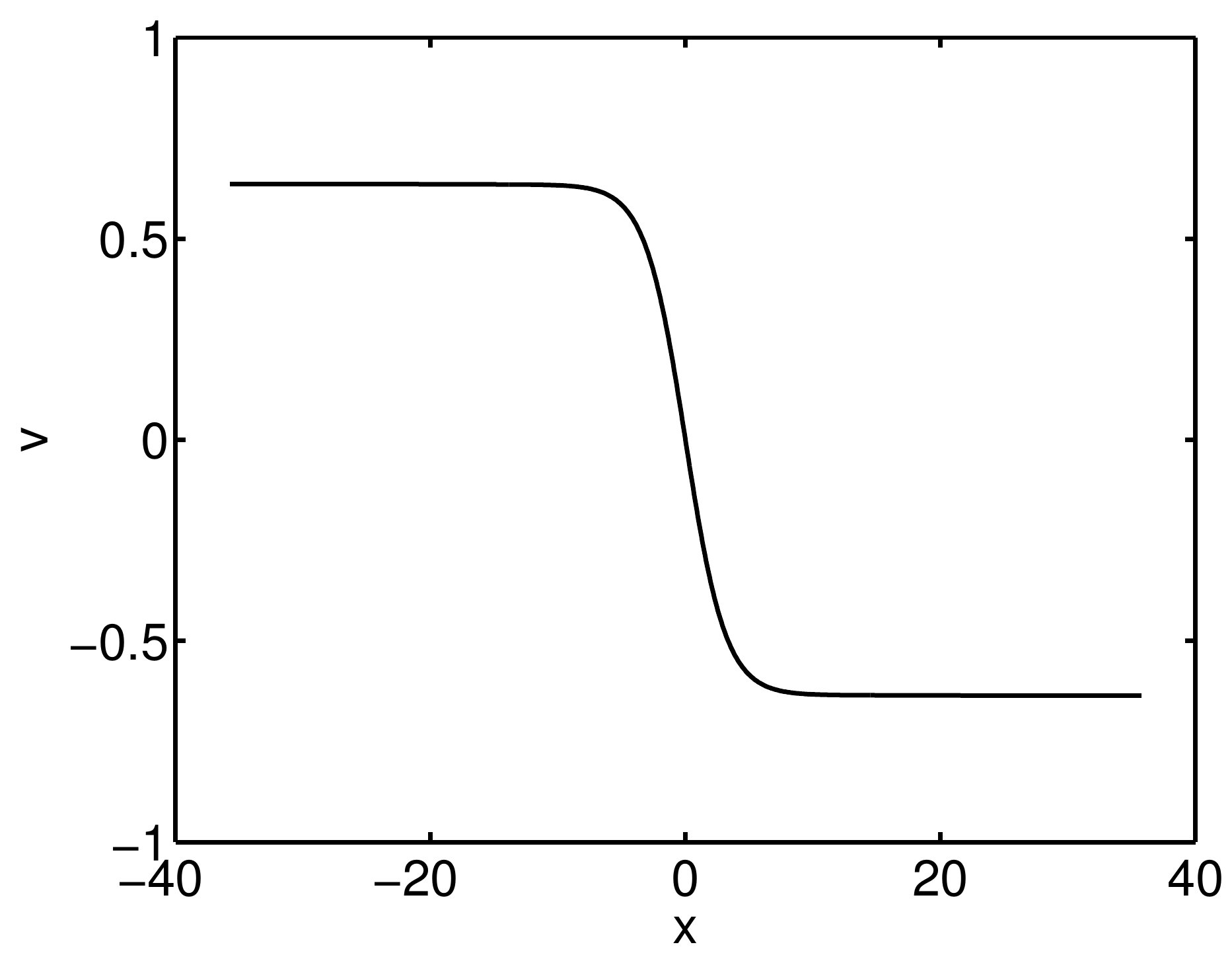}& (b) \includegraphics[scale=0.15]{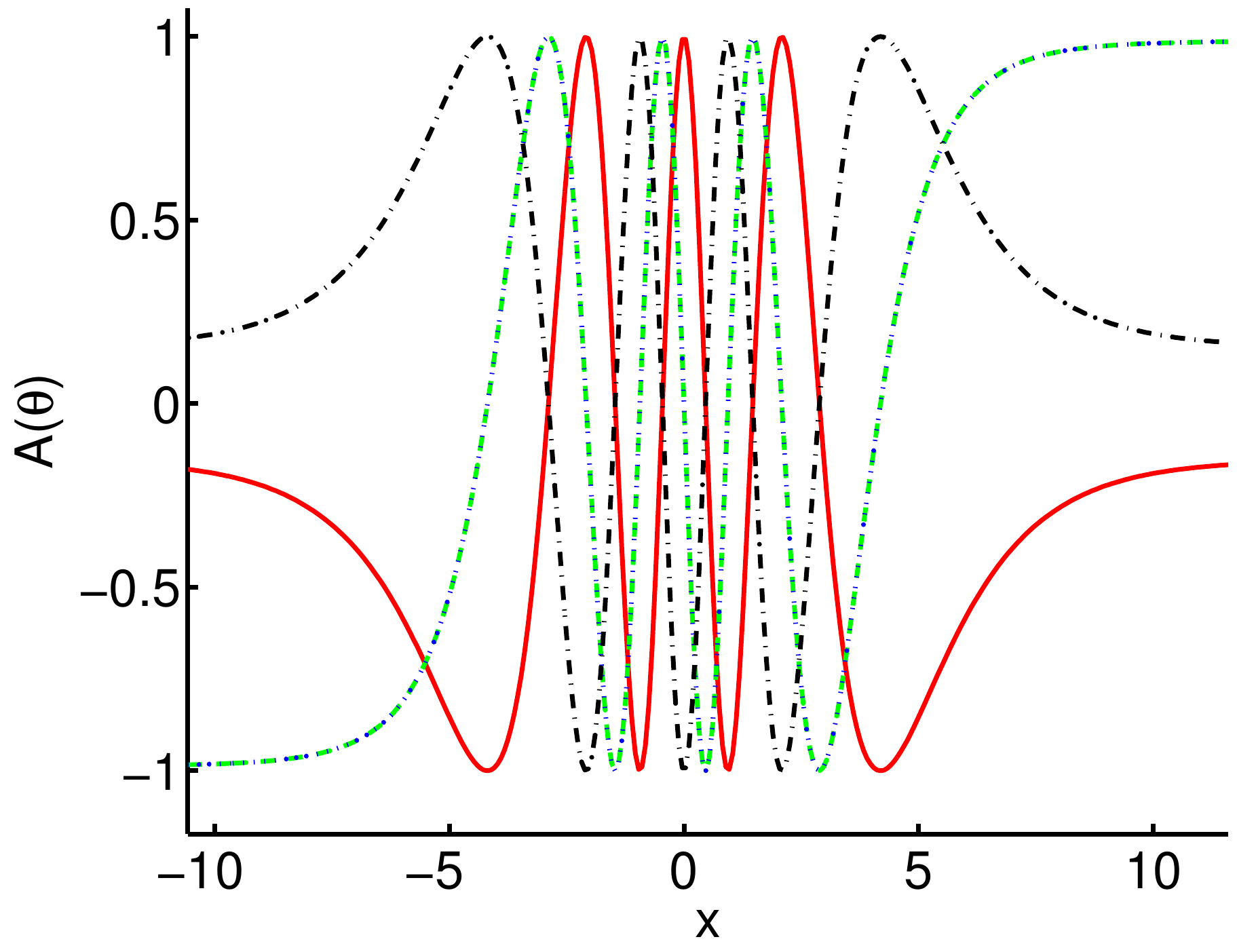}& (c) \includegraphics[scale=0.15]{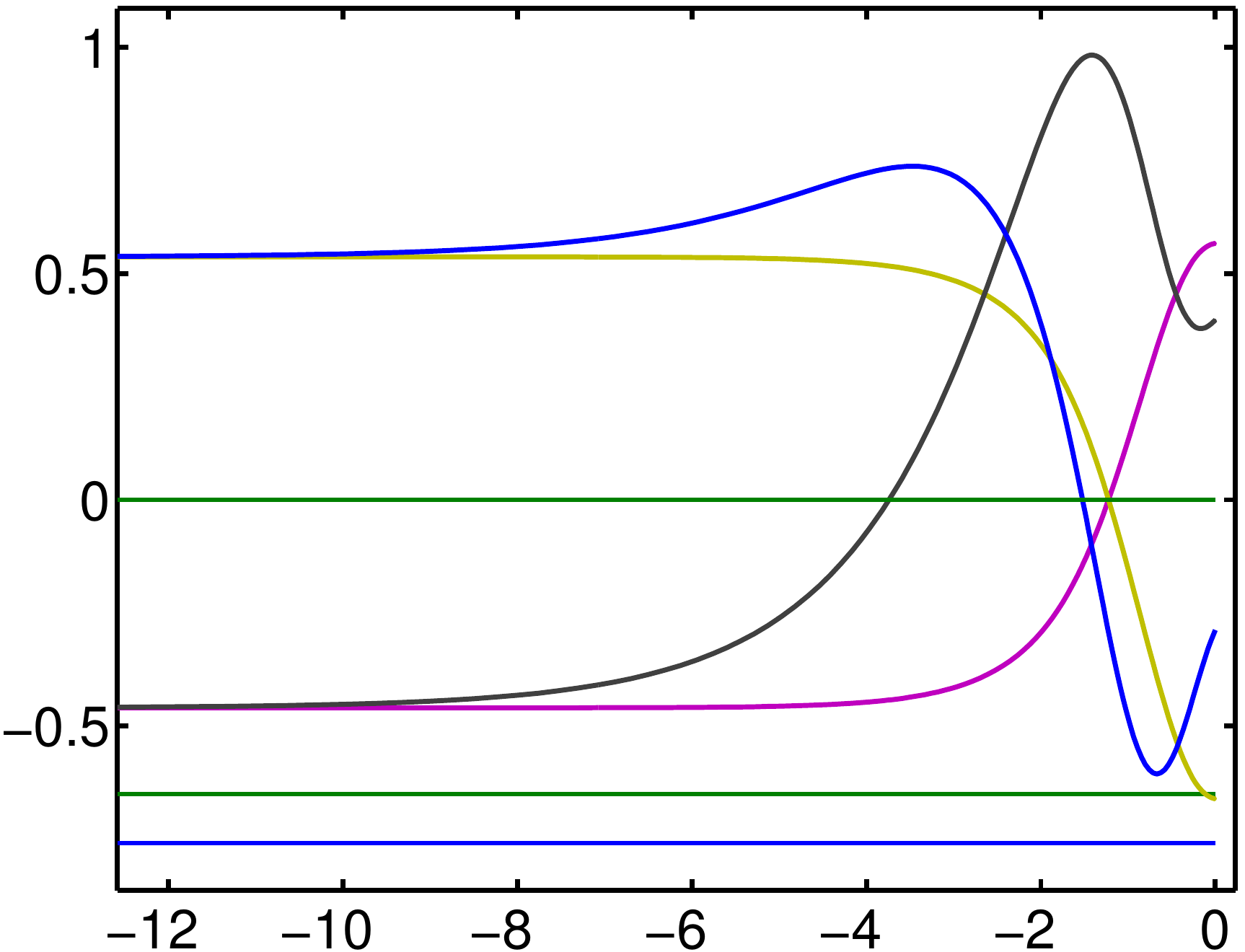}& (d) \includegraphics[scale=0.15]{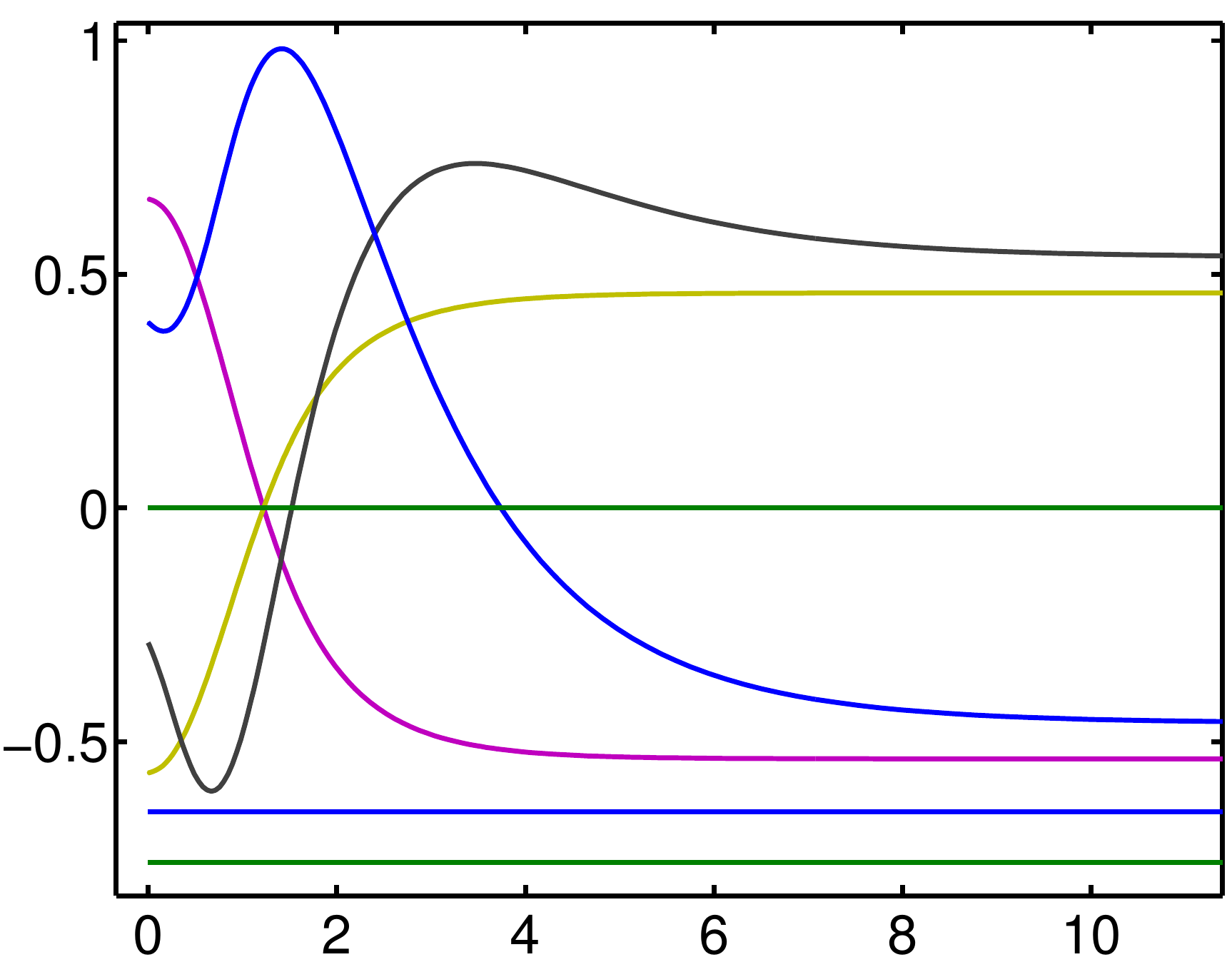}
\end{array}
$
\end{center}
\caption{Rotating model with $M = 2.7174$, $\gamma  = 0.635$. (a) Plot of profile $\bar v_{\gamma}(x) = -\gamma \tanh(\gamma x/2)$ against $x$. (b) Zoomed in plot of components of $A(\theta(x))$ against $x$. (c)-(d) Zoomed in plot of components of the evolved manifolds in the Evans function computations for $\lambda = 27/2$. 
}
\label{fig662}\end{figure}

\begin{figure}[htbp]
 \begin{center}
$
\begin{array}{lccr}
(a) \includegraphics[scale=0.18]{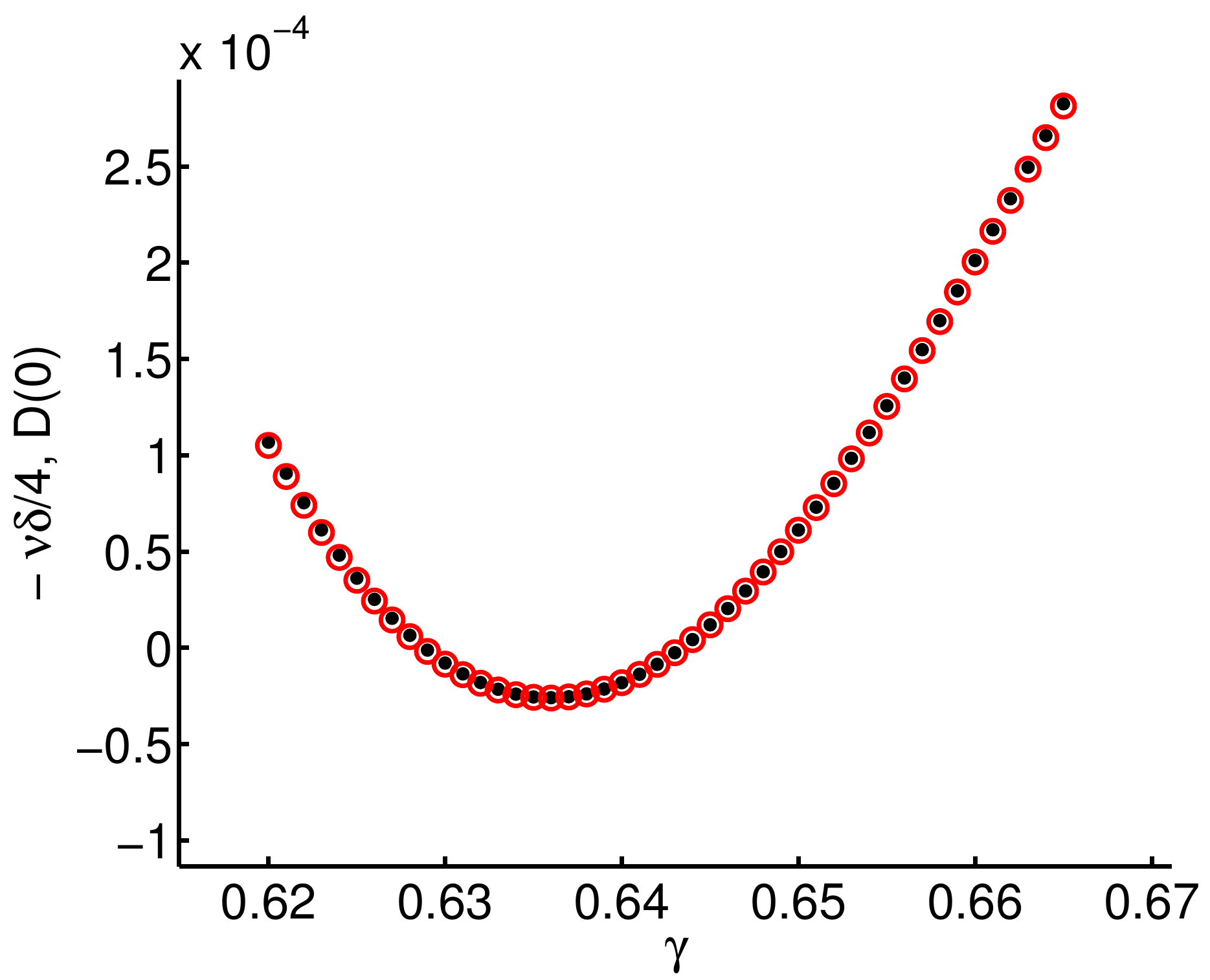}&(b) \includegraphics[scale=0.18]{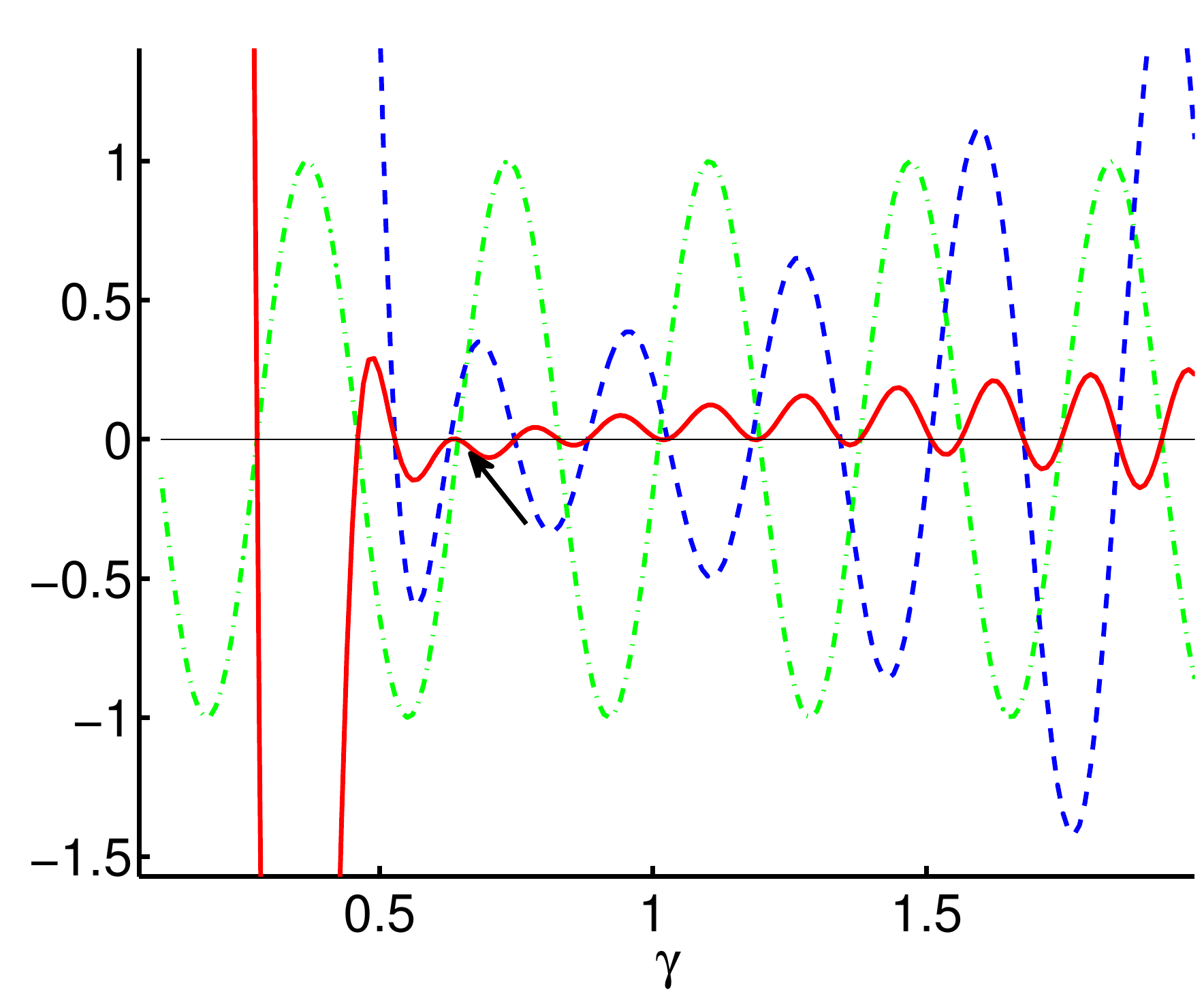} & (c) \includegraphics[scale=0.18]{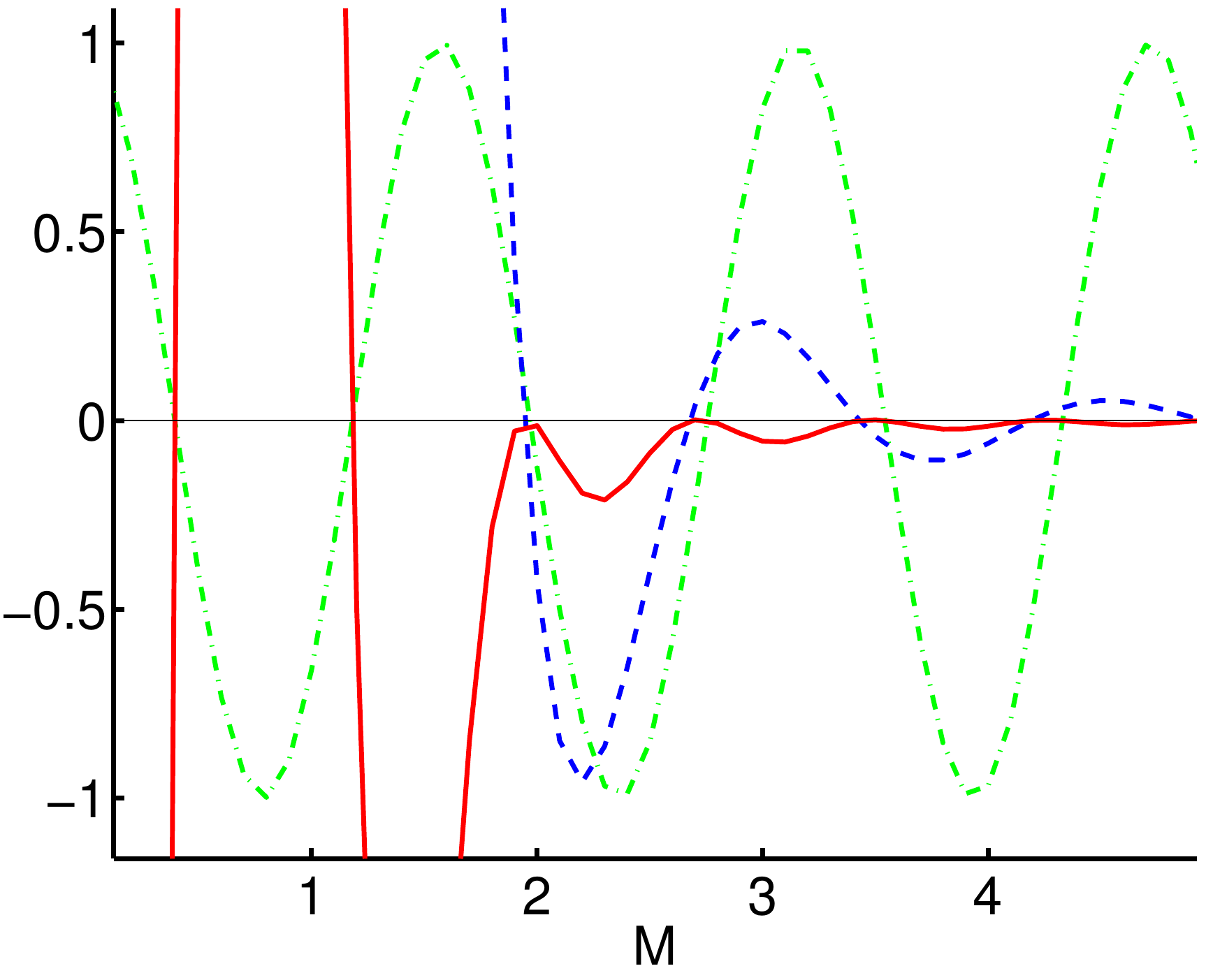}&(d) \includegraphics[scale=0.18]{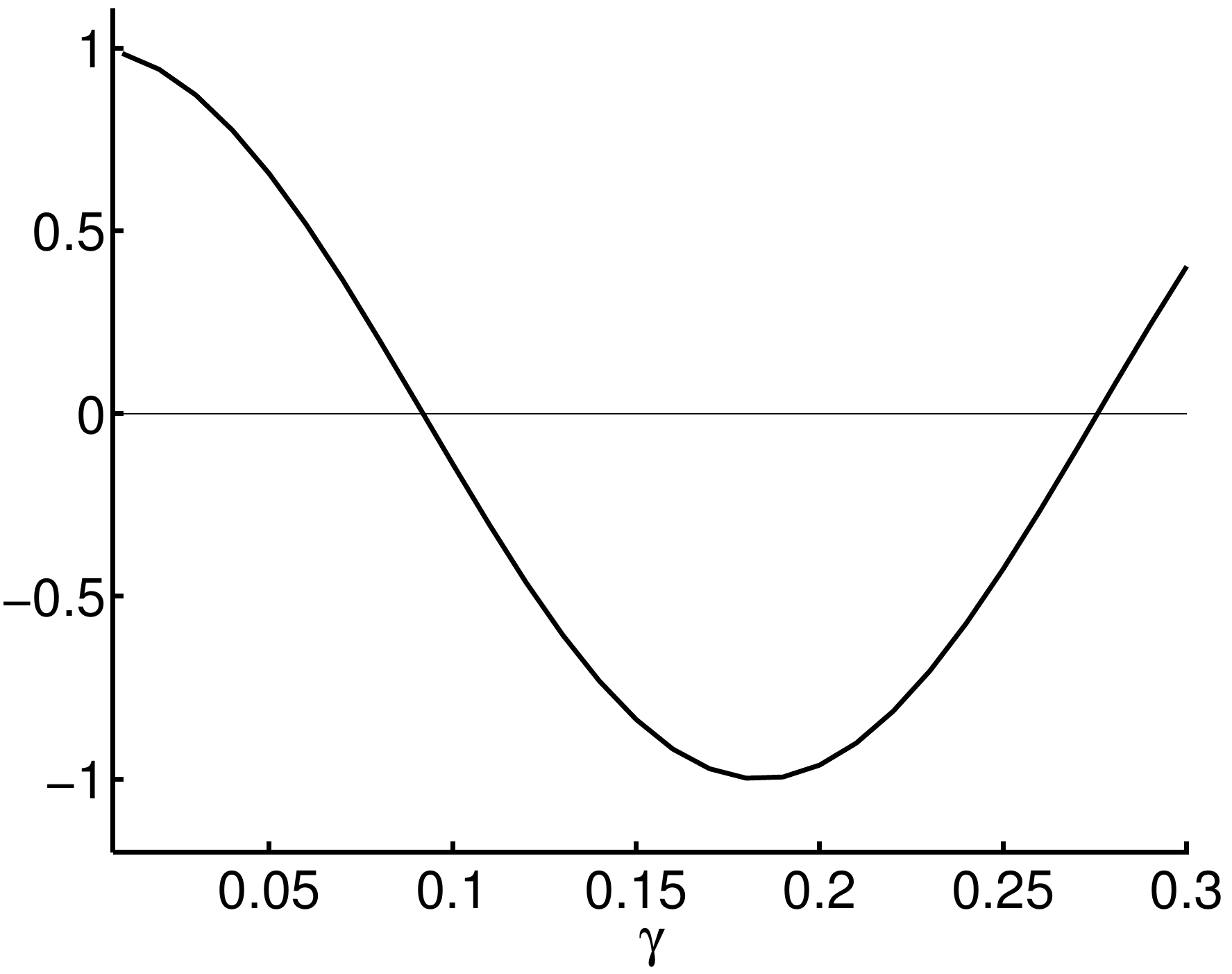} 
\end{array}
$
\end{center}
\caption{Rotating model. (a) Plot of $D(0)$ (solid black dots) and $-\delta \nu/4$ (red circles) against $\gamma$ for $M = 2.72$ demonstrating the identity $D(0) =- \nu \delta/4$. (b) Plot of $D(0)$ (solid red), $\delta$ (dot-dashed green), and $\nu$ (dashed blue) against $\gamma$ for $M = 2.72$. We mark with an arrow the general region where the Hopf bifurcation occurs. (c) Plot of $D(0)$ (solid red), $\delta$ (dot-dashed green), and $\nu$ (dashed blue) against $M$ with $\gamma = 0.635$. (d) Plot of $\delta(\gamma)$ for $M = 2.72$. Note that the first root of $\delta$ occurs at $\delta_* = 1/(4M) \approx 0.09$. }
\label{fig763}\end{figure}


\subsubsection{Evans function computations}
We now detail the results of our numerical study for the designer model. For all our Evans function computations we used the integrated form \ref{e:evans} of the eigenvalue problem. We used MATLABs implicit ODE solver ode15s with relative and absolute error bounds set respectively to 1e-6 and 1e-8. We used winding number computations and a divide and conquer scheme  described in \ref{s:num} to find roots of the Evans function. To obtain a global picture of roots, we 
used Lemma \ref{deslem}
to ensure that the radius of the contour on which the Evans function is computed
 was sufficiently large to enclose all possible unstable eigenvalues. 
We also sometimes used a smaller radius for speed when looking just for the existence of 
Hopf bifurcations. 

To search for complex roots and or Hopf bifurcation, we held fixed one of the parameters
$\gamma$, $M$, and $M\gamma$ while letting a second parameter vary, examining the resulting
motion of eigenvalues along these curves.
The results of such experiments can be fairly complicated, and can involve a large number of roots
(up to $16$ in the experiment of Figure \ref{fig814}(b)).
However, most of these roots are real for most of the time, and, when they do collide and split
into a complex pair, tend quickly to rejoin and become again real.
Thus, to find complex roots, and in particular to identify Hopf bifurcations (roots that are not
only complex, but cross the imaginary axis as a conjugate pair) involves
finding a rather
small-measure subset of parameter space, and requires a systematic search;
they were not easy to find!
An approach that proved quite useful in identifying Hopf bifurcations was to count the number of unstable
roots at a selection of mesh points in the $\gamma$-$M\gamma$ plane and look for jumps of
$2$ across boundaries of different regions in the number of roots, then check with a refined
study that (as generically should hold) these correspond to crossing of the imaginary
axis by a non-real conjugate pair.

We note in Figure \ref{fig763} that for $M = 2.72$ the first root of $\delta(\gamma)$ occurs at $\delta_*= 1/(4M)\approx 0.09$. 
Figure \ref{fig814} (a) demonstrates that Hopf bifurcations exist and occur roughly periodically in
the parameter $M\gamma$, with $\gamma$ near special values.
We used $R = 16$ for the contour radius and checked that
 $\delta(\gamma)$ had no roots in the vicinity of the jumps by 2 to verify that indeed they correspond to Hopf bifurcations. 
Figure \ref{fig814} (b) shows that for $M\gamma$ fixed the number of Hopf bifurcations occurring increases as $\gamma \to 0$. 
Figure \ref{fig814} (c) shows a 
Hopf bifurcation occurring as $M$ increases, for $\gamma $ fixed.
Figure \ref{fig814b} (c) shows a 
Hopf bifurcation occurring as $\gamma$ increases, for $M$ fixed.
(The latter corresponds to variation in shock amplitude for a fixed set of 
equations, hence is our true interest.) 
Figure \ref{fig814b} (d) shows crossing of the two eigenvalues with largest
real part, with a double-multiplicity eigenvalue occurring at the point of crossing, violating simplicity of the top eigenvalue.

\begin{figure}[htbp]
 \begin{center}
$
\begin{array}{lcr}
(a) \includegraphics[scale=0.25]{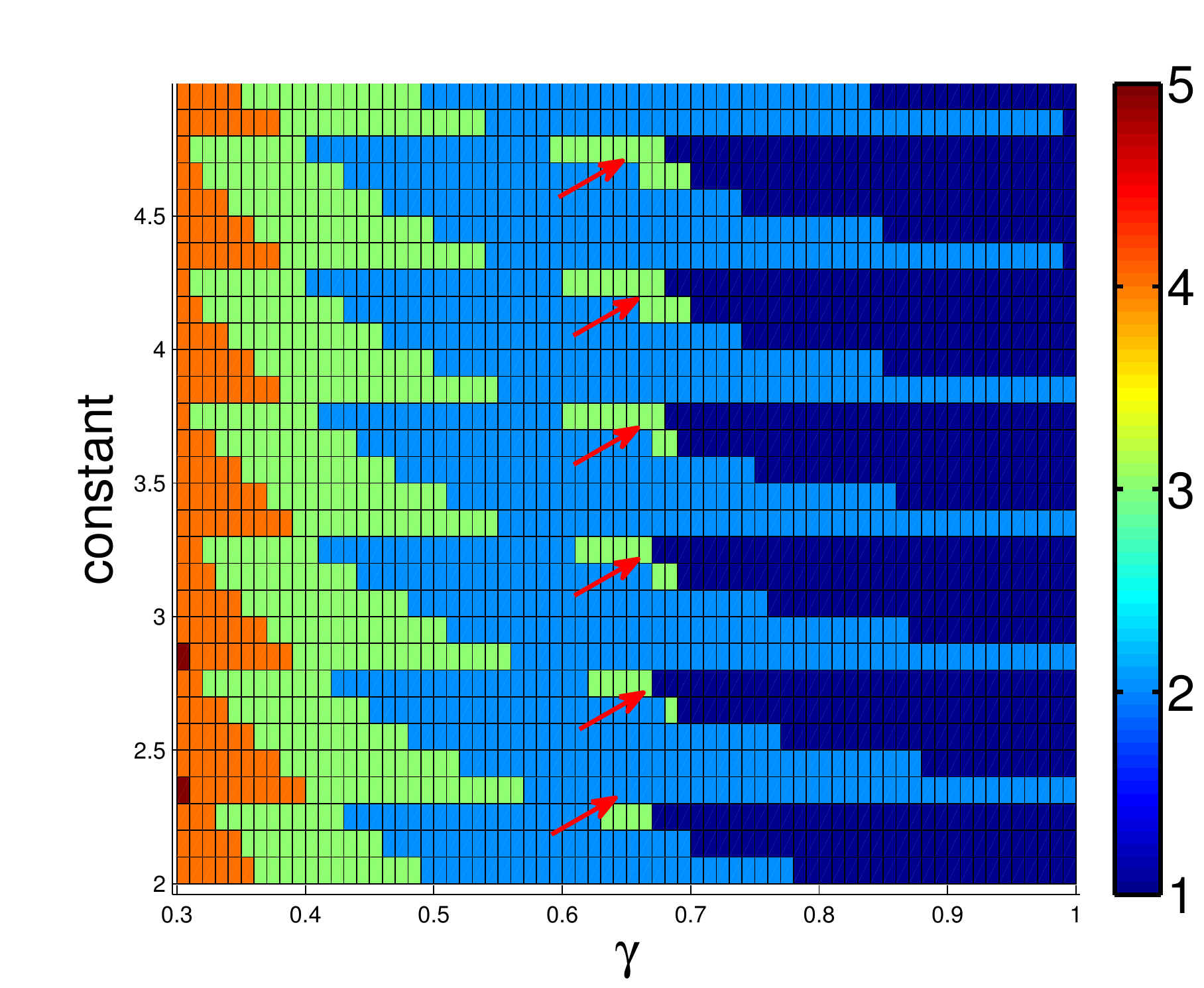}& (b) \includegraphics[scale=0.25]{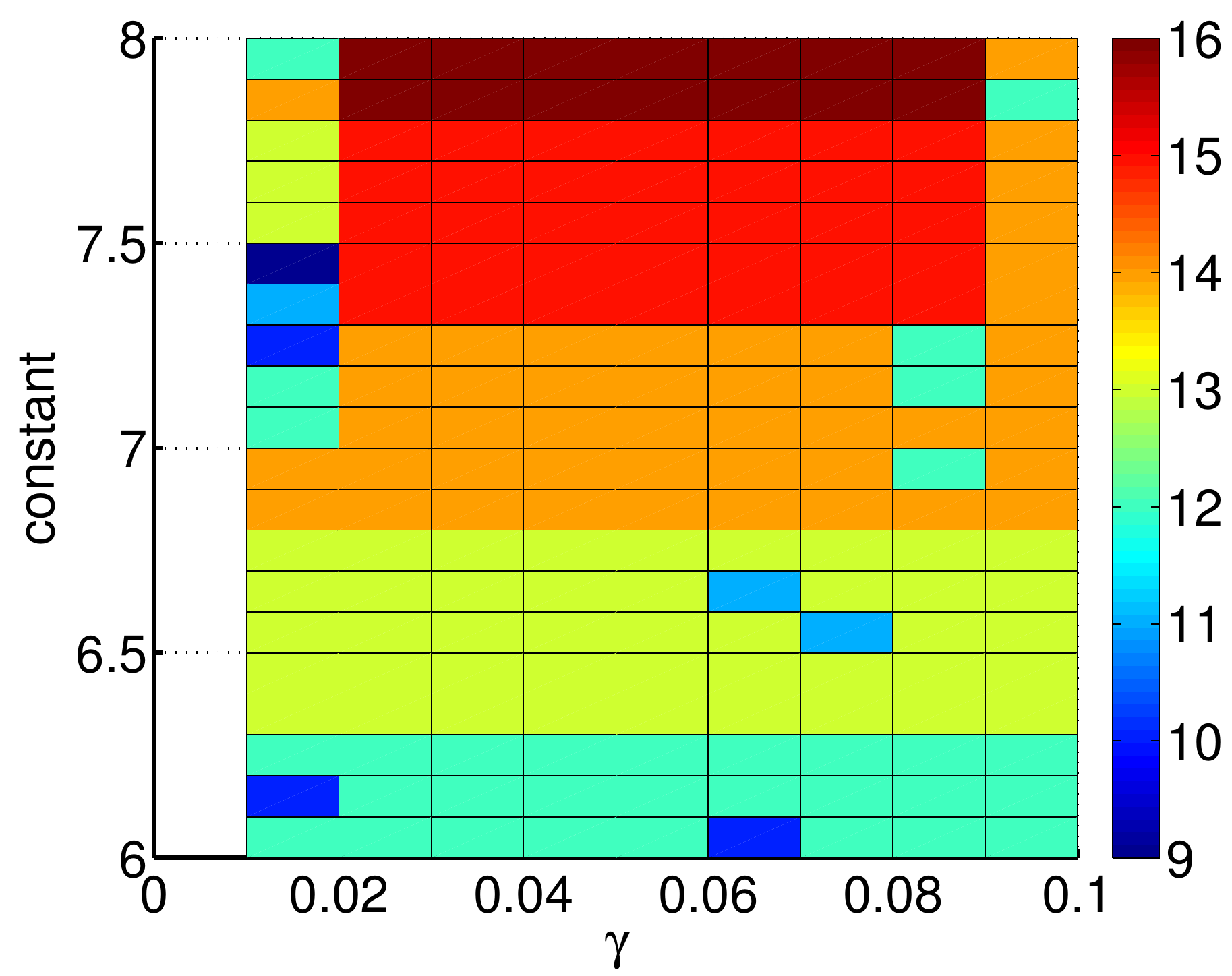}& (c) \includegraphics[scale=0.3]{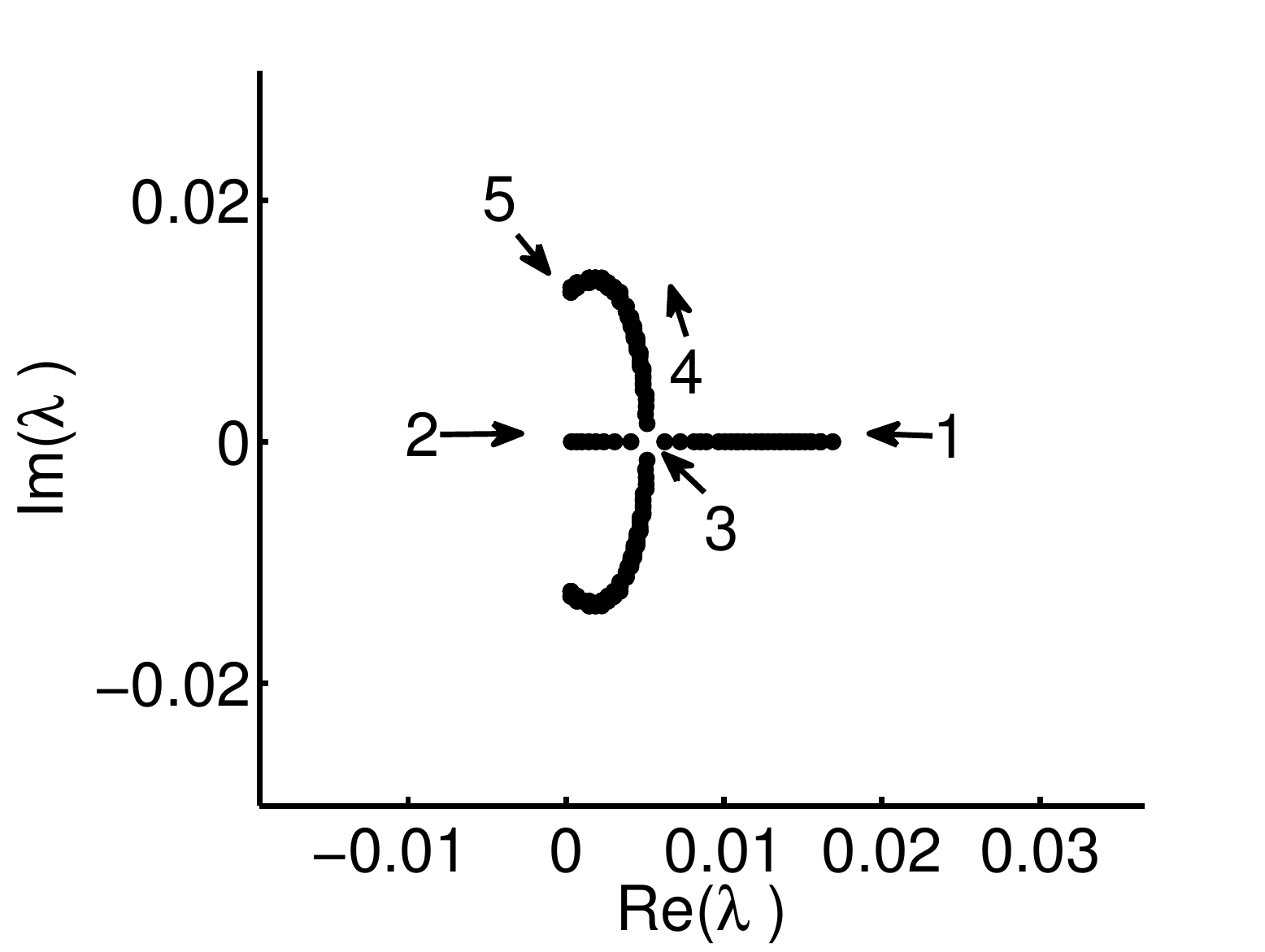}
\end{array}
$
\end{center}
\caption{Designer system. (a) Color plot of the number of roots present for $(\gamma,M\gamma = \mathrm{constant})$. Arrows indicate where the number of roots jumps by 2 corresponding to a Hopf bifurcation. (b)  Color plot of the number of roots present for $(\gamma,M\gamma = \mathrm{constant})$.   (c) Plot in the complex plane of the path two roots of the Evans function take as they cross the imaginary axis. Here $\gamma = 0.65$ is fixed and $M$ varies. There is a root at about $\lambda = 0.09$ not in the viewing window. (1) When $M = 2.57$ there is one root in the viewing window traveling left, (2) at $M = 2.5815$ a root passes through the origin traveling right , (3) the two roots collide at $M \approx 2.585$, (4) upon collision, the roots split off of the real axis, (5) at $M \approx 2.661$ the roots pass through the imaginary axis.    }
\label{fig814}\end{figure}


\begin{figure}[htbp]
 \begin{center}
$
\begin{array}{lccr}
(a) \includegraphics[scale=0.14]{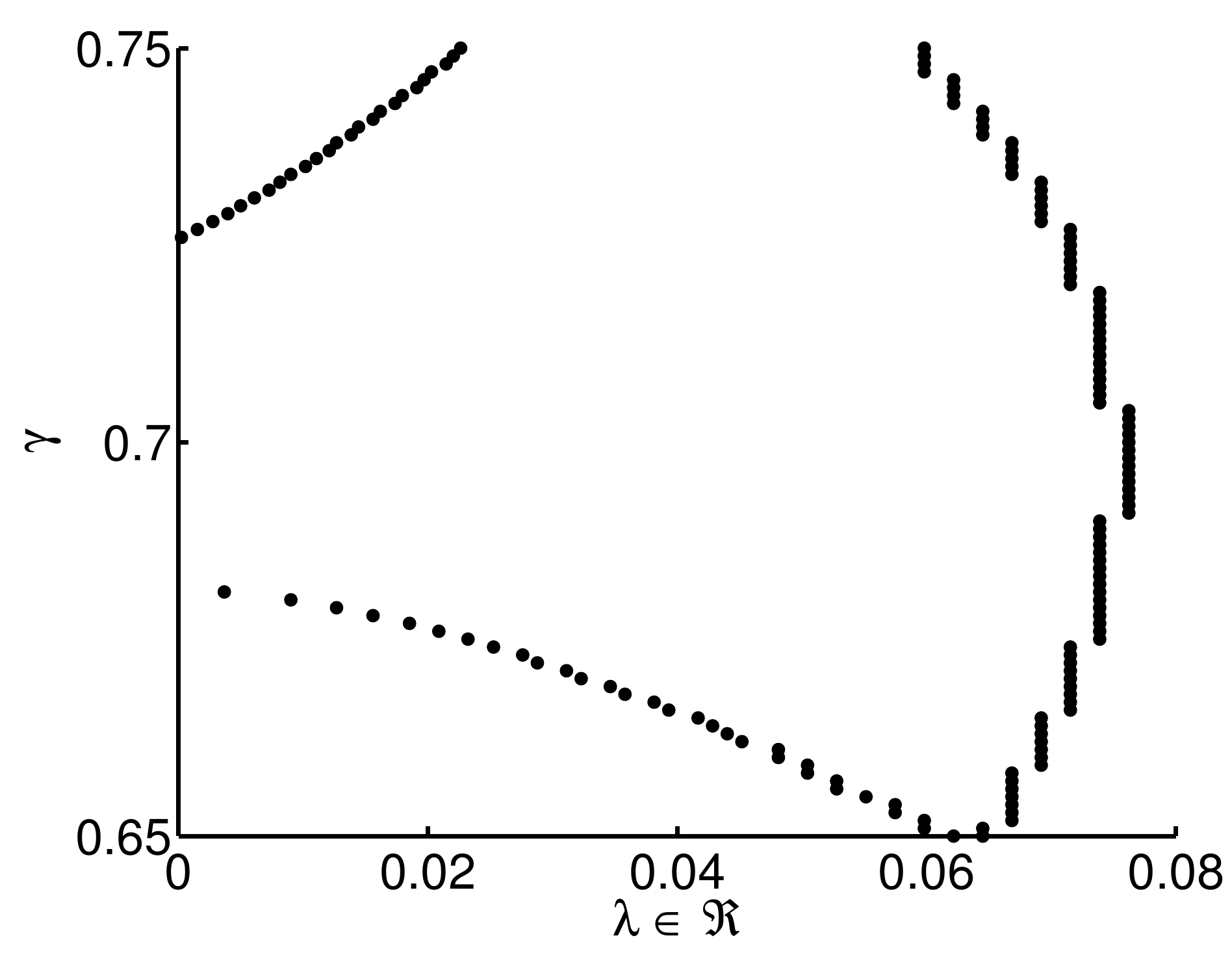}& (b) \includegraphics[scale=0.2]{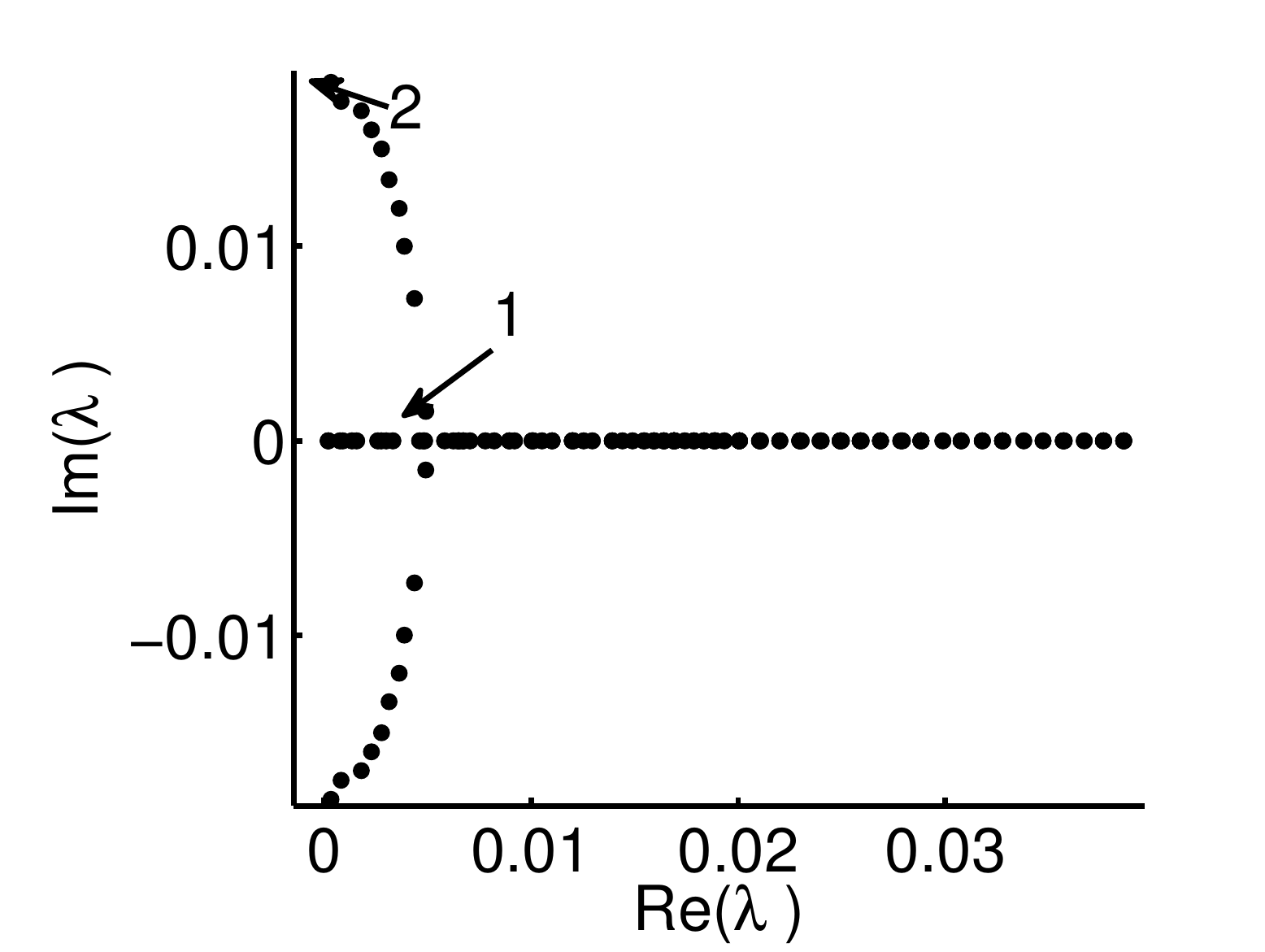}& (c) \includegraphics[scale=0.2]{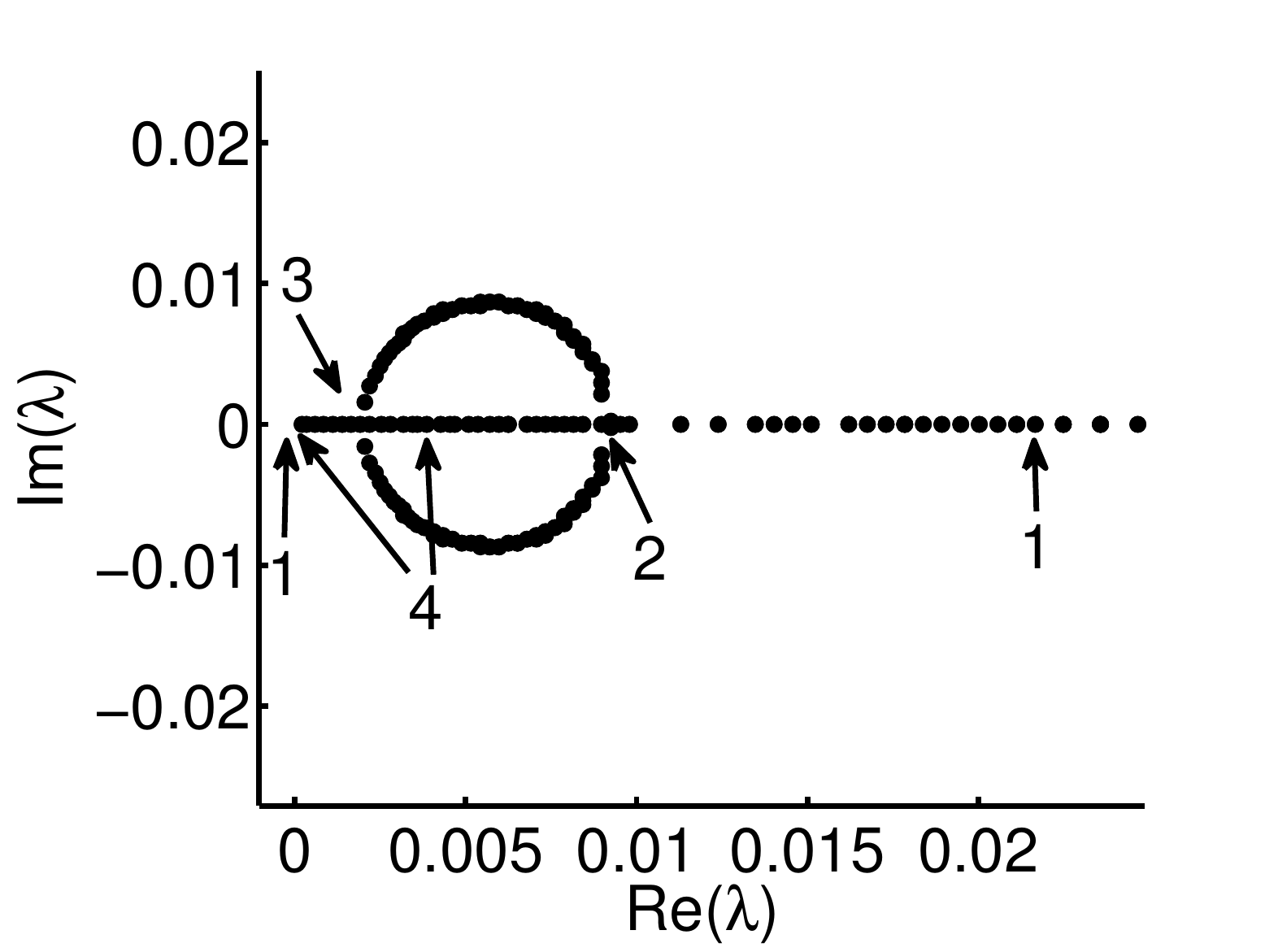}
& (d) \includegraphics[scale=0.14]{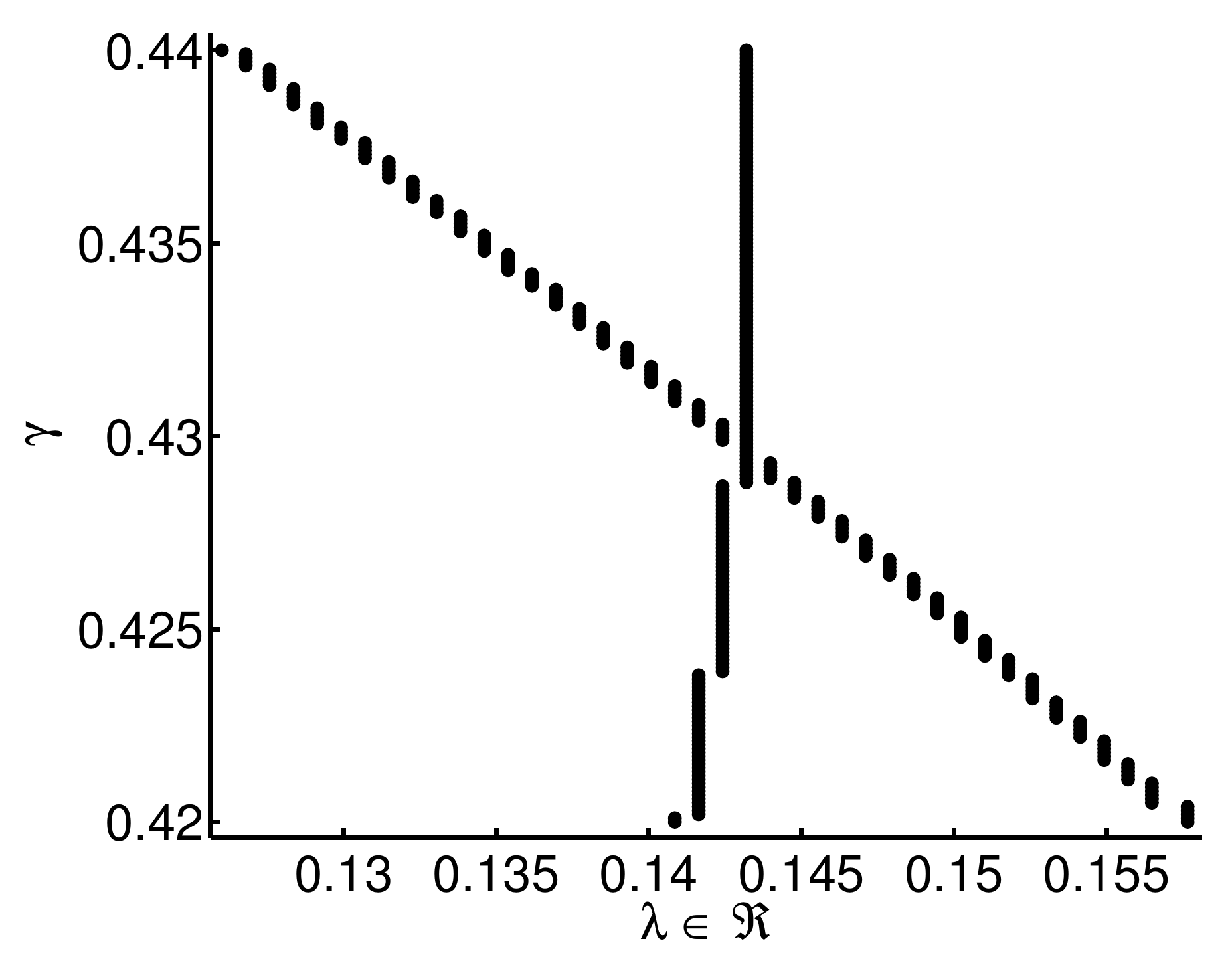}
\end{array}
$
\end{center}
\caption{Designer system. In plots (a)-(c) there is a real root to the right of the viewing window.
(a) Plot of $\gamma$ against root location for $M = 3.1$. The two small modulus roots of the Evans function do not collide and split in the right half plane, but presumably they do in the left half-plane. 
(b)  For $M=3.2836$ two roots of the Evans function collide and split into a non real complex conjugate pair which cross the imaginary axis indicating a Hopf bifurcation. The arrows correspond to the location of the roots of the Evans function for the following values of $\gamma$, (1) $\gamma \approx 0.6545$, (2) $\gamma \approx 0.664$.
(c)  For $M =3.51$ two roots collide and split into a non real complex conjugate pair, but they rejoin in the right half plane. The arrows correspond to the location of roots for the following values of $\gamma$, (1) $\gamma \approx 0.6252$, (2) $\gamma \approx 0.629$, (3) $\gamma \approx 0.6409$, (4) $\gamma = 0.641$.    
(d) Plot of the location of the two largest modulus roots of the Evans function against $\gamma$ for $M = 3.2836$. Note that the roots pass through each other.
 }
\label{fig814b}\end{figure}

\begin{figure}[htbp]
\begin{center}
$
\begin{array}{lcr}
(a) \includegraphics[scale=0.18]{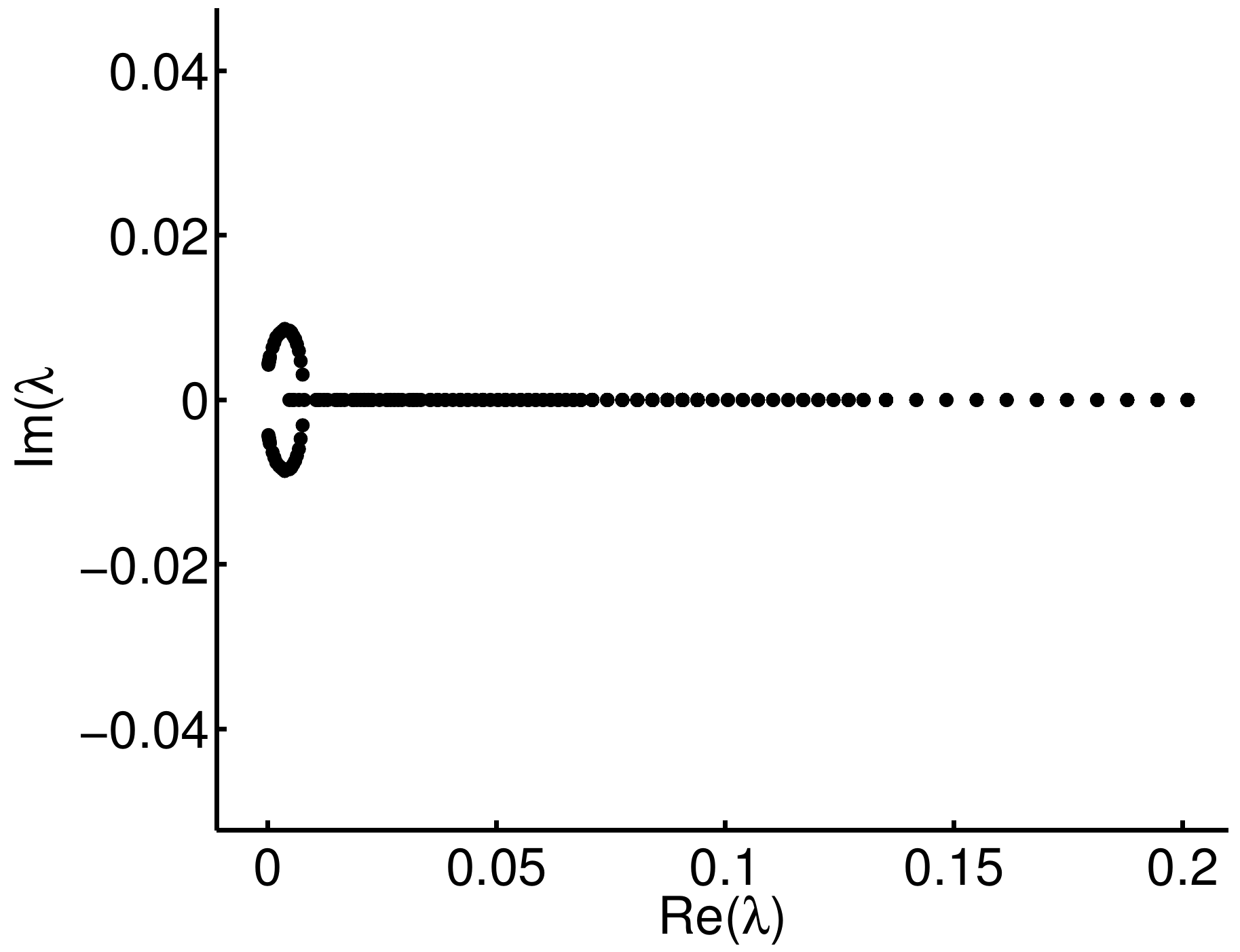} \quad (b) \includegraphics[scale=0.18]{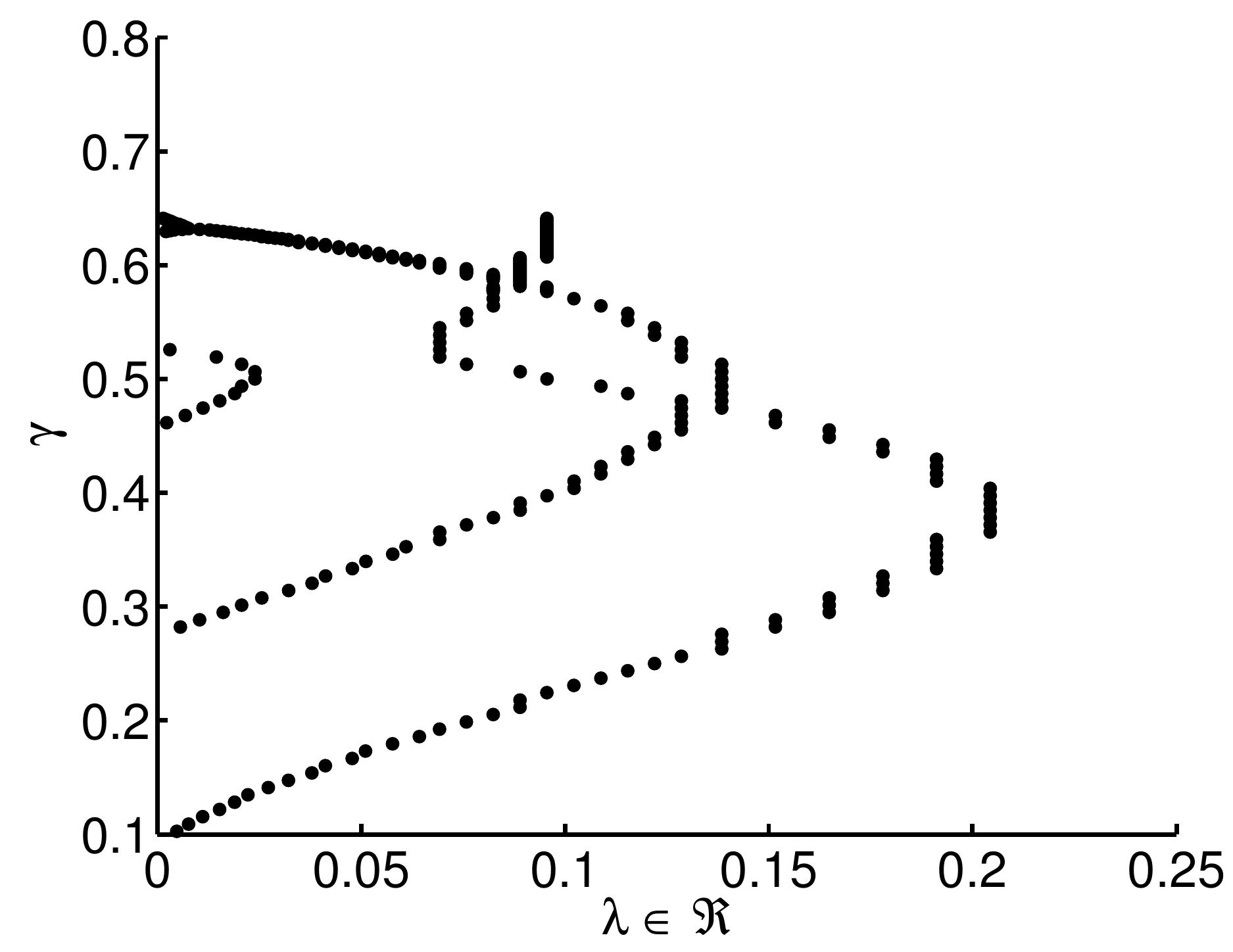} \quad (c) \includegraphics[scale=0.18]{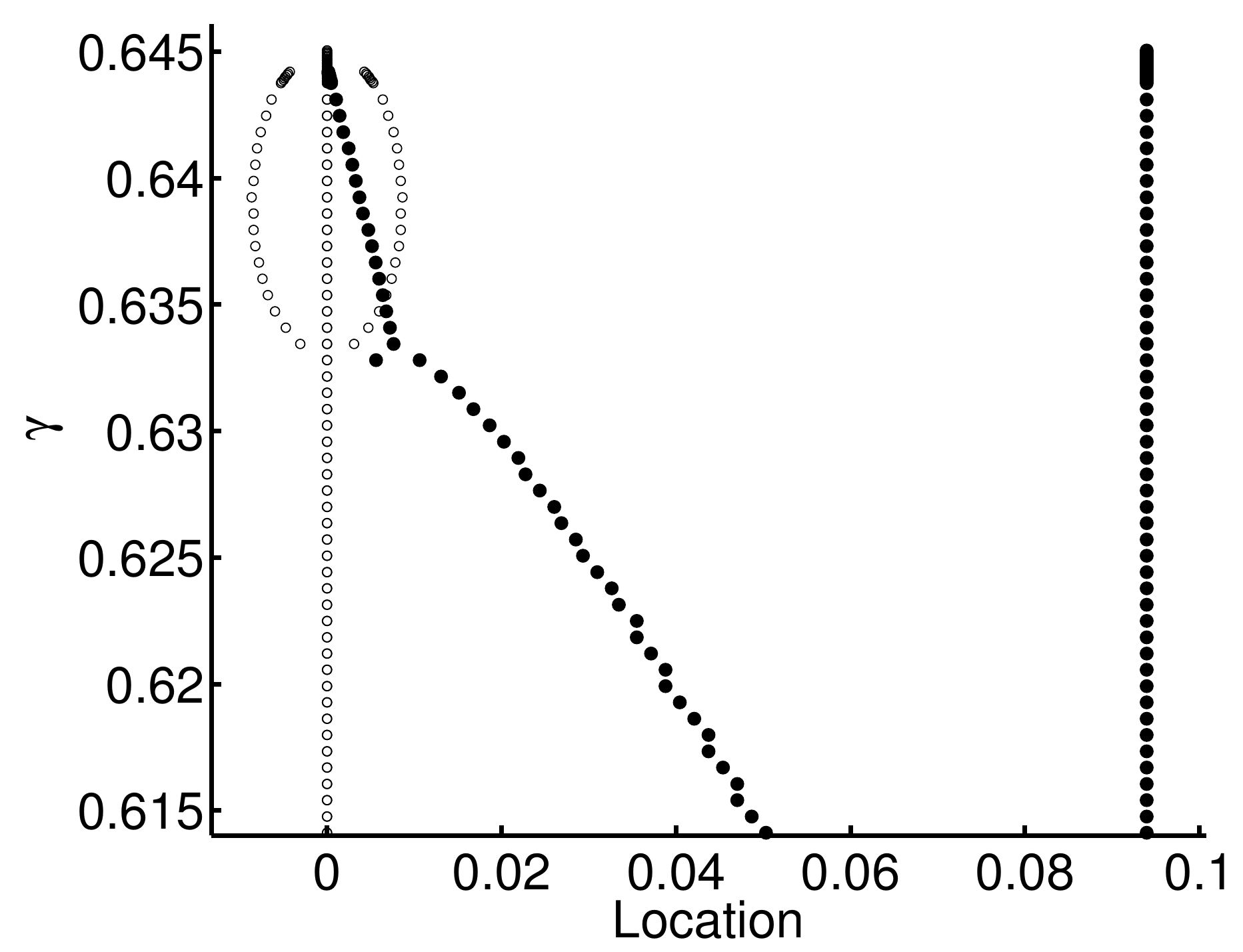}\\
(d) \includegraphics[scale=0.18]{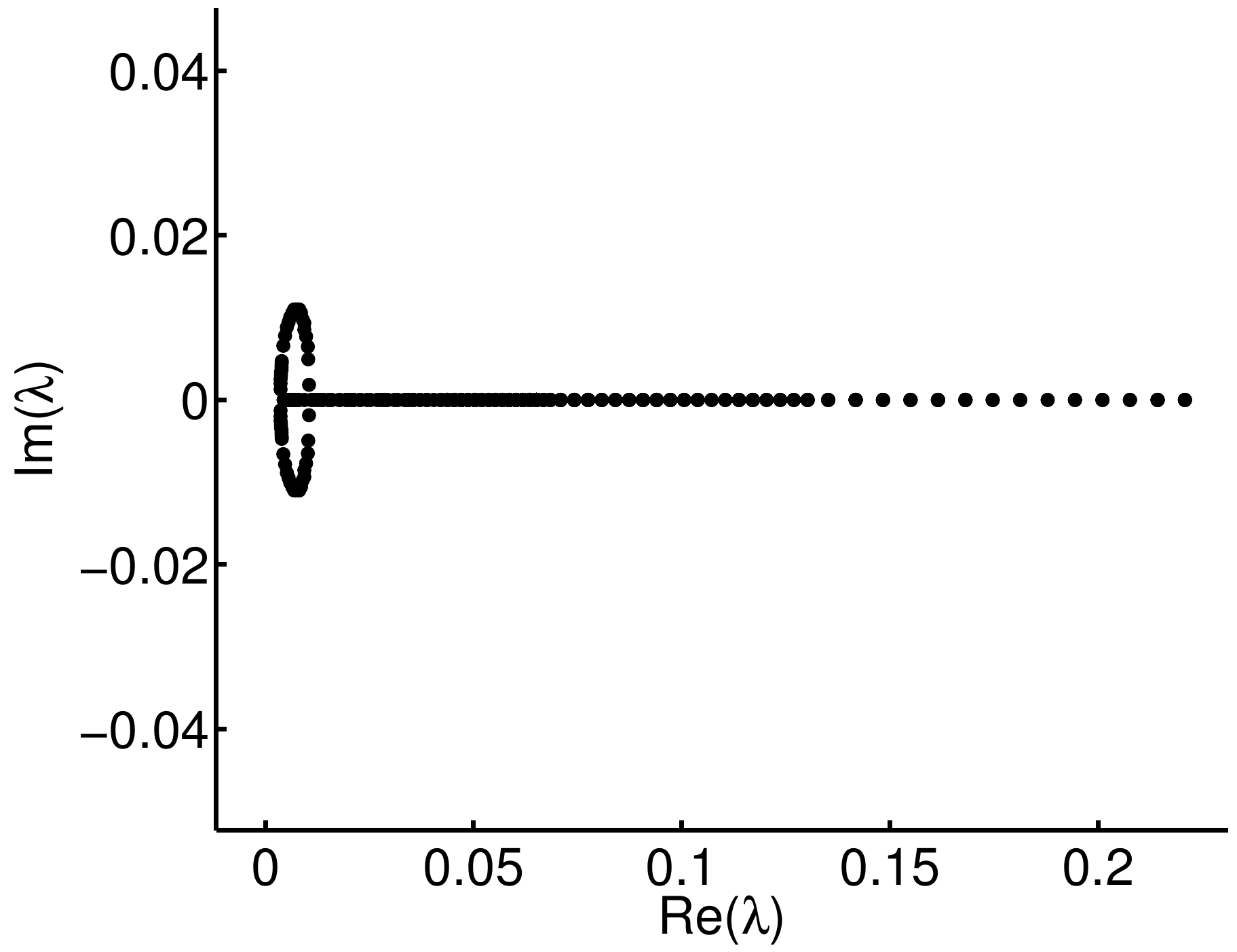} \quad (e) \includegraphics[scale=0.18]{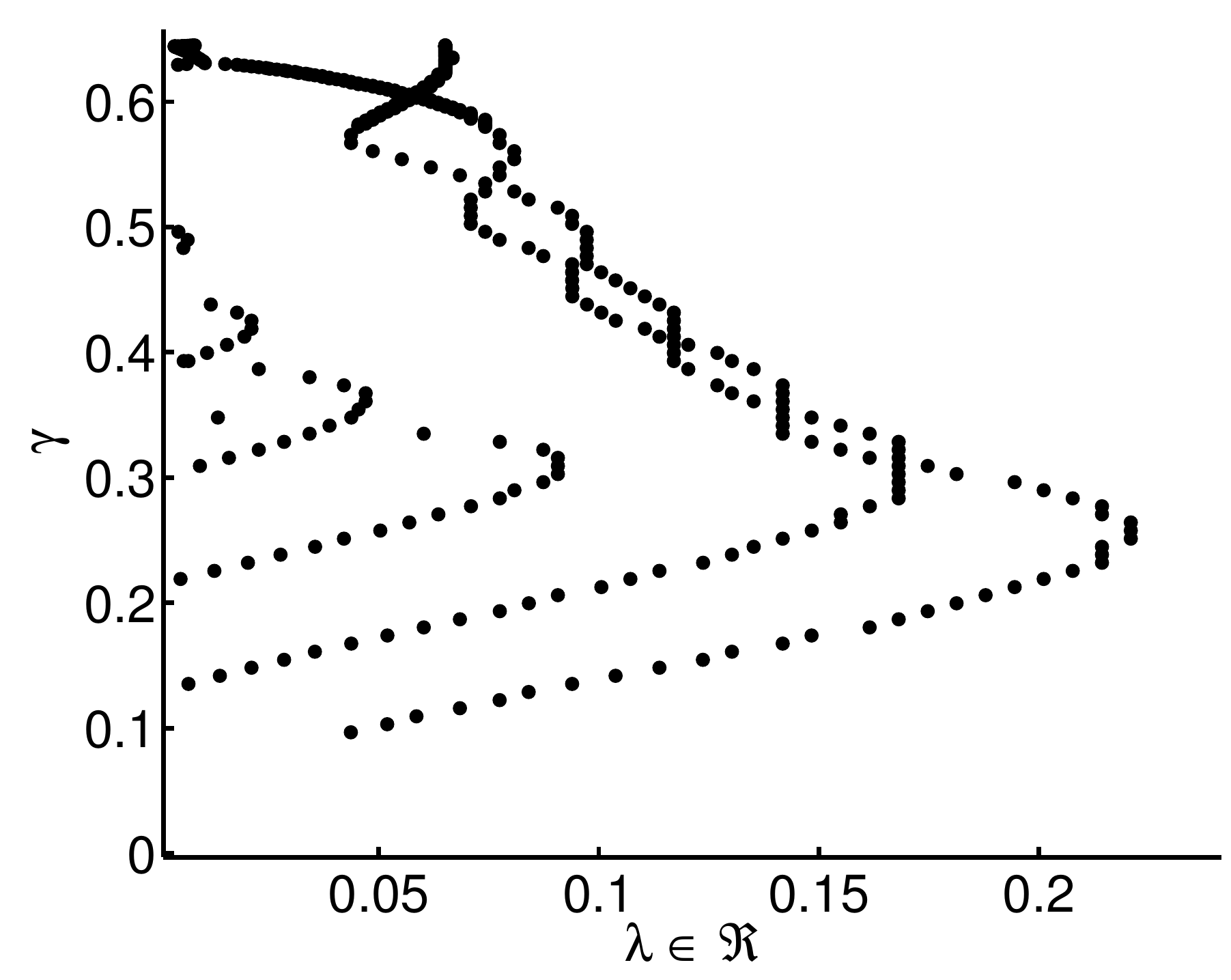} \quad (f) \includegraphics[scale=0.18]{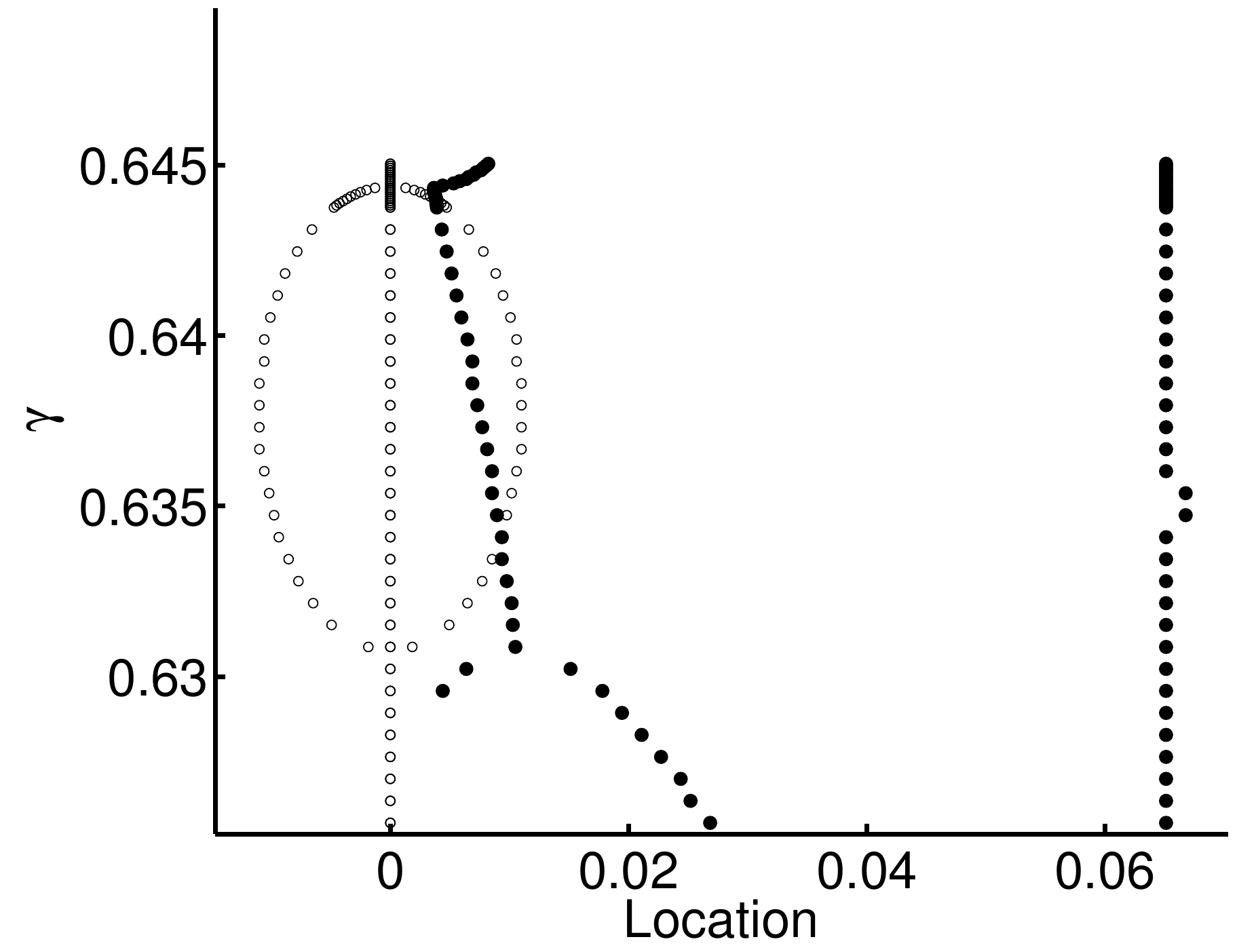}\\
(g) \includegraphics[scale=0.18]{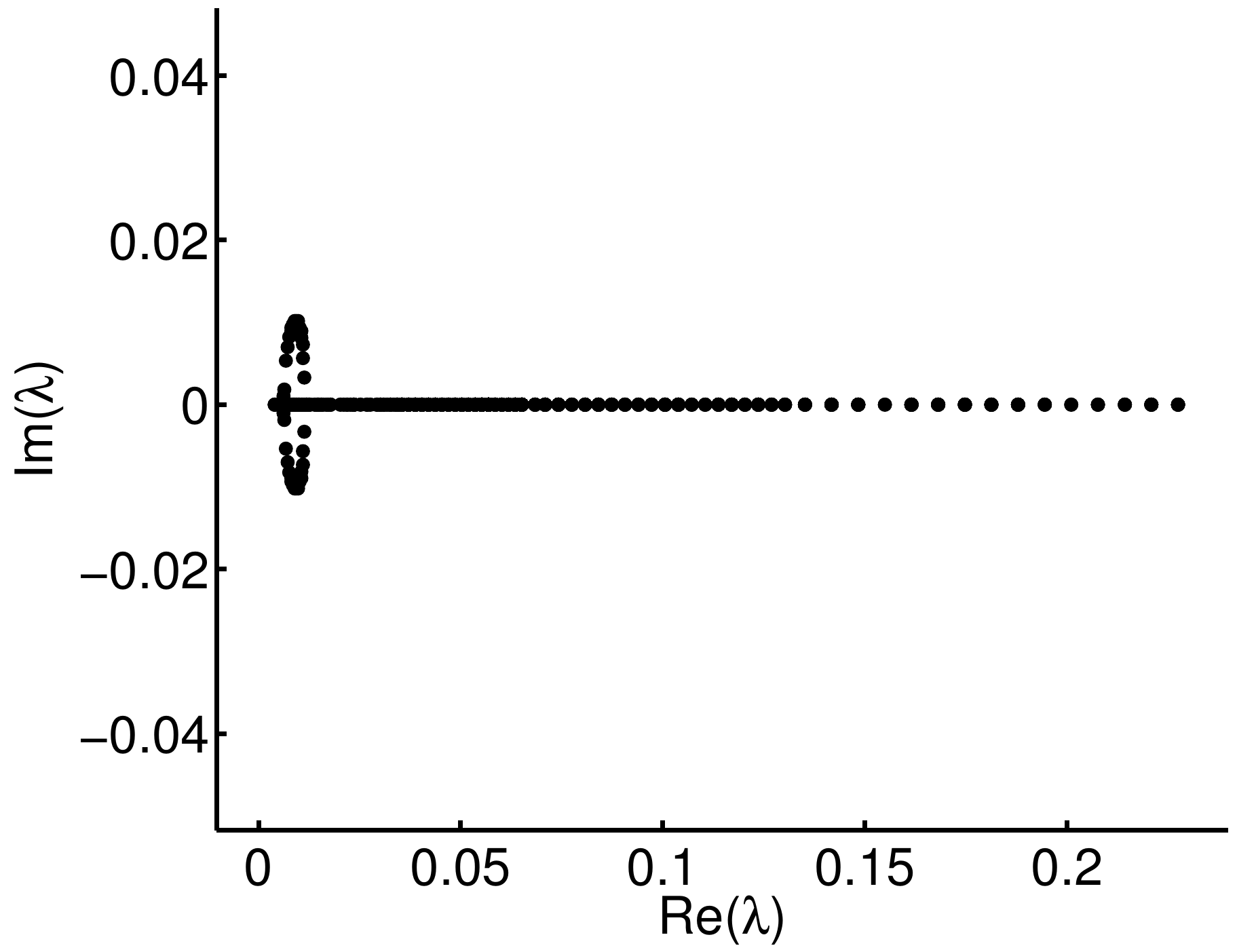} \quad (h) \includegraphics[scale=0.18]{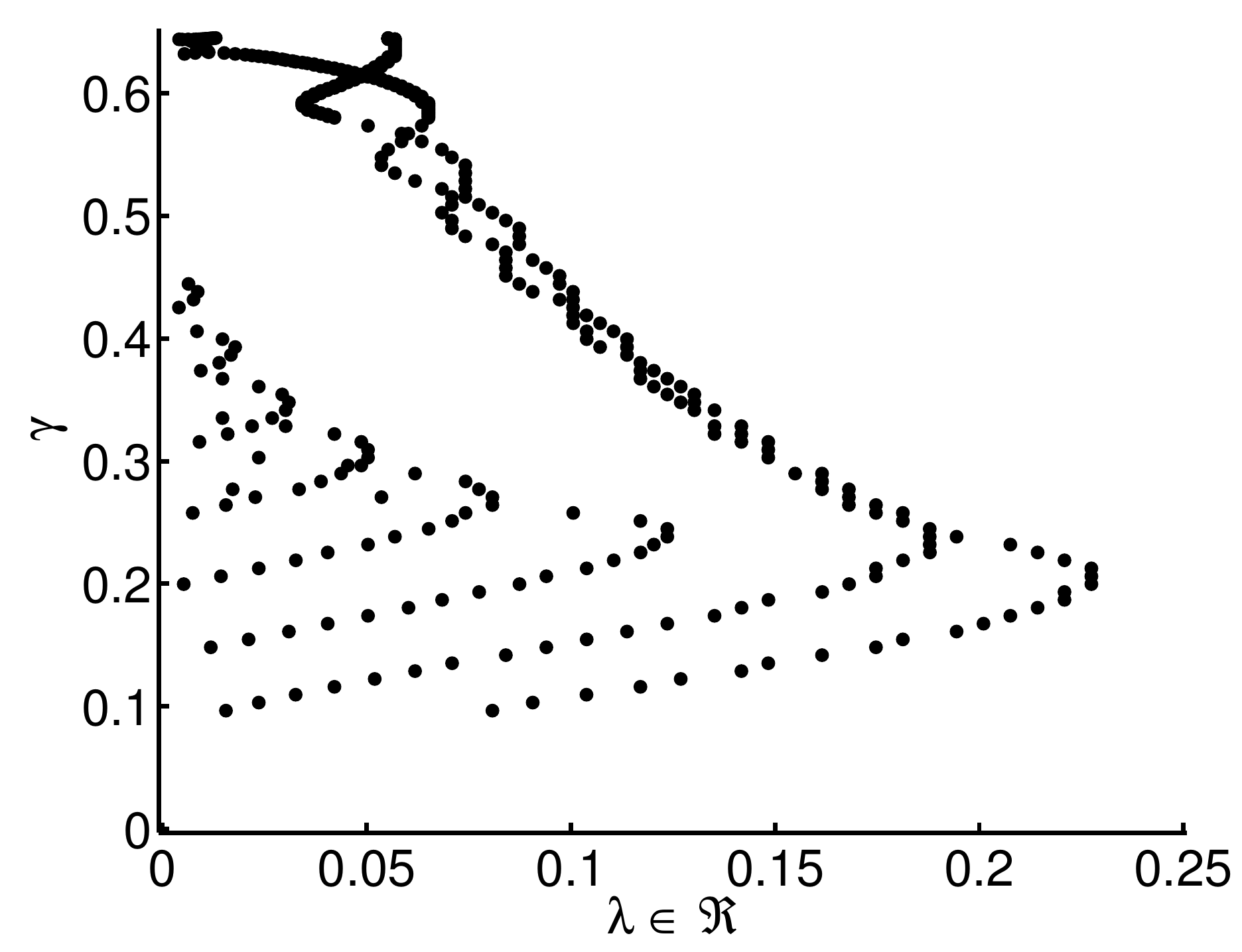} \quad (i) \includegraphics[scale=0.18]{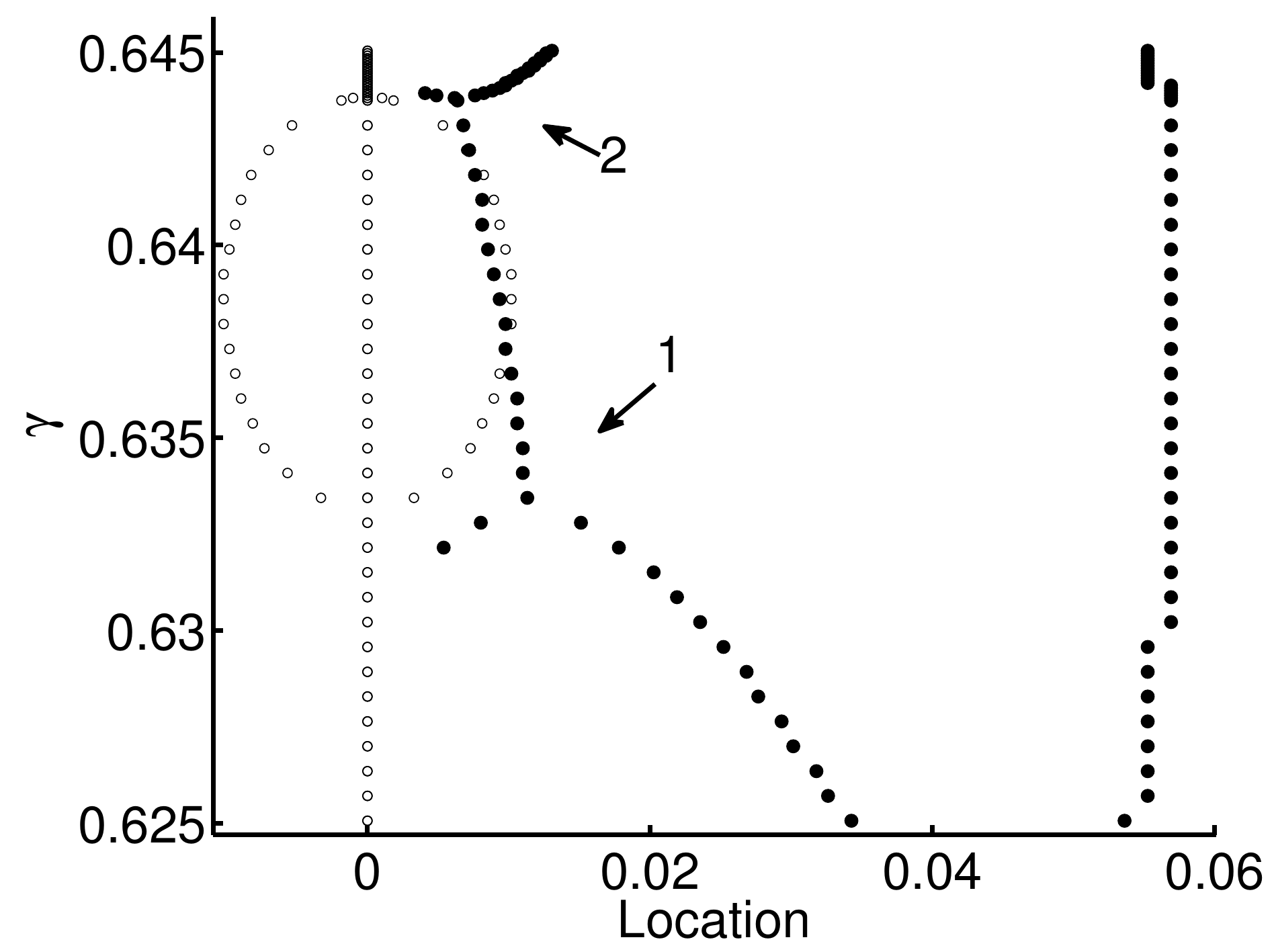}
\end{array}
$
\end{center}
\caption{ For (a)-(c), $M \approx 2.730$. For (d)-(f), $ M \approx 5.8194$,  and for (g)-(i), $M \approx 8.9230$. We have in (a), (d), and (g) plots of roots in the complex plane as $\gamma$ varies. In (b), (e), and (h) we plot $\gamma$ against location. In (c), (f), and (i) we plot $\gamma$ against location where a closed dot corresponds to the real part of the root and an open circle corresponds to the imaginary part.
In (i) arrows 1 and 2 emphasize respectively where two roots on the real axis collide and split to form a complex conjugate pair and then collide again and split along the real axis, similarly as in (c) and (f). There are 3 roots in (a)-(c), 4 in (d)-(f), and 5 in (g)-(i).
}
\label{fig646}\end{figure}

\br\label{nurmk}
Vanishing of $\nu$ appears to be closely related to both Hopf bifurcation
and appearance of complex roots; hence,
unlike the situation for gas dynamics, we have evidently violation of transversality
of profiles within the parameter regime under study.
\er

%
%
%
%

\medskip

{\bf Acknowledgement:} The numerical Evans function 
computations performed in this paper were carried out
using STABLAB \cite{BHZ2}, a MATLAB-based numerical stability package
developed by Jeffrey Humpherys together with Barker and Zumbrun.
We gratefully acknowledge his contribution.
Thanks also to Stefano Bianchini, Fedja Nazarov, and Ben Texier
for interesting conversations on this topic, and in particular
to Bianchini for his contribution through system \eqref{2bsys}
\cite{Bi}.
The second author thanks the University of Indiana, Bloomington
and the third author thanks the \'Ecole Normale Sup\'erieure, Paris,
the University of Paris 13, and the
Foundation Sciences Math\'ematiques de Paris for 
their hospitality during visits in which this work was partially carried out.
Thanks finally to Indiana University
University Information Technology Services (UITS)
for providing the QUARRY supercomputing environment used 
in our large-scale parallel Evans computations.

\appendix
\section{Lopatinski computations for partial equation of state}\label{s:pcond}
To compute the Lopatinski condition in terms of the more
general relation $p=\hat p(\tau,e)$, not necessarily connected
with a full equation of state,  
rewrite \eqref{euler} for smooth solutions as
\ba\label{aeuler}
\tau_t-v_x&=0,\\
v_t+ p_x&=0,\\
e_t+pv_x&=0,\\
\ea
so that, for $w=(\tau,v,e)^T$, 
$
A_\pm=
\begin{pmatrix}
0 & -1 & 0\\
\hat p_\tau& 0 &\hat  p_e\\
0&p&0
\end{pmatrix}_\pm,
$
$a_1=-c$, $a_2=0$, $a_3=c$, $c=\sqrt{ -\hat p_\tau + p \hat p_e}$,
and
\be\label{ars}
r_1^+=
\begin{pmatrix}
1\\ c \\ -p\\
\end{pmatrix},
\quad
r_2^+=
\begin{pmatrix}
-\hat p_e\\ 0\\ \hat p_\tau
\end{pmatrix},
\quad
r_3^+=
\begin{pmatrix}
1\\ -c \\ -p\\
\end{pmatrix},
\ee
so that hyperbolicity corresponds to $\hat p_\tau - p \hat p_e <0$.
Computing
$
A_\pm^0=
\begin{pmatrix}
1 &  0&0\\
0 & 1&0\\
0 &v& 1
\end{pmatrix}_\pm,
$
we obtain
$$
A_+^0 r_1^+=
\begin{pmatrix}
1\\ c \\ vc - p\\
\end{pmatrix},
\quad
A_+^0 r_2^+=
\begin{pmatrix}
-\hat p_e\\ 0\\ \hat p_\tau
\end{pmatrix},
\quad
A_+^0 r_3^+=
\begin{pmatrix}
1\\ -c \\  -vc-p\\
\end{pmatrix},
$$

Using \eqref{gasRH} to compute $[U]=[\tau](1,-\sigma, -p)$, we
obtain for a $1$-shock
\be\label{ac2}
\check \delta=
[\tau]
\det
\begin{pmatrix}
1& -\hat p_e & 1\\
-\sigma & 0 & -c\\
-p & \hat p_\tau
 & -vc - p\\
\end{pmatrix},
\ee
where all quantities are evaluated at $U_+$, or,
using $p\hat p_e-\hat p_\tau=c^2$,
\ba\label{althatp}
\check\delta [\tau]^{-1}&=
(\sigma -c)c^2 + \hat p_e\sigma vc
= (\sigma -c)c^2 + \hat p_ec [p],
\ea
recovering \eqref{dcond3} (noting by homotopy to
the small amplitude case that $\check \delta>0$ $\sim$ stability) 
in the form
$
\frac{\hat p_e}{ c^2}
< \frac{-\frac{\sigma}{c}+1} {  [p] } ,
$
equivalent to \eqref{dcond3} by the relation 
$\hat p_e=\frac{\bar p_S}{\bar e_S}=
\frac{-\bar e_{\tau S}}{\bar e_S}$,
where $\bar p(\tau,S):=-\bar e_\tau(\tau,S)$.
This yields the alternate forms 
$\frac{\hat p_e}{ c^2} < {2}{  [p] } $
and
$\frac{\hat p_e}{ c^2} > {1}{  [p] } $
of \eqref{strongintro} and \eqref{weakintro}, convenient for computations
in 
form \eqref{euler} involving only a  pressure law $\hat p$
and not a complete equation of state.
When $\bar e_{\tau S}<0$,
this is equivalent to Majda's condition \eqref{mcond} with
$\Gamma= \frac{-\tau \bar e_{\tau s}}{\bar e_s}$
rewritten as $\Gamma= \tau \hat p_e$.

\br
As we have seen, the computations carried out here are simpler, if anything,
than the computations in entropy coordinates \eqref{ent_euler}.
However, what is not clear from this formulation is why the computations
simplify so, involving only the quadratic formula for a $3\times 3$
determinant, a circumstance that originates from the decoupled form
of the equations in entropy variables.
\er

\section{Helmholtz energy and associated energy functions}\label{s:helmholtz}
The completion of a pressure law $p=P(T,\tau)$
by a compatible energy law $e=E(T,\tau)$
may be obtained more systematically by referring to a related 
circle of ideas involving the Helmholtz energy
\be\label{Adef2}
A:=e -TS,
\ee
where $e$, $T$, and $S$ have their usual meanings.
Here, $A$ is the Legendre transform of $e$ with respect
to $S$, $T=e_S$; that is, the relation between $e$ and $A$
is analogous to the relation between Hamiltonian and Lagrangian
energies in classical mechanics. 
See \cite{W1} or \cite{MP} for background discussion.
For a nice dicscussion of the Legendre transform and relations to
the principle of least action, see \cite{W2}.

The key properties of the Helmoltz energy from our point of 
view are as follows.

\bl\label{hlem}
Considering $ A=\tilde A(T,\tau)$, we have
\be\label{checkA}
p=-\tilde A_\tau,
\qquad
S=-\tilde A_T.
\ee
\el

\begin{proof}
By direct computation, $\tilde A_T= \bar e_S (dS/dT)
-  T (dS/dT) - S= -S$, verifying the second relation.
Similarly, 
$$
\tilde A_\tau= \tilde e_\tau 
-  T (\partial \tilde S/\partial \tau) 
=
\bar e_\tau + \bar e_S (\partial \tilde S/\partial \tau)
-  T (\partial \tilde S/\partial \tau) 
=
\bar e_\tau,
$$
verifying the first.
\end{proof}

Evidently, using \eqref{checkA}, we may recover $A$ from the pressure
law $p=\tilde p(T,\tau)$, up to an arbitrary function of $T$, by
antidifferentiation in $\tau$.
More fundamentally
\be\label{AQ}
A=-R T\ln \tilde Q,
\ee 
where $\tilde Q(T,\tau)$ is the
(rescaled) {\it partition function} coming from the 
statistical thermodynamic derivation of the equation of state.

\medskip
{\bf Exactness.}  To check whether a pair of gas laws
$e=\tilde e(T,\tau)$ and $p=\tilde p(T,\tau)$ are ``exact,'' that is,
originate from a complete equation of state $A=\tilde A(T,\tau)$,
or, equivalently, $e=\bar e(\tau, S)$,
we may note that differentiating $S=\frac{e-A}{T}$ with respect to $\tau$ 
yields $S_\tau=\frac{e_\tau + p}{T}$, hence exactness, $S_\tau=p_T$,
is equivalent to
\be\label{test}
Tp_T= e_\tau+ p.
\ee

We can use this to determine also exactness of a pair of gas laws in
the alternative form $p=\hat p(\tau,e)$, $T=\hat T(\tau,e)$, using
the Implicit Function Theorem to rewrite \eqref{test} as
\be\label{alttest}
T\hat p_e = p \hat T_e -\hat T_\tau.
\ee
(Here, we are implicitly assuming $\hat T_e>0$, so that we may invert
$T=\hat T(e,\tau)$ to solve for $e=\tilde e(T,\tau)$.)

\medskip
{\bf Hyperbolicity and convexity.}  
The hyperbolicity condition may be deduced using
$\bar p(S,\tau)=\tilde p(\tau, e_S(S,\tau))$ to obtain
$
\bar p_\tau=  \tilde p_\tau + \tilde p_T \bar e_{S\tau}
=  \tilde p_\tau -\tilde p_T \bar p_{S}
=  \tilde p_\tau -\tilde p_T \tilde p_{T} (\partial T/\partial S)
,
$
or, by $S=-\tilde A_T$, hence $\partial S/\partial T= -\tilde A_{TT}$,
$$
\bar p_\tau =  
\tilde p_\tau + \frac{\tilde p_T \tilde p_{T}}{\tilde A_{TT}} 
=
\frac{\tilde A_{T\tau}^2-\tilde A_{TT}\tilde A_{\tau \tau}}{\tilde A_{TT}}.
$$
Thus, hyperbolicity, $\bar p_\tau <0$, holds when
$\tilde A_{TT}>0$ and $\tilde A$ is convex as a function of $(\tau, T)$,
or else $\tilde A_{TT}= -\partial S/\partial T=-1/\bar e_{SS} <0$ 
(as usual) and $\tilde A$ is {\it not concave}
.
Convexity of $\bar e$
may be computed similarly to be equivalent to
$\tilde A_{TT}<0$ and
$
\bar e_{SS}\bar e{\tau \tau}- \bar e_{S\tau}^2=
\frac{\tilde A_{\tau \tau} } {\tilde A_{TT}}
>0,
$
or $\tilde A_{\tau \tau}>0 $ and $\tilde A_{TT}<0$.

\begin{example}[Ideal gas]
Here, we have 
$A=-\tilde RT\ln \tau -C\tilde RT,$
from which we find as usual
$\check e(T,\tau)=\frac32\tilde RT $ (monatomic gas law, $\gamma=5/3$.)
\end{example}

\begin{example}[van der Waals gas]
From the reference \cite{W1}, we find for $\Gamma=2/3$,
\be\label{A}
A=-\tilde R T \ln \Big((\tau -\tilde b)T^{3/2}\Big) - \frac{\tilde a}{\tau }
\ee
where the rescaled partition function $\tilde Q(T,\tau)$ satisfies
$\ln \tilde Q(\tau, T)=
\ln (\tau -\tilde b) + \frac{\tilde a}{\tilde R \tau T} + C\ln T^{3/2} ,$
$C=\const,$
and
$S=-A_T =\tilde R\big( \ln \big((\tau -\tilde b)T^{3/2}\big) +\frac 52\big).$
Thus, to obtain $\check e$, we can just compute 
$\check e(T,\tau)=A+Ts =
\frac32 \tilde R T -\frac{\tilde a}{\tau} .$
\end{example}

\begin{example}[Redlich--Kwong gas]
Here \cite{W1}, 
$A=-\tilde R T\Big( \ln (\tau -\tilde b)T^{3/2}\Big) 
- \tilde a T^{-1/2}\ln(1+\tilde b/\tau).$
\end{example}

\section{Numerical Evans function protocol}\label{s:protocol}
For completeness,
we briefly describe in this appendix the procedure followed in
the numerical Evans function studies of part II.
For further details, see, e.g., \cite{BHZ1,BHRZ,BLeZ,Z4,Z5}.

\subsection{Profile and Evans function equations}\label{s:egen}
We first recall the convenient general forms derived in \cite{BHLyZ1}
for profile and Evans function equations of a general class of
systems \eqref{vgencon} 
arising in continuum mechanics: specifically, systems for which $B$ has 
block-diagonal structure
\be\label{Bstruct}
B(w)=\bp 0 & 0\\
0 & b(w)\ep,
\quad
b\in \RR^{r\times r},
\quad
\Re \sigma(b)>0,
\ee
and, along the traveling-wave profile $\bar w$,
\beq\label{eq:a11inv}
\det (df^1_{11} -\sigma df^0_{11}) (\bar w) \neq 0\, 
\; \hbox{\rm (hyperbolic noncharacteristicity)}.
\eeq

\begin{example}
Navier--Stokes equations \eqref{NS}--\eqref{nslaws} 
have form \eqref{vgencon}, \eqref{Bstruct}
with
$w=(\tau,v,T)$,
$f^0(w)=(\tau, v, \check e(\tau,T)+v^2/2)^T$,
$f^1(w)=(-u, \check p(\tau,T), v\check p(\tau,T))^T$,
and $b(w)= \tau^{-1} \bp   \mu & 0\\ v\mu& \kappa\ep$.
\end{example}

\begin{example}\label{arteg}
Designer system \eqref{2bsys} 
has form \eqref{vgencon}, \eqref{Bstruct} with 
$w=(u,v)^T=f^0(w)$, $r=n$.
\end{example}

\subsubsection{Traveling-wave system}
By the block structure assumption \eqref{Bstruct}, 
\eqref{twode} decomposes into
\begin{subequations}\label{eq:tw}
\begin{align}
0 & = f^1_1(\bar w)-f^1_1(U_-)- \sigma (f^0_1(\bar w)-f^0_1(U_-))\,, \label{eq:alg}
 \\
b(\bar w) \bar w_2' &=\tilde f^1_2(\bar w)-\tilde f^1_2(U_-)
- \sigma (f^0_2(\bar w)-f^0_2(U_-)) \,,\label{eq:dif}
\end{align}
\end{subequations}
where $f^j_k$ denotes the $k$th block of $f^j$.
By the Implicit Function Theorem, \eqref{eq:a11inv} guarantees that
the first equation may be solved, locally, for 
$\bar w_1$ as a function of $\bar w_2$, either analytically
(preferable) as done here for gas dynamics
or numerically, so that the second
equation defines a flow in $\bar w_2$. 
Viscous shock profiles are computed using the reduced system \eqref{eq:dif}.

\subsubsection{Integrated Evans system}\label{ssec:flux1d}
Changing to co-moving coordinates $x \to x-\sigma t$ in which
$\bar w$ is a steady solution, we obtain the linearized eigenvalue
equation
\beq\label{eq:eval}
\lambda A^0w+(A^1 w)'=(B w')'\,,
\eeq
where
\be\label{coeffs}
A^0:=df^0(\bar w)\,,
\quad A^1w:=df^1(\bar w)w-\sigma df^0(\bar w)-\dif B(\bar w)(w,\bar w_x)\,,
\quad B:=B(\bar w)\,.
\ee

Defining $W':=\bar A^0U$, we thus have, integrating \eqref{eq:eval},
\be\label{ieig}
\lambda W +A^1 (A^0)^{-1} W' =B( (A^0)^{-1} W')',
\ee
or, setting 
$ Z:=\bp  W\\ (0,I_{n-r}) (\bar A^0)^{-1}W'\ep\in \CC^{n+r}, $
and solving for $(I,0)(\bar A^0)^{-1}W'$ in the first block of \eqref{ieig}
using again the assumption \eqref{eq:a11inv} (see \cite{BHLyZ1}),
the \emph{first-order integrated Evans system}
\be\label{Zeq}
Z'=\mathbf{A}_\mathrm{int}Z,
\quad
\mathbf{A}_\mathrm{int}:=
\bp
-\lambda A^0_{11}(A^1_{11})^{-1} &  0 &  A^0_{12}-A^0_{11}(A^1_{11})^{-1} A^1_{12} \\
-\lambda A^0_{21}(A^1_{11})^{-1} & 0  & A^0_{22}-A^0_{21}(A^1_{11})^{-1} A^1_{12} \\
-\lambda b^{-1}A^1_{21}(A^1_{11})^{-1} & \lambda b^{-1} &
b^{-1}(A^1_{22}-A^1_{21}(A^1_{11})^{-1} A^1_{12})
\ep.
\eeq
The {\it integrated Evans function} $\tilde D(\lambda)$ is computed as a 
Wronskian of \eqref{Zeq}.

\subsection{Winding number computations}\label{s:num}
The Evans function is analytic in the closed
right-half plane where,
away from $\lambda=0$, zeros of the Evans function match eigenvalues of the system in both location and multiplicity. Thus we can determine spectral stability of the system by numerically computing the winding number of the Evans function on a semicircle $B(0,R)\cap \{\Re(\lambda)\geq 0\}$ in the right-half plane with $R$ chosen large enough to ensure 
that any possible unstable eigenvalues are inside the semicircle. 
We determine the value of $R$ either analytically
by energy estimates as in  \cite{Br} or 
numerically by convergence to theoretically-predicted high-frequency
asymptotics
\be\label{hfas}
D(\lambda)\sim C_1e^{C_2\sqrt{\lambda}}
\ee
as in \cite{HLyZ1,BLZ,BLeZ,BZ1,BZ2}.
The former method is carried out on a case-by-case basis,
the latter method supported automatically in the MATLAB-based 
numerical Evans function package STABLAB \cite{BHZ2}.

Here, we follow the standard practice of solving,
not the eigenvalue equations \eqref{eq:eval}, but
their integrated form \eqref{ieig} \cite{Go,ZH,HuZ1,HLZ}, thus removing the
zero of $D$ occuring at $\lambda=0$ as a result of translational
invariance of the underlying system \eqref{vgencon} (see \eqref{evcon}).
That is, we compute the interated Evans function $\tilde D$ instead
of $D$.
This allows us to perform well-conditioned winding number computations
directly through $\lambda=0$ throughout most of the computational
domain, greatly speeding performance- the exception being at or near
the precise parameter values at which stability transitions appear,
at which this difficulty is unavoidable.

To locate unstable roots when they occur,
we use the method of moments as described in \cite{Bro,BJNRZ,BZ1}. 
computing moments
of roots within a contour $\Gamma$ by the generalized winding number
computation
$
\sum_{r_j\in \Gamma^{\rm int} }r_j^l
=\frac{1}{2\pi i}\oint_{\Gamma}\frac{\lambda^l\partial_\lambda \tilde D(\lambda)}
{\tilde D(\lambda)}~d\lambda,
$
recovering $r_j$ by solution of an resulting polynomial equation, and
using a divide-and-conquer strategy to keep the number of roots per
subcontour small enough for optimal accuracy vs. speed 
(in practice, one or two).

Standard practice (see, e.g., \cite{BHRZ,HLZ,BLeZ}), 
is to approximate the traveling wave profile 
using \textsc{MATLAB}'s boundary-value solver
{\tt bvp6c} \cite{HM},
an adaptive Lobatto quadrature scheme, with projective boundary conditions,
and on a truncated, finite computational domain
$[-L,L]$,  where the values of approximate plus and minus spatial infinity $L$ 
are determined experimentally by the requirement that the absolute error 
between the computed profile values at $x=\pm L$ be within a prescribed
tolerance of the actual endstates $U_\pm$ at plus and minus spatial infinity.
See \cite{Be1,Be2} for theoretical justification
and rigorous error bounds for this approach.
We follow this approach for the designer system \eqref{2bsys}.
However, for our gas-dynamical computations near the large-amplitude
limit $S_-\to -\infty$, this becomes impractical due to stiffness/presence
of multiple scales in the traveling-wave ODE, and we depart from this
method in favor of a simple shooting algorithm with a stiff ODE solver
(specifically, the adaptive mesh solver 
\textit{ode15s} supported in MATLAB).

We then approximate the 
Evans function $\tilde D$ numerically following the approach described in detail in \cite{BHZ1,BLeZ,BLZ} . To preserve analyticity of the Evans function, we use the method of Kato \cite{BHZ1, BrZ,GZ,HuZ1,K,Z4,Z7} to obtain a holomorphic initializing basis with respect to $\lambda$.  
To evolve in $x$ the manifolds whose determinant yields the Evans function, 
we use the polar coordinate (``analytic orthogonalization'') method described in \cite{HuZ1}. This method evolves an orthonormal solution basis along with a complex radial equation that maintains the property of analyticity. 
See \cite{HuZ1,Z7} for theoretical justification and error bounds.
To ensure an accurate winding number count, an adaptive mesh in $\lambda$ 
is used requiring that the relative error of change in $\tilde D$ for each step be less than 0.2. Recall, by Rouch\'e's Theorem, that if the relative variation of $D$ is less than 1.0, then winding number accuracy is preserved.
These routines are supported in STABLAB \cite{BHZ2}.


For gas-dynamical systems \eqref{ceg2} and \eqref{canon}, 
the Evans function is poorly conditioned due to multiple scales/stiffness 
of the eigenvalue ODE.
In this case, we find it necessary 
to use MATLAB's adaptive mesh stiff differential equation solver 
\textit{ode15s} both to solve for the profile and to perform 
the Evans function computations, 
rather than the Runge-Kutta-Fehlberg method encoded in \textit{ode45} 
that is typically used in Evans function computations.

\subsection{Hardware and computational statistics}\label{s:stats}

In this short section we detail computational statistics for some typical
values in the numerical studies reported in previous sections. All numerical
computations were performed on either a Mac Pro with 2 Quad-Core
Intel Xeon processors with speed 2.26 GHz, the super computer Quarry
at Indiana University, or a MacBook with an Intel Core 2 Duo processor
with speed 2.0 GHz. 

%
The Evans function computations were carried out using MATLAB's implicit ODE solver,
 ode15s, with relative error tolerance set at 1e-6 and absolute error set at 1e-8. 
For most of the parameters in our studies, we used MATLAB's boundary value solver bvp5c with relative
 error set at 1e-6, absolute error set at 1e-8, and the boundary error set at 1e-6. 
For the few parameters
 for which the boundary value solver failed, we 
used a shooting method as described above, with ode15s.
 
We now relate computational statistics, run on the MacBook, for each example system.
In each case, a minimum of 40 $\lambda$ points were used on a semi-circular mesh with 
radius $R = 10$. Additional
points were added to the mesh adaptively so that the relative change in consecutive 
contour points varied by no more than 0.2. The full contour was constructed by using the 
conjugate symmetry of the Evans function.
For $C = 10$ and $S_- = - 4$ in the global model, it took 23.6 seconds to compute the 
profile using a boundary value solver with a $\tanh$ solution as an initial guess. It took 176 seconds
compute the Evans function on 52 points. 
For $S_- = -4$ in the local model it took 62.1 seconds to compute the profile using a boundary value
solver and a $\tanh$ solution for an initial guess. It took 95.5 seconds to compute the Evans function on 63 points.
For $S_- = -4$ in the stable model it took 7.98 seconds to solve the profile with a boundary value solver initialized
with a $\tanh$ solution. It took 28.9 seconds to compute the Evans function on 40 points.

To give some perspective on the computational difficulty of the example models, we include data for the isentropic
gas model considered in \cite{BHRZ} with $\gamma = 5/3$ and $v_+ = 1e-4$. The boundary value solver took 0.899 seconds using a $\tanh$ solution for an initial guess. The Evans function took 27.7 seconds to compute on 40 points. As the above data indicates, solving the profile in the example gas systems considered here is much more difficult than in isentropic gas. The profile was amenable to solve using a boundary value solver for $S_-$ as small as $S_- = -200$ for the global model with $C = 10$, $S_- = -8.1$ for the local model, and $S_- = -11.5$ for the stable model. When using continuation to solve, we took steps of $0.1$ in $S_-$.
  
We also compare computational performance on the Mac Pro with its 8 cores, which is able to utilize parallel computing.
For $C = 10$ and $S_- = - 4$ in the global model, it took 16.9 seconds to compute the 
profile using a boundary value solver with a $\tanh$ solution as an initial guess. It took 28.75 seconds
compute the Evans function on 52 points. 
For $S_- = -4$ in the local model it took 37.78 seconds to compute the profile using a boundary value
solver and a $\tanh$ solution for an initial guess. It took 17.97 seconds to compute the Evans function on 63 points.
For $S_- = -4$ in the stable model it took 5.63 seconds to solve the profile with a boundary value solver initialized
with a $\tanh$ solution. It took 6.81 seconds to compute the Evans function on 40 points.
For the isentropic model with $\gamma = 5/3$ and $v_+ = 1e-4$, the boundary value 
solver took 0.5833 seconds using a $\tanh$ solution for an initial guess. 
The Evans function took  6.0929 seconds to compute on 40 points.

\end{document}